%% file: manuscript.tex
\pdfoutput=1 
\documentclass[
    a4paper,
    twocolumn
]{article}

\usepackage[english]{babel}
\usepackage[T1]{fontenc}
\usepackage[utf8]{inputenc}
\usepackage{eurosym}  

\usepackage{amstext} 
\DeclareRobustCommand{\officialeuro}{%
  \ifmmode\expandafter\text\fi
  {\fontencoding{U}\fontfamily{eurosym}\selectfont{}e}}

\usepackage{ifpdf}
\makeatletter
\ifpdf\ifdefined\pdffontattr
  \immediate\pdfobj stream {
    /CIDInit /ProcSet findresource begin
    12 dict begin
    begincmap
    /CIDSystemInfo
    << /Registry (TeX)
    /Ordering (Euro)
    /Supplement 0
    >> def
    /CMapName /TeX-Euro-0 def
    /CMapType 2 def
    1 begincodespacerange
    <00> <FF>
    endcodespacerange
    1 beginbfchar
    <65> <20AC>
    endbfchar
    endcmap
    CMapName currentdict /CMap defineresource pop
    end
    end
  }
  {\edef\@tempa#1#2{%
     \noexpand\fontseries{#1}\noexpand\fontshape{#2}\noexpand\selectfont
     \pdffontattr\font{/ToUnicode \the\pdflastobj\space 0 R}}
   \fontencoding{U}\fontfamily{eurosym}%
   \@tempa{m}{n}\@tempa{m}{sl}\@tempa{m}{ol}
   \@tempa{bx}{n}\@tempa{bx}{sl}
  }
\fi\fi
\makeatother

\usepackage{microtype}

\usepackage[
  a4paper,
  left=1.5cm,
  right=1.5cm,
  top=3cm,
  bottom=4cm
]{geometry}
\usepackage{fancyhdr}
\usepackage{lastpage}
\fancypagestyle{firststyle}{
   \fancyhead[R]{Page \thepage\ of \pageref*{LastPage}}
   \fancyfoot[L]{
    \footnotesize \rule{0.25\textwidth}{0.4pt} \newline
     $^{*}$\correspondingauthor
   }
   \fancyfoot[C]{}  
}

{
  \fancyhead[R]{Page \thepage\ of \pageref*{LastPage}}
}

\renewenvironment{abstract}{\noindent\textbf{Abstract:}}{}

\newcommand{\eg}{\text{e.g.}}
\newcommand{\ie}{\text{i.e.}}
\newcommand{\cf}{\text{cf.\;}}

\usepackage{amsfonts}
\usepackage{amsmath}
\usepackage{amssymb}
\usepackage{bm} 
\usepackage{mathtools}
\DeclarePairedDelimiter\card{\lvert}{\rvert}

\newcommand{\st}{\text{s.\,t.\;}}
\usepackage{xifthen}
\renewcommand{\min}[1][]{
	\ifthenelse{\isempty{#1}}{\operatorname{min}}{\ensuremath{\underset{#1}{\text{min}\,}}}
}
\newcommand{\der}[1]{\frac{\text{d}#1}{\text{d}t}} 

\usepackage{graphicx} 
\usepackage{float} 
\usepackage{scalerel} 
\usepackage{subcaption}
\usepackage[pdftex,dvipsnames]{xcolor} 
\definecolor{fzjblue}{RGB}{2,61,107} 
\colorlet{color1}{fzjblue}
\definecolor{fzjlightblue}{RGB}{173,189,227} 
\colorlet{color2}{fzjlightblue}
\definecolor{fzjgray}{RGB}{235,235,235} 
\colorlet{fzjgrey}{fzjgray}
\colorlet{color3}{fzjgray}
\definecolor{fzjred}{RGB}{235, 95, 115}  
\colorlet{color4}{fzjred}
\definecolor{fzjgreen}{RGB}{185, 210, 95}  
\colorlet{color5}{fzjgreen}
\definecolor{fzjyellow}{RGB}{250, 235, 90}  
\colorlet{color6}{fzjyellow}
\definecolor{fzjviolet}{RGB}{175, 130, 185}  
\colorlet{color7}{fzjviolet}
\definecolor{fzjorange}{RGB}{250, 180, 90}  
\colorlet{color8}{fzjorange}
\definecolor{fzjwhite}{RGB}{255,255,255}

\usepackage{pgfplots} 
\pgfplotsset{width=7cm, compat=newest}
\pgfkeys{/pgf/number format/.cd,
         1000 sep={\,},              
         min exponent for 1000 sep=4 
        }
\usepgfplotslibrary{groupplots}  
\usepgfplotslibrary{fillbetween} 
\usepgfplotslibrary{external}  

\newcommand{\setpgfexternalcounter}[1]{ 
  \makeatletter%
  \pgfkeysgetvalue{/tikz/external/figure name}\myexternalname
  \expandafter\gdef\csname c@tikzext@no@\myexternalname\endcsname{#1}%
  \makeatother
}

\usetikzlibrary{decorations.text} 
\usetikzlibrary{shapes} 
\usetikzlibrary{arrows.meta}  
\usepackage{jigsaw} 
\usepackage{adjustbox} 

\usepackage{enumitem} 
\usepackage{lscape, booktabs, tabularx, makecell} 
\usepackage{multicol} 
\newlist{mylist}{itemize}{1}
\setlist[mylist]{
  label=\textbullet,
  nosep, wide, leftmargin=*,
  before=\vspace{-0.57\baselineskip},
  after =\vspace{-0.85\baselineskip}
}

\usepackage{siunitx}
\usepackage{upquote}
\usepackage{etoolbox} 
\robustify{\texttt}
\let\originaltexttt\texttt
\begingroup
\catcode`'=\active
\catcode``=\active
\makeatletter
\renewrobustcmd{\texttt}[1]{%
   {%
   \everyeof{\noexpand}\endlinechar-1
   \expandafter\catcode\string``=\active
   \expandafter\catcode\string`'=\active
   \let'\textquotesingle
   \let`\textasciigrave
   \ifx\encodingdefault\upquote@OTone
    \ifx\ttdefault\upquote@cmtt
     \def'{\char13 }\def`{\char18 }%
    \fi
   \fi
   \scantokens{\originaltexttt{#1}}%
   }%
}%
\endgroup

\def\IEK10{
  Forschungszentrum Jülich GmbH,
  Institute of Energy and Climate Research,
  Energy Systems Engineering (IEK-10),
  Jülich 52425,
  Germany
}
\def\RWTH{
  RWTH Aachen University
  Aachen 52062,
  Germany
}
\def\ETH{
  ETH Zürich,
  Energy \& Process Systems Engineering,
  Zürich 8092,
  Switzerland
}
\def\JARA{
  JARA-ENERGY,
  Jülich 52425,
  Germany
}
\def\EBC{
  RWTH Aachen University,
  E.ON Energy Research Center,
  Institute for Energy Efficient Buildings and Indoor Climate,
  Aachen 52056,
  Germany
}
\def\LTT{
  RWTH Aachen University,
  Institute of Technical Thermodynamics,
  Aachen 52056,
  Germany
}
\def\SVT{
  RWTH Aachen University,
  Process Systems Engineering (AVT.SVT),
  Aachen 52074,
  Germany
}

\newcommand{\mytitle}{COMANDO: A Next-Generation Open-Source Framework for Energy Systems Optimization}

\newcommand{\affil}{
  \begin{itemize}[leftmargin=3mm, itemsep=-3mm]
    \item[$^a$] \IEK10 \\
    \item[$^b$] \RWTH \\
    \item[$^c$] \ETH \\
    \item[$^d$] \JARA \\
    \item[$^e$] \EBC \\
    \item[$^f$] \LTT \\
    \item[$^g$] \SVT
  \end{itemize}
}
\newcommand{\myauthor}{
  Marco Langiu$^{a, b}$,
  David Yang Shu$^{a, c}$,
  Florian Joseph Baader$^{a, b}$,
  Dominik Hering$^{a, b}$,
  Uwe Bau$^a$,
  André Xhonneux$^a$,
  Dirk Müller$^{d, a, e}$,
  André Bardow$^{d, a, f, c}$,
  Alexander Mitsos$^{d, a, g}$,
  Manuel Dahmen$^{a, *}$}

\newcommand{\correspondingauthor}{Manuel Dahmen, \IEK10 \newline E-mail: m.dahmen@fz-juelich.de}
\author{\myauthor}

\usepackage[
  colorlinks,
  linkcolor=blue,
  citecolor=blue,
  urlcolor=blue,
  bookmarks=true,
  pdftitle={\mytitle},
  pdfauthor={Marco Langiu}
]{hyperref}

\usepackage[capitalise, nameinlink]{cleveref}
\crefname{table}{Tab.}{Tab.}

\usepackage[round]{natbib}
\usepackage{doi}

\def\dv{\bm{x}}
\def\ov{\bm{y}_{s}(\cdot)}
\def\ovfunc{\bm{y}_{s}(t)}
\def\derov{\dot{\bm{y}}^\text{d}_s(t)}
\def\roc{\bm{f}}
\def\iv{\bm{y}^\text{d}_{s}(t=0)}
\def\design{I}
\def\operation{{II}}
\def\of{F}
\def\dobj{\of_\design}
\def\oobj{\of_{\operation,s}}
\def\moobj{\dot{\of}_\operation}
\def\myargs{\big(\dv, \ovfunc, \bm{p}_s(t)\big)}
\def\PGA{
  \begin{alignedat}{3}
    & \min[\dv \in \mathcal{X}]
      & \,
        & \dobj(\dv) + \textstyle{\sum\limits_{s\in \mathcal{S}}} w_s \, \oobj^*(\dv) \\
    & \underset{\hphantom{\dv, \bm{y}}}{\!\!\text{s.\,t.}}
      & & \oobj^*(\dv) = \hspace*{-4.5mm} \min[\ov \in \mathcal{Y}_s(\dv, \cdot)] \hspace*{-4mm} \oobj(\dv, \ov)\; \forall s \in \mathcal{S}
  \end{alignedat}
}
\def\bc{black!4}

\newcommand\acronym{\underline{c}omponent-\underline{o}riented \underline{m}od\-el\-ing and optimiz\underline{a}tion for \underline{n}onlinear \underline{d}esign and \underline{o}peration\,}

\begin{document}
  \twocolumn[
    \begin{@twocolumnfalse}
      \thispagestyle{firststyle}

      \begin{center}
        \begin{large}
          \textbf{\mytitle}
        \end{large} \\
        \myauthor
      \end{center}

      \vspace{0.5cm}

      \begin{footnotesize}
        \affil
      \end{footnotesize}

      \vspace{0.5cm}

      \begin{abstract}
        \noindent
        \begin{small}
          Existing open-source modeling frameworks dedicated to energy systems optimization typically utilize (mixed-integer) linear programming ((MI)LP) formulations, which lack modeling freedom for technical system design and operation.
          We present COMANDO, an open-source Python package for \acronym of integrated energy systems.
          COMANDO allows to assemble system models from component models including nonlinear, dynamic and discrete characteristics.
          Based on a single system model, different deterministic and stochastic problem formulations can be obtained by
          varying objective function and underlying data, and by applying automatic or manual reformulations.
          The flexible open-source implementation allows for the integration of customized routines required to solve challenging problems, e.g., initialization, problem decomposition, or sequential solution strategies.
          We demonstrate features of COMANDO via case studies, including automated linearization, dynamic optimization, stochastic programming, and the use of nonlinear artificial neural networks as surrogate models in a reduced-space formulation for deterministic global optimization.

          \vspace{0.5cm}


          \vspace{0.5cm}
        \end{small}
      \end{abstract}

      \begin{center}
        \input{figures/graphical_abstract}
      \end{center}

      \vspace{0.3cm}

      \begin{flushleft}
        \leavevmode{\parindent=17.5mm\indent}
        \textbf{Highlights:}
        {\begin{itemize}[leftmargin=22.5mm]
          \item Open-source framework for optimization of energy systems design and operation
          \item Component-oriented modeling, allowing for hybrid mechanistic/data-driven models
          \item Optimization considering nonlinearity, dynamics and parametric uncertainty
          \item Four case studies, demonstrating flexibility and wide range of application
        \end{itemize}}
      \end{flushleft}
      \vspace*{3mm}

      Published in \textit{Computers and Chemical Engineering}, doi:\href{https://doi.org/10.1016/j.compchemeng.2021.107366}{10.1016/j.compchemeng.2021.107366} \\
      \copyright{} 2021. This manuscript version is made available under the \href{https://creativecommons.org/licenses/by-nc-nd/4.0/}{CC-BY-NC-ND 4.0 license}
      \newpage

    \end{@twocolumnfalse}
  ]

  \input{sections/1_introduction.tex}
  \input{sections/2_M_AND_O.tex}
  \input{sections/3_COMANDO.tex}

  \input{sections/4_case_study.tex}
  \input{sections/5_conclusion.tex}
  \input{sections/6_acknowledgement.tex}
  \input{sections/7_appendices.tex}
  \bibliographystyle{abbrvnat}
  \bibliography{literature}

\end{document}

%% file: figures/graphical_abstract.tex
\begin{tikzpicture}
  \draw[thick]  (0, 0) rectangle (13, 5);
  \def\backgroundcolor{black!10}

  \node[anchor=north, align=center] at (4.25cm, 5cm) {\textbf{component-oriented} \\ \textbf{modular modeling}};
  \begin{scope}[
    shift={(2.5cm, 2.25cm)},
    scale=1.75, transform shape, very thick]
    \begin{scope}
      \color{black}
      \piece[fzjblue!10]{-1}{-1}{0}{0}
      \node at (0.7, 0.38) {\includegraphics[scale=0.04]{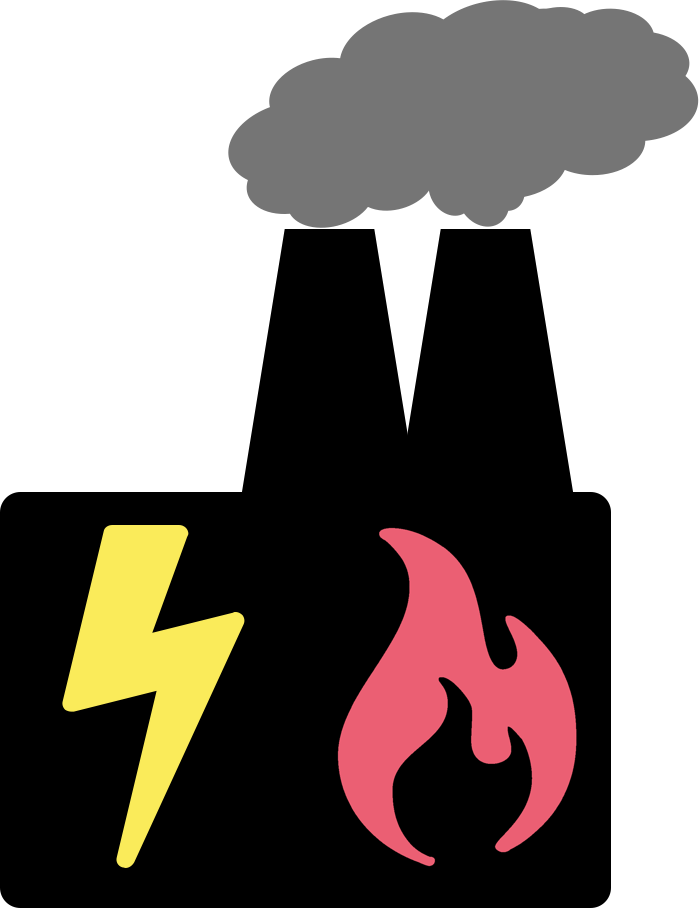}};
      \node at (0.4, 0.7) {\includegraphics[scale=0.025]{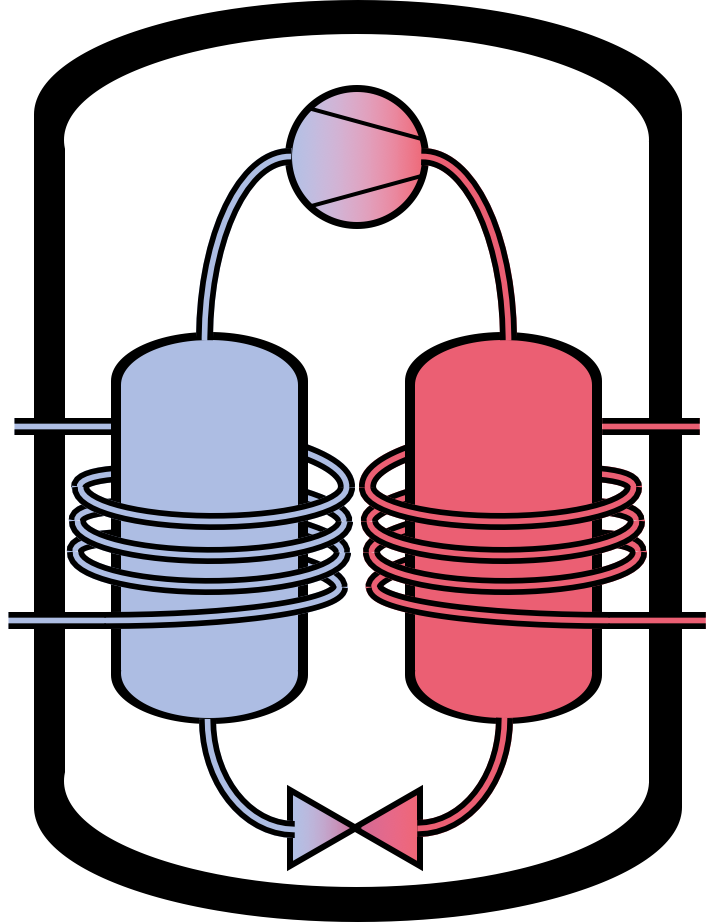}};
      \node at (0.25, 0.3) {\includegraphics[scale=0.03]{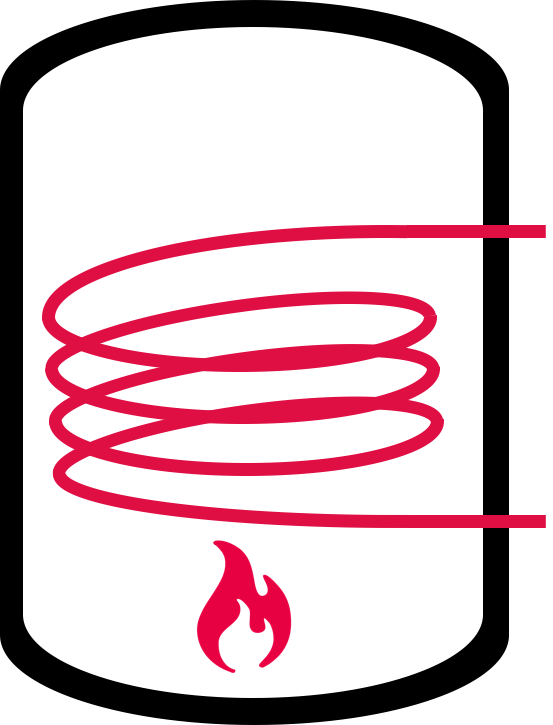}};
    \end{scope}
    \begin{scope}[xshift=1cm]
      \color{black}
      \piece[fzjblue!50]{1}{0}{0}{1}
      \node at (0.45, 0.7) {\includegraphics[scale=0.08]{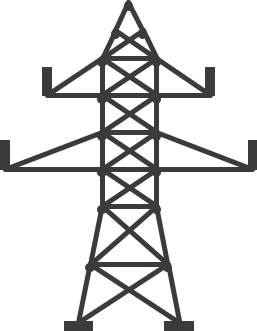}};
      \node at (0.8, 0.5) {\includegraphics[scale=0.08]{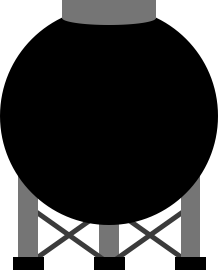}};
    \end{scope}
    \begin{scope}[yshift=-1cm]
      \color{black}
      \piece[fzjblue!30]{0}{1}{1}{0}
      \node at (0.25, 0.5) {\includegraphics[scale=0.08]{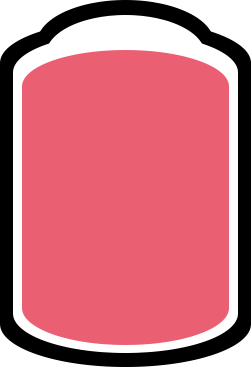}};
      \node at (0.58, 0.3) {\includegraphics[scale=0.12]{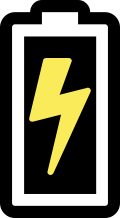}};
    \end{scope}
    \begin{scope}[xshift=1.15cm, yshift=-1.15cm]
      \color{black}
      \piece[fzjblue!10]{0}{0}{-1}{-1}
            \node at (0.44, 0.56) {\includegraphics[scale=0.06]{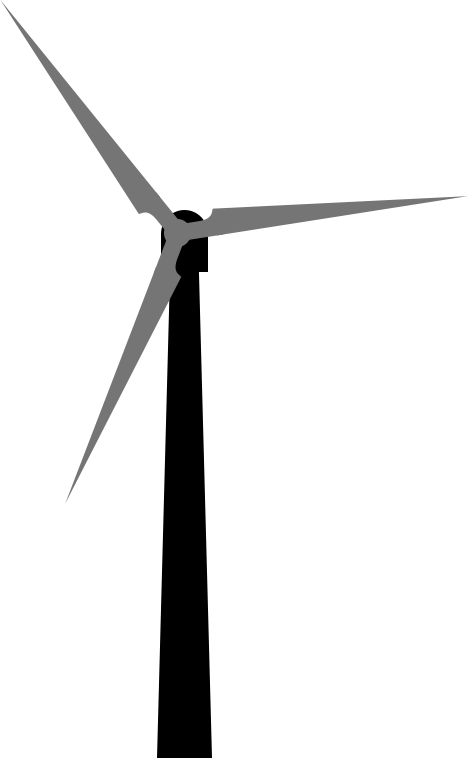}};
      \node at (0.7, 0.46) {\includegraphics[scale=0.04]{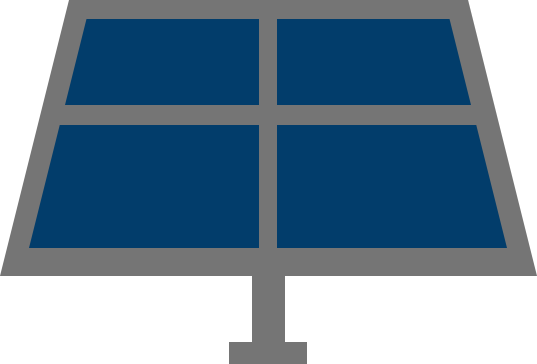}};
      \node at (0.5, 0.36) {\includegraphics[trim=122 0 0 0,clip, scale=0.08]{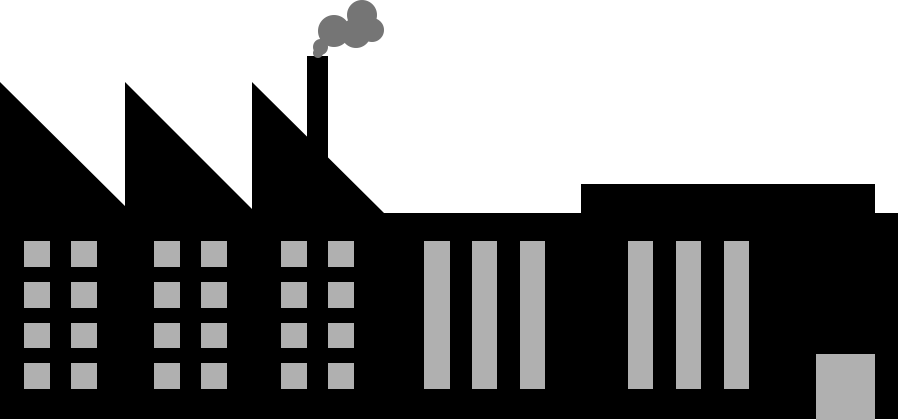}};
    \end{scope}
  \end{scope}

  \node[anchor=north, align=center] at (1.2cm, 4.8cm) {\textbf{nonlinearity}};
  \newcommand\expr[2]{#2 * (1.6 * #1 - #1^2)}
    \begin{axis}[
      shift={(0.5cm, 3cm)},
      scale=0.25,
      transform shape,
      clip=false,
      axis line style = {draw=none},
      line width=1.5pt,
      ticks=none,
      axis lines=middle,
      xmin=0, xmax=1,
      restrict x to domain = 0.15:1,
      ymin=0, ymax=0.8,
      colormap name=my_colormap,
      enlargelimits  = false,
      xlabel         = {part load},
      xlabel style   = {below left},
      ylabel         = {output},
      ylabel style   = {left=0.25cm, rotate=90},
      samples=200
    ]
      \addplot[]{\expr{x}{0.2}};
      \addplot[]{\expr{x}{0.4}};
      \addplot[]{\expr{x}{0.6}};
      \addplot[]{\expr{x}{0.8}};
      \addplot[]{\expr{x}{1}};
      \draw[->, >={Stealth[scale=0.5]}] (axis cs:0.8, 0.05) -- +(axis cs:-0.1, 0.75) node[left, rotate=30] {input};
      \coordinate (ymax) at (axis cs:0,0.8);
      \coordinate (origin) at (axis cs:0,0);
      \coordinate (xmax) at (axis cs:1,0);
   \end{axis}
   \draw[<->, >={Stealth[scale=0.75]}, ultra thick] (ymax) -- (origin) -- (xmax);

  \node[anchor=north, align=center] at (1.2cm, 2.5cm) {\textbf{dynamics}};
  \begin{scope}[shift={(0.5cm, 0.25cm)}, scale=0.4, >={Stealth[scale=0.75]}, ultra thick]

    \draw[shift={(0.5, 0.5)}, smooth,variable=\x, line join=round, fzjblue]
      (-0.5, 0) -- (0, 0) --
  plot[domain=0:1] ({\x}, {2 * (1 - exp(-\x/0.2))}) --
  plot[domain=1:1.75] ({\x}, {1.9865 - 2 * (1 - exp(-(\x - 1)/0.6))}) --
  plot[domain=1.75:2] ({\x}, {0.5595 + 2 * (1 - exp(-(\x - 1.75)/0.2))}) --
  plot[domain=2:4] ({\x}, {1.9865 - 2 * (1 - exp(-(\x - 2)/0.6))});
  \draw[very thick, dotted, fzjred, shift={(0.5, -0.5)}] (2.25, 2) -- node[right] {$\substack{\dot{y}(t)\,=\\f(y, t)}$} +(-0.5, 2.6);
    \draw[<->] (4.5, 0) node[right] {$t$} --
               (0, 0) -- (0, 3.8) node[shift={(-0.25cm, -0.6cm)}, rotate=90] {$y(t)$};
  \end{scope}

  \node[align=center] at (7.25cm, 2.3cm) {\textbf{uncertain} \\ \textbf{data}};
    \def\demand{
        -0.3*cos(x*360)
    }
  \newcommand\randomweather[4]{      \pgfmathsetseed{#1}%
    \addplot[draw=none, no markers, smooth, name path=TD] {#2 + \demand + #3 * (1.5 + 0.5 * rand)};
    \pgfmathsetseed{#1}%
    \addplot[draw=none, no markers, smooth, name path=BD] {#2 + \demand - #3 * (1.5 - 0.5 * rand)};
    \pgfmathsetseed{#1}%
    \addplot[no markers, smooth,  #4!30] fill between[of=TD and BD];
    \addplot[very thick, no markers, smooth, #4] {#2 + \demand + #3 * rand};
  }
    \newcommand\randompowerdemand[4]{
      \pgfmathsetseed{#1}
      \addplot[draw=none, no markers, smooth, name path=T] coordinates {
        (0.010, #2 + #3 * 0.3 + #3 * 0.1 * rand + #3 * 0.433)
        (0.054, #2 + #3 * 0.3 + #3 * 0.1 * rand + #3 * 0.302)
        (0.128, #2 + #3 * 0.3 + #3 * 0.1 * rand + #3 * 0.261)
        (0.208, #2 + #3 * 0.3 + #3 * 0.1 * rand + #3 * 0.260)
        (0.273, #2 + #3 * 0.3 + #3 * 0.1 * rand + #3 * 0.360)
        (0.312, #2 + #3 * 0.3 + #3 * 0.1 * rand + #3 * 0.515)
        (0.382, #2 + #3 * 0.3 + #3 * 0.1 * rand + #3 * 0.595)
        (0.462, #2 + #3 * 0.3 + #3 * 0.1 * rand + #3 * 0.566)
        (0.542, #2 + #3 * 0.3 + #3 * 0.1 * rand + #3 * 0.588)
        (0.621, #2 + #3 * 0.3 + #3 * 0.1 * rand + #3 * 0.541)
        (0.689, #2 + #3 * 0.3 + #3 * 0.1 * rand + #3 * 0.646)
        (0.731, #2 + #3 * 0.3 + #3 * 0.1 * rand + #3 * 0.791)
        (0.804, #2 + #3 * 0.3 + #3 * 0.1 * rand + #3 * 0.864)
        (0.870, #2 + #3 * 0.3 + #3 * 0.1 * rand + #3 * 0.816)
        (0.926, #2 + #3 * 0.3 + #3 * 0.1 * rand + #3 * 0.715)
        (0.970, #2 + #3 * 0.3 + #3 * 0.1 * rand + #3 * 0.573)
        (1.002, #2 + #3 * 0.3 + #3 * 0.1 * rand + #3 * 0.466)
      };
      \addplot[draw=none, no markers, smooth, name path=B] coordinates {
        (0.010, #2 - #3 * 0.25 + #3 * 0.433)
        (0.054, #2 - #3 * 0.25 + #3 * 0.302)
        (0.128, #2 - #3 * 0.25 + #3 * 0.261)
        (0.208, #2 - #3 * 0.25 + #3 * 0.260)
        (0.273, #2 - #3 * 0.25 + #3 * 0.360)
        (0.312, #2 - #3 * 0.25 + #3 * 0.515)
        (0.382, #2 - #3 * 0.25 + #3 * 0.595)
        (0.462, #2 - #3 * 0.25 + #3 * 0.566)
        (0.542, #2 - #3 * 0.25 + #3 * 0.588)
        (0.621, #2 - #3 * 0.25 + #3 * 0.541)
        (0.689, #2 - #3 * 0.25 + #3 * 0.646)
        (0.731, #2 - #3 * 0.25 + #3 * 0.791)
        (0.804, #2 - #3 * 0.25 + #3 * 0.864)
        (0.870, #2 - #3 * 0.25 + #3 * 0.816)
        (0.926, #2 - #3 * 0.25 + #3 * 0.715)
        (0.970, #2 - #3 * 0.25 + #3 * 0.573)
        (1.002, #2 - #3 * 0.25 + #3 * 0.466)
      };
      \addplot[no markers, smooth,  #4!30] fill between[of=T and B];
      \addplot[very thick, no markers, smooth, #4] coordinates {
        (0.010, #2 + #3 * 0.433)
        (0.021, #2 + #3 * 0.411)
        (0.029, #2 + #3 * 0.356)
        (0.026, #2 + #3 * 0.391)
        (0.041, #2 + #3 * 0.324)
        (0.054, #2 + #3 * 0.302)
        (0.068, #2 + #3 * 0.293)
        (0.082, #2 + #3 * 0.280)
        (0.096, #2 + #3 * 0.273)
        (0.112, #2 + #3 * 0.268)
        (0.128, #2 + #3 * 0.261)
        (0.144, #2 + #3 * 0.260)
        (0.160, #2 + #3 * 0.267)
        (0.176, #2 + #3 * 0.287)
        (0.192, #2 + #3 * 0.274)
        (0.208, #2 + #3 * 0.260)
        (0.224, #2 + #3 * 0.263)
        (0.240, #2 + #3 * 0.281)
        (0.256, #2 + #3 * 0.307)
        (0.265, #2 + #3 * 0.333)
        (0.273, #2 + #3 * 0.360)
        (0.280, #2 + #3 * 0.393)
        (0.288, #2 + #3 * 0.422)
        (0.295, #2 + #3 * 0.452)
        (0.302, #2 + #3 * 0.486)
        (0.312, #2 + #3 * 0.515)
        (0.321, #2 + #3 * 0.544)
        (0.335, #2 + #3 * 0.572)
        (0.351, #2 + #3 * 0.589)
        (0.367, #2 + #3 * 0.595)
        (0.382, #2 + #3 * 0.595)
        (0.398, #2 + #3 * 0.591)
        (0.414, #2 + #3 * 0.584)
        (0.430, #2 + #3 * 0.576)
        (0.446, #2 + #3 * 0.566)
        (0.462, #2 + #3 * 0.566)
        (0.478, #2 + #3 * 0.575)
        (0.494, #2 + #3 * 0.575)
        (0.510, #2 + #3 * 0.587)
        (0.526, #2 + #3 * 0.595)
        (0.542, #2 + #3 * 0.588)
        (0.558, #2 + #3 * 0.576)
        (0.573, #2 + #3 * 0.561)
        (0.589, #2 + #3 * 0.548)
        (0.605, #2 + #3 * 0.542)
        (0.621, #2 + #3 * 0.541)
        (0.638, #2 + #3 * 0.548)
        (0.652, #2 + #3 * 0.561)
        (0.667, #2 + #3 * 0.583)
        (0.680, #2 + #3 * 0.616)
        (0.689, #2 + #3 * 0.646)
        (0.696, #2 + #3 * 0.674)
        (0.703, #2 + #3 * 0.704)
        (0.711, #2 + #3 * 0.734)
        (0.721, #2 + #3 * 0.763)
        (0.731, #2 + #3 * 0.791)
        (0.742, #2 + #3 * 0.817)
        (0.756, #2 + #3 * 0.841)
        (0.772, #2 + #3 * 0.858)
        (0.788, #2 + #3 * 0.863)
        (0.804, #2 + #3 * 0.864)
        (0.819, #2 + #3 * 0.862)
        (0.834, #2 + #3 * 0.855)
        (0.846, #2 + #3 * 0.843)
        (0.860, #2 + #3 * 0.835)
        (0.870, #2 + #3 * 0.816)
        (0.884, #2 + #3 * 0.804)
        (0.894, #2 + #3 * 0.781)
        (0.906, #2 + #3 * 0.763)
        (0.915, #2 + #3 * 0.740)
        (0.926, #2 + #3 * 0.715)
        (0.938, #2 + #3 * 0.687)
        (0.948, #2 + #3 * 0.658)
        (0.958, #2 + #3 * 0.629)
        (0.966, #2 + #3 * 0.601)
        (0.970, #2 + #3 * 0.573)
        (0.978, #2 + #3 * 0.551)
        (0.988, #2 + #3 * 0.522)
        (0.997, #2 + #3 * 0.490)
        (1.002, #2 + #3 * 0.466)
      };
    }
  \newcommand\peak[2]{        #2 * exp(-10000*(x-#1)^2)
    }
    \newcommand\randompowerprice[4]{
      \pgfmathsetseed{#1}
      \addplot[draw=none, no markers, smooth, name path=TP] {#2 + #3 * 0.1 * rand + .7};
      \addplot[draw=none, no markers, smooth, name path=BP, samples=2] {0};
      \addplot[no markers, smooth,  #4!30] fill between[of=TP and BP];
      \addplot[very thick, no markers, smooth, samples=150, #4] {#2 + #3 * 0.1 * rand + \peak{0.2}{0.4} - \peak{0.3}{0.6 * #2} - \peak{0.5}{1 * #2} + \peak{0.85}{0.4}};
    }
  \begin{axis}[
               shift={(7.25cm, .25cm)},
               scale=0.25,
               transform shape,
               clip=false,
               axis on top = true,
               line width=1.5pt,
               ticks= none,
               domain=0:1,
               width=7cm,
               ymin=0, ymax=4,
               enlargelimits=false,
               axis lines*=middle,
               axis line style={draw=none},
               samples=30,
      ]        \randompowerprice{2}{3}{1.8}{fzjblue};
      \randomweather{2}{2}{0.3}{fzjred};
      \randompowerdemand{0}{-0.1}{1.8}{fzjyellow};
      \coordinate (ymax) at (axis cs:0,4);
      \coordinate (origin) at (axis cs:0,0);
      \coordinate (xmax) at (axis cs:1,0);
   \end{axis}
   \draw[<->, >={Stealth[scale=0.75]}, ultra thick] (ymax) -- node[left, xshift=-5mm, yshift=6mm, rotate=90] {$\substack{
                      \text{prices}\\
                      \text{weather}\\
                      \text{demands}
                     }$} (origin) -- (xmax) node[xshift=2mm] {$t$};

  \begin{scope}[shift={(8.5, 1.75)}, scale=0.3,
                every node/.style={circle, draw, thick, minimum size=1cm, inner sep=1pt, outer sep=0pt, transform shape}
  ]
    \node {\Large $\dv$} [grow'=right, level distance=2cm, sibling distance=1.3cm]
    child {node {\Large $\bm{y}_1(t)$}
         edge from parent node [draw=none, above] {\Large $s_1$}
        }
    child {node {\Large $\bm{y}_2(t)$}
         edge from parent node [draw=none, above, xshift=1.5mm, yshift=-1.5mm] {\Large $s_2$}
        }
    child {node[draw=none, label={[label distance=-7.5mm]90:{\huge $\vdots$}}] {} edge from parent[draw=none]}
    child {node {\Large $\bm{y}_{\card{\mathcal{S}}}(t)$}
         edge from parent node [draw=none, right, yshift=0.5mm] {\Large $s_{\card{\mathcal{S}}}$}
        };
  \end{scope}

  \coordinate (start) at (8.75cm, 0.6cm);
  \coordinate (end) at (9.5cm, 2.6cm);

  \def\outgoing{-15}
  \def\incoming{300}
  \draw[-{Stealth[angle=120:1pt 1]}, fzjblue!50, line width=2ex] (start) to[out=\outgoing,in=\incoming] (end);

  \path [postaction={
            decoration={
                text align={left indent=0.5cm},
                text along path,
                text={|\bfseries|scenarios},
                raise=1.5pt
            },
            decorate
         }
        ]
        [out=\outgoing] (start) ++(0.25cm, 0) to [in=
        \incoming+30] (end);

    \coordinate (algo) at (11.25cm, 1cm);

     \draw[-{Stealth[angle=120:1pt 1]}, fzjblue!50, line width=2ex] (end) ++(0.75cm, 0) to[out=\incoming-45,in=195] (algo);

      \path [postaction={
            decoration={
                text align={left indent=0.5cm},
                text along path,
                text={|\bfseries|solution},
                raise=-3.5pt
            },
            decorate
         }
        ]
        [out=\incoming-45] (end) ++(0.75cm, 0) to [in=195] (algo);

    \node[anchor=north] at (9.75, 4.25) (OP) {
        \!\!$\PGA$
    };
    \node[anchor=north, align=center] at (9.5cm, 5cm) {\textbf{optimization of design \& operation} \\ \textbf{under uncertainty}
  };

  \node[align=center] at (11.75cm, 2.3cm) {\textbf{algorithm} \\ \textbf{development}};
  \begin{scope} [shift={(12, 0.7)}, scale=0.15, >={Stealth[black, scale=0.4]}, thick, transform shape]
    \node[circle, draw, minimum size=1.5cm] (c1) at (0, 7) {};
    \node[circle, draw, minimum size=1cm] (cc1) at (0, 7) {};
    \node (d1) at (0, 5) [draw, diamond, aspect=2, scale=4, fill=fzjorange] {};

    \node[rectangle, draw, fill=fzjblue!50, minimum width=2cm, minimum height=1cm] (r1) at (0, 3) {};
    \node[rectangle, draw, fill=fzjblue!50, minimum width=2cm, minimum height=1cm] (r2) at (-3, 3) {};
    \node[rectangle, draw, fill=fzjblue!50, minimum width=2cm, minimum height=1cm] (r3) at (0, 1) {};
    \node[rectangle, draw, fill=fzjblue!50, minimum width=2cm, minimum height=1cm] (r4) at (3, 1) {};

    \node (d2) at (-3, 1) [draw, diamond, aspect=2, scale=4, fill=fzjorange] {};

    \node (d3) at (0, -1) [draw, diamond, aspect=2, scale=4, fill=fzjorange] {};

     \node[circle, draw, minimum size=1.5cm] (c2) at (0, -3) {};
     \node[circle, draw, fill, minimum size=1cm] (cc2) at (0, -3) {};

     \draw[->] (c1) -- (d1);
     \draw[->] (d1) -- (r1);
     \draw[->] (r1) -- (r3);
     \draw[->] (r3) -- (d3);
     \draw[->] (d3) -- (c2);

     \draw[->] (d1) -| (r2);
     \draw[->] (r2) -- (d2);
     \draw[->] (d2) |- (c2);

     \draw[->] (d2) -- (r3);

     \draw[->] (d3) -| (r4);
     \draw[->] (r4) |- (r1);
  \end{scope}
\end{tikzpicture}

%% file: sections/1_introduction.tex

\section{Introduction}
\label{sec:introduction}

Energy systems are networks of interconnected components that generate and transform energy using a set of renewable or fossil resources to satisfy various kinds of demands \citep{beller1976reference}.
The economic and ecologic performance of energy systems is strongly influenced by system design and operation.
The design comprises all choices regarding the configuration, \ie, the selection and interconnection of components (discrete variables), as well as sizing and other technical specifications (continuous variables).
The operation comprises commitment (discrete variables) and dispatch (continuous variables) of individual components, i.e., how their activity and output levels are chosen at different points in time.
The prospective operation also needs to be taken into account during system design \citep{pistikopoulos1995uncertainty, frangopoulos2002brief}.
However, energy demands, prices, weather and other operational aspects can be highly variable and their future values are inherently uncertain, rendering
the design and operation of energy systems a challenging decision process.
To ensure optimal economic and ecologic performance, it is common to cast these decision processes into mathematical optimization problems \citep[\eg][]{papoulias1983structural, ghobeity2012optimal, gunasekaran2014optimal, andiappan2017state, frangopoulos2018recent, demirhan2019energy, sass2020model}.
This is typically done via general purpose algebraic modeling languages (AMLs), \eg, GAMS \citep{bussieck2004general} or Pyomo \citep{hart2011pyomo}, or via specialized energy system modeling frameworks (ESMFs), \eg, OSeMOSYS \citep{howells2011osemosys} or oemof \citep{hilpert2018open}.
While AMLs offer flexibility in the choice of algebraic formulation and solution approach, %
ESMFs employ a component-oriented modeling approach, \ie, system models are created by specifying connections between component models.
This approach simplifies the modeling process, model maintenance, and model re-use.

Established ESMFs typically employ linear programming (LP) \citep{schrattenholzer1981energy, fishbone1981markal, loulou2007etsap, bakken2007etransport, howells2011osemosys, hunter2013modeling, dorfner2016open} or mixed-integer linear programming (MILP) \citep{pfenninger2015renewables, hilpert2018open, atabay2017open, brown2018pypsa, johnston2019switch} formulations, well-suited for techno-economic analysis of large-scale systems \citep{connolly2010review, pfenninger2014energy, beuzekom2015review}.
In contrast, technical system design and operation must consider more detailed system behavior, often giving rise to nonlinearities and dynamic effects that are difficult or impractical to represent with MILP formulations, \citep[see \eg,][]{li2011stochastic, goderbauer2016adaptive, schaefer2019reduced, schaefer2019economic}.
To address the challenges of technical design and operation, we propose a next-genera\-tion ESMF for \acronym (COMANDO), an open source Python package \citep{COMANDO_REPO}.
COMANDO borrows a generic, nonlinear representation of mathematical expressions and features for algorithm development from AMLs, and the representation of differential equations and more general system model aggregation from differential-algebraic modeling frameworks (DAMFs) such as gPROMS \citep{gPROMS}, MODELICA \citep{elmqvist1997modelica}, or DAE Tools \citep{nikolic2016dae}.
With this combination of features, COMANDO incorporates flexible nonlinear and dynamic modeling into the modularity of an ESMF.
Additionally, COMANDO enables the simultaneous consideration of multiple operating scenarios through a two-stage stochastic programming formulation, allowing for rigorous optimization of energy system design and operation under uncertainty and/or variability of operating conditions.
While the vast majority of existing ESMFs is implemented as a layer on top of an AML, COMANDO is based on the computer algebra system SymPy \citep{meurer2017sympy}.
SymPy provides data structures for representing generic mathematical expressions
and corresponding methods to analyze and manipulate expressions.
These features facilitate the creation of automatic reformulation routines (\eg, automatic linearization), custom interfaces to AMLs or solvers, and user-defined solution algorithms.

This paper is structured as follows:
In \cref{sec:M_AND_O}, we give a brief review of the state of the art in optimization-based energy-system design and operation and identify the lack of an open-source tool dedicated specifically to the technical design and operation of different types of energy systems.
To this end, we present COMANDO in \cref{sec:COMANDO}.
In \cref{sec:case_study}, we present four case studies highlighting important features of COMANDO.
\cref{sec:conclusion} concludes the work.

%% file: sections/2_M_AND_O.tex

\section{Optimization-based energy system design and operation}
\label{sec:M_AND_O}

In \cref{sec:problem} we introduce a generic mathematical programming problem for the optimal design and operation of energy systems.
In \cref{sec:tools}, we briefly summarize advantages and disadvantages of the three major classes of tools that can be used to formulate and tackle variants of this problem, namely algebraic modeling languages (AMLs), energy system modeling frameworks (ESMFs) and differential-algebraic modeling frameworks (DAMFs).

\subsection{Problem formulation} \label{sec:problem}
  Realizing an optimal energy system requires optimal decisions at both the design stage and the operational stage.
  Due to the variability and uncertainty associated to energy system operation, there can be many relevant operational scenarios that need to be considered to obtain a reliable design.
  A suitable modeling approach for this setting is two-stage stochastic programming \citep{dantzig1955linear, birge2011introduction, li2015optimal, yunt2008designing}.
  It allows for the simultaneous consideration of multiple operating scenarios $s \in \mathcal{S}$, resulting in the following problem structure:
  \makeatletter
  \newcases{subproblem}{\quad}{%
    \hfil$\m@th\displaystyle{##}$}{$\m@th\displaystyle{##}$\hfil}{.}{\rbrace}
  \makeatother
  \begin{equation}
  \label{prob:P} \tag{P}
    \scriptsize
    \begin{alignedat}{3}
      & \min[\dv]
        && \dobj(\dv) + \sum_{s \in \mathcal{S}} w_s \, \oobj^*(\dv) \\[-1mm]
      & \underset{\hphantom{\dv, \bm{y}}}{\text{s.\,t.}}
        && \bm{g}_\design(\dv) \leq \bm{0}\\[-2mm]
      & && \bm{h}_\design(\dv) = \bm{0}\\[-3mm]
      & && \!\!\!\begin{subproblem}
              \oobj^*(\dv) = \min[\bm{y}_s(\cdot)] \!\!\!\!
                & \oobj(\dv, \ov) = \!\!
                \int_{\mathcal{T}_s} \!\!\!\!\moobj\myargs \; \text{d}t \\[-2mm]
              \!\st \!\!\!\!
                & \iv = \bm{y}^\text{d}_{s, 0} \\
                & \begin{aligned}[t]\\[-\dimexpr\baselineskip+\fontdimen22\textfont2\relax]
                  \begin{rcases}
                    \derov = \roc\myargs \!\!\!&\\
                    \bm{g}_\operation\myargs \leq \bm{0} &\\
                    \bm{h}_\operation\myargs = \bm{0} &\\[-1mm]
                    \bm{y}_{s}(t) = [\bm{y}^\text{d}_{s}(t), ...] \\
                    \bm{y}_{s}(t) \in \mathcal{Y}_{s}(t) \subset \mathbb{R}^{n_y} \!\times \mathbb{Z}^{m_y} \!\!\!\!\!\!\!\!\! \\
              \end{rcases} \; \forall t \in \mathcal{T}_s
              \end{aligned} \\[-1mm]
              & \mathcal{T}_s = \left[0, T_s\right] \\[-2mm]
          \end{subproblem} \; \forall s \in \mathcal{S} \\[1mm]
      & && \dv \in \mathcal{X} \subset \mathbb{R}^{n_x} \!\times \mathbb{Z}^{m_x} \\[-1mm]
      & && \mathcal{S} = \{s_1, s_2, \cdots, s_{\card{\mathcal{S}}}\}
    \end{alignedat}
  \end{equation}
  %
  The two-stage structure of \eqref{prob:P} distinguishes between design- and operation-related variables, constraints, and objectives.
  We group design decisions into the vector $\dv$ and operational decisions into one vector $\ov$ for each scenario $s$, with associated probability of occurrence $w_s$.
  Further, the operational decisions are functions of time $t$ from a continuous operating horizon $\mathcal{T}_s = [0, T_s]$
  (in general, each scenario may consider a different time horizon).
  Likewise, for different scenarios $s$ the input data, i.e., the values of model parameters $\bm{p}_s(\cdot)$, may be functions of time $t$.
  The objective function of the first stage is comprised of design costs $\dobj$ and the expected value of the optimal operating costs.
  For a given design $\dv$ and scenario $s$, the optimal operating costs $\oobj^*$ correspond to the optimal objective value of the second stage.
  The operating costs are described by an integral over the operating horizon $\mathcal{T}_s$ of the momentary operating costs $\moobj$.
  The set of feasible design and operational decisions is described via constraints $\bm{g}_{\design}$, $\bm{g}_{\operation}$, $\bm{h}_{\design}$, and $\bm{h}_{\operation}$ (with an appropriate number of elements in $\bm{h}_{\operation}$, allowing for degrees of freedom), as well as bounds and integrality restrictions in the form of $\mathcal{X}$ and $\mathcal{Y}_{s}(t)$, with $n$ and $m$ corresponding to the number of continuous and discrete decisions, respectively.
  Additionally, for the subset of operational variables that correspond to differential states (identified via the superscript $^{\text{d}}$), an initial state $\bm{y}^\text{d}_{s, 0}$ and the right hand side $\roc$ of a corresponding differential equation are given.

  Formulation \eqref{prob:P} covers both mixed design and operation problems, as well as pure operational problems (with fixed design decisions $\dv$).
  If the values of $w_s$ are interpreted as frequencies of occurrence for a certain operational setting, the corresponding scenarios can also be interpreted as typical operating points or periods, as done \eg, in \citet{yunt2008designing} and \citet{baumgaertner2019rises3}, respectively.
  Such scenarios can be derived from standardized reference load profiles, or via clustering of historical data \citep[see, \eg,][]{schuetz2018comparison}.
  If constraints coupling different scenarios are added to formulation \eqref{prob:P}, problems considering long-term effects such as seasonal storage can also be considered \citep[see, \eg,][]{gabrielli2018optimal, baumgaertner2019rises4}.

  \def\Tdisc{\widehat{\mathcal{T}}}
  \def\Tdiscsub{\Tdisc}
  The two-stage formulation \eqref{prob:P} can be cast into an equivalent single-stage formulation, also referred to as the \emph{deterministic equivalent}, see \eg, \citep{yunt2008designing}, that can be solved with general-purpose solvers.
  While solvers interfaced from DAMFs, as well as some specialized dynamic optimization solvers, \eg, DyOS \citep{caspari2019dyos}, directly accept continuous-time problem formulations and take care of time-discretization internally, almost all solvers available via AMLs and ESMFs require discrete-time formulations as input.
  To obtain a discrete-time formulation, a particular discretization scheme is chosen, and $\ov$, $\oobj(\dv, \ov)$, and $\roc\myargs$ are replaced by corresponding discrete-time counterparts.

  An alternative to solving the deterministic equivalent is to employ an algorithm capable of exploiting the special constraint structure of the two-stage formulation \eqref{prob:P}.
  Such an algorithm decomposes \eqref{prob:P}
  into multiple subproblems that are solved iteratively to obtain increasingly tighter bounds on the solution of \eqref{prob:P}.
  Different decomposition algorithms are applicable, depending on the presence and location of nonlinearity, nonconvexity and integrality; for a concise overview, see \citet{li2019generalized}.

\subsection{Tools} \label{sec:tools}
  Both deterministic equivalent formulations as well as suitable decomposition algorithms can be implemented in AMLs such as AMPL \citep{fourer1990modeling}, GAMS \citep{bussieck2004general}, or AIMMS \citep{bisschop2006aimms}.
  In recent years, several AML extensions have been developed that can be leveraged for
  energy system
  modeling.
  In particular, stochastic programming related functionality has been incorporated widely, both in commercial AMLs such as AMPL \citep{fourer1990modeling} and GAMS \citep{bussieck2004general} (through SAMPL \citep{valente2009extending} and Extended Mathematical Programming \citep{ferris2009extended}, respectively), as well as in the open-source AMLs Pyomo \citep{hart2011pyomo} and JuMP \citep{dunning2017jump} (through PySP \citep{watson2012pysp} and Struct\-JuMP, formerly StochJuMP, \citep{huchette2014parallel} or StochasticPrograms.jl \citep{biel2019efficient}, respectively).
  Further modeling constructs tailored towards special problem structures have been incorporated through block-oriented modeling \citep{friedman2013block} in Pyomo and through Plasmo.jl \citep{jalving2017graph, jalving2019graph} in JuMP.
  Finally, Pyomo.DAE \citep{nicholson2018pyomo.dae} enables the direct representation of differential equations within optimization problems expressed in Pyomo and provides various options for automatic discretization.
  Through the combination of features offered by these extensions, newer AMLs are in principle well suited to model and optimize energy system design and operation.
  However, their abstract nature can complicate implementation, code maintenance, and re-use, and renders the resulting problem formulations difficult to comprehend.
  Development on the Pyomo AML has resulted in the modeling tool IDAES \citep{miller2018next}, which employs methodologies from process systems
  engineering with the aim of advancing fossil energy systems \citep{IDAEShomepage}.
  In particular, IDAES provides models for thermal power plants and associated components.
  These systems are considered in the form of process flowsheets, \ie, components are modeled as control volumes with in- and outflows, whose steady-state and dynamic behavior can be specified via so-called property packages.

  Compared to AMLs and their various extensions,
  ESMFs provide an even higher level of abstraction, allowing to model generic energy systems comprised of utilities for generating, converting, or storing different energy forms.
  This higher level of abstraction is commonly achieved via an interface layer on top of an AML that separates component and system modeling from problem formulation.
  In a first modeling step, models of energy system components, \eg, boilers, combined-heat-and-power units, or heat pumps are created.
  These component models contain variables, parameters and constraints specifying possible in- and outputs as well as the internal component behavior.
  In a second modeling step, system models are aggregated by specifying the connections between different components.
  Finally, component- and system-level constraints are combined with an objective, \eg, the minimization of total annualized cost (TAC) or global warming impact (GWI), yielding a problem formulation that can be passed to an appropriate solver.

  The modular, object-oriented nature of modern ESMFs such as oemof \citep{hilpert2018open} allows component and system models to be implemented as classes, inheriting reoccurring functionality, \eg, from generic models representing generation, transformation, storage or consumption of different energy commodities.
  Such inheritance allows for more structured modeling, thereby simplifying model maintenance and re-use compared to AMLs, \eg, through the creation of component libraries.
  However, the vast majority of ESMFs is based on either linear programming (LP) or mixed-integer linear programming (MILP) formulations, \ie, all participating functions must be linear in the decision variables $\bm{x}$ and $\bm{y}$.
  In such ESMFs, the user must provide linear approximations for all nonlinear expressions.
  While this is usually not considered a limitation in the context of \emph{system analysis}, \ie, the principal focus of most ESMFs \citep[\cf][]{pfenninger2014energy}, problems concerned with \emph{technical design and operation} need to represent systems in more detail, often giving rise to nonlinearities that are difficult or impractical to linearize.
  In the presence of such nonlinearities, it is
  often sensible to use the original nonlinear equations or nonlinear surrogate models such as artificial neural networks (ANNs), as, \eg, in \citet{schaefer2020wavelet}, which however is not possible in MILP-based ESMFs.

  \begin{table*}[htb!]
    \centering
    \caption{Overview of the three tool classes that inspired COMANDO: algebraic modeling languages (AMLs), energy system modeling frameworks (ESMFs), and differential-algebraic modeling frameworks (DAMFs)}
    \scriptsize
    \label{tab:tool_classes}
    \begin{tabularx}{1\linewidth}{
        %
        >{\hsize=0.2\hsize\linewidth=\hsize}X
        >{\hsize=1.55\hsize\linewidth=\hsize}X
        >{\hsize=1.1\hsize\linewidth=\hsize}X
        >{\hsize=1.15\hsize\linewidth=\hsize}X
    }
      \toprule
      \thead{tool \\ class}
        & \thead{representative \\ examples}
          & \thead{typical domain \\ of application}
            & \thead{features adopted \\ in COMANDO} \\
      \midrule
      AMLs
        & \begin{mylist}
            \item AMPL 
                  \citep{fourer1990modeling}
            \item GAMS 
                  \citep{bussieck2004general}
            \item AIMMS 
                  \citep{bisschop2006aimms}
            \item PYOMO 
                  \citep{hart2011pyomo}
            \item JuMP 
                  \citep{dunning2017jump}
          \end{mylist}
        & development of detailed, app\-li\-ca\-tion-spe\-ci\-fic (MI)LP/(MI)NLP problem formulations for arbitrary applications and specialized solution routines
        & \begin{mylist}
            \item free choice of modeling approach
            \item possibility to specify alternative problem formulations
            \item development of user-defined algorithms
          \end{mylist} \\
      \midrule
      ESMFs
        & \begin{mylist}
            \item MESSAGE 
                  \citep{schrattenholzer1981energy}
            \item MARKAL/TIMES 
                  \citep{fishbone1981markal,loulou2007etsap}
            \item eTransport 
                  \citep{bakken2007etransport}
            \item OSeMOSYS 
                  \citep{howells2011osemosys}
            \item Temoa 
                  \citep{hunter2013modeling}
            \item calliope 
                  \citep{pfenninger2015renewables}
            \item urbs 
                  \citep{dorfner2016open}
            \item ficus 
                  \citep{atabay2017open}
            \item oemof 
                  \citep{hilpert2018open}
            \item PyPSA 
                  \citep{brown2018pypsa}
            \item Switch 2.0 
                  \citep{johnston2019switch}
          \end{mylist}
        & system analysis (superstructure optimization, capacity expansion planning) for large scale national/international energy systems, typically using (MI)LP formulations
        & \begin{mylist}
            \item component-oriented modeling
            \item focus on energy systems
            \item separation of modeling and problem formulation
            \item open-source availability
          \end{mylist} \\
        \midrule
        DAMFs
        & \begin{mylist}
            \item gPROMS \citep{gPROMS}
            \item MODELICA 
                  \citep{elmqvist1997modelica}
            \item DAE Tools 
                  \citep{nikolic2016dae}
          \end{mylist}
        & detailed operational simulation; varying degrees of optimization capabilities, typically local solutions to NLP formulations
        & \begin{mylist}
            \item modeling with differential equations
            \item generic bidirectional connectivity
            \item modularity for definition of subsystems
          \end{mylist} \\
      \bottomrule
    \end{tabularx}
  \end{table*}

  Besides AMLs and ESMFs, differential-algebraic modeling frameworks (DAMFs) constitute a third class of tools that can be used to model energy systems.
  DAMFs also employ a component-oriented modeling approach, which, however, is more general than in a typical ESMF:
  In DAMFs, components may correspond to actual physical machinery or to a particular physical phenomenon (\eg, heat transfer) and can constitute subsystems, which are themselves composed of other components.
  Additionally, the information exchanged between components is not restricted to a particular kind of quantity, such as energy.
  DAMFs are particularly focused on detailed operational aspects, allowing for differential equations and nonlinear expressions within component models.
  They provide powerful features for operational \emph{simulation} of the resulting models for which a fixed design is assumed.
  Design \emph{optimization} is also possible in several DAMFs, \citep[\eg,][]{smith1997optimal, pfeiffer2012optimization}, and
  the commercial tool gPROMS \citep{gPROMS} even allows for the direct consideration of parametric uncertainty using formulations similar to Problem \eqref{prob:P}, \citep[see, \eg,][]{bansal2000simultaneous}.
  In contrast, noncommercial, open-source tools such as Open-Modelica \citep{thieriot2011towards} or Optimica \citep{akesson2010modeling} are currently limited to a single set of operational data, impeding design optimization under uncertainty.
  Furthermore, DAMFs usually offer less freedom in the choice of problem formulation, solver or algorithm in comparison to AMLs.
  In particular, many tools employ gradient-based methods, \citep[\eg,][]{
      pfeiffer2012optimization,
      navarro2014computer,
      magnusson2015dynamic}
  yielding only local solutions, or heuristic global optimization methods, \eg, random search, genetic algorithms, or simulated annealing \citep{
      thieriot2011towards,
      pfeiffer2012optimization,
      kim2018efficient},
  which treat the system model as a black box and cannot reliably locate global solutions.

  AMLs, ESMFs and DAMFs each exhibit strengths related to a particular aspect of modeling and optimizing energy system design and operation.
  ESMFs are tailored to energy systems modeling and offer a component-oriented approach that benefits model maintenance and re-usability.
  However, their principal focus is on system analysis.
  In particular, their restriction to LP or MILP formulations makes them less suited for applications concerned with technical design and operation.
  Both AMLs and DAMFs lift the restriction to (MI)LP formulations, but AMLs lack high-level component-oriented abstractions for generic energy systems and DAMFs lack control over the choice of problem formulation and optimization algorithm.
  We therefore propose a next-generation ESMF that allows for flexible, component-oriented modeling, including nonlinear and differential-algebraic formulations, parametric uncertainty, and the possibility to specify specialized solution algorithms.
  Its basic structure is presented in the following Section.
  A summary of the above discussion, highlighting the roles of each tool class and their influence on COMANDO is given in \cref{tab:tool_classes}.

  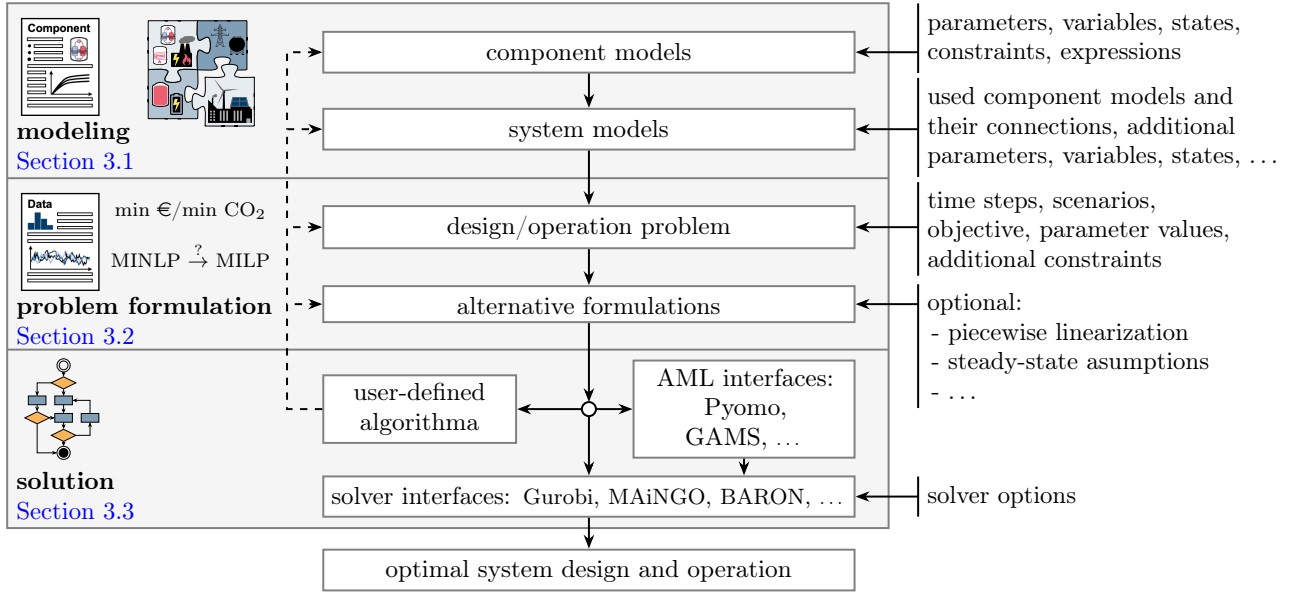
\begin{figure*}[!bt]
    \centering
    \hspace*{1cm}
    \begin{adjustbox}{width=\textwidth}
       \input{figures/workflow}
    \end{adjustbox}
    \vspace*{-2.7mm}
    \caption{Workflow for modeling, problem formulation, and optimization using COMANDO.}
    \label{fig:workflow}
  \end{figure*}

%% file: figures/workflow.tex
\begin{tikzpicture}[thick, >={Stealth[black, length=2mm]}]
  \def\top{2.5}
  \def\lef{-6.5}
  \def\rig{6}
  \def\h{2.5}
  \def\dx{8.25}
  \def\tw{7.3}
  \def\tow{6.8}
  \def\ws{0}
  \def\tc{black}
  \def\fc{gray}

  \fill[fill=\bc]  (\rig, \top) rectangle (\lef, \top-\h) node (v1) {};
  \node[align=left, anchor=south west, color=\tc] at (v1) {\bf{modeling} \\ \cref{sec:modeling}};
  \node[inner sep=0pt] (russell) at ($(v1) + (0.75, 1.6)$)
    {\includegraphics[height=1.4cm]{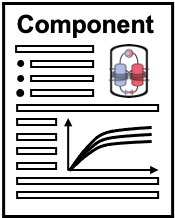}};
  \node[draw=\fc, rectangle, fill=white, text width=\tw cm, text height=0.28cm, align=center, xshift=\dx cm, yshift=1.8cm] (CM) at (v1)
  {component models};
  \node[right of=CM, xshift=\tow cm, yshift=0.2cm, text width=6cm] (CMT) {parameters, variables, states, \\ constraints, expressions};
  \draw (CMT.south west) -- (CMT.north west);

  \node[draw=\fc, rectangle, fill=white, text width=\tw cm, text height=0.28cm, align=center, xshift=\dx cm, yshift=0.7cm] (SM) at (v1)
  {system models};
  \node[anchor=west, right of=SM, xshift=\tow cm, yshift=0.025cm, text width=6cm] (SMT) {used component models and \\ their connections, additional \\ parameters, variables, states, \dots};
  \draw  (SMT.south west) -- (SMT.north west);

  \fill[fill=\bc] (\rig, \top - \h - \ws) rectangle (\lef, \top - 2 * \h - \ws) node (v2) {};
  \node[align=left, anchor=south west, color=\tc] at (v2) {\bf{problem formulation} \\
  \cref{sec:problem_formulation}};
    \node[inner sep=0pt] (russell) at ($(v2) + (0.75, 1.6)$)
    {\includegraphics[height=1.4cm]{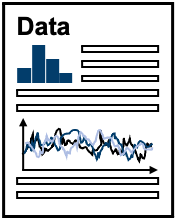}};
  \node[align=center, anchor=north] at ($(v2) + (2.6, 2.3)$) {\footnotesize min \euro\slash min CO$_2$ \\[3mm]
  \footnotesize
  MINLP $\overset{?}{\rightarrow}$ MILP};
  \node[draw=\fc, rectangle, fill=white, text width=\tw cm, text height=0.28cm, align=center, xshift=\dx cm, yshift=1.8cm] (PF) at (v2)
  {design\slash operation problem};
  \node[right of=PF, xshift=\tow cm, yshift=-0.025cm, text width=6cm] (PFT) {time steps, scenarios, \\ objective, parameter values, \\ additional constraints};
  \draw  (PFT.south west) -- (PFT.north west);

  \node[draw=\fc, rectangle, fill=white, text width=\tw cm, text height=0.28cm, align=center, xshift=\dx cm, yshift=0.7cm] (AF) at (v2)
  {alternative formulations};
  \node[anchor=west, right of=AF, xshift=\tow cm, yshift=-0.6cm, text width=6cm] (AFT) {optional: \\ \,- piecewise linearization \\ \,- steady-state asumptions \\ \,- \dots};
  \draw  (AFT.south west) -- (AFT.north west);

  \fill[fill=\bc] (\rig, \top - 2 * (\h + \ws) rectangle (\lef, \top - 3 * \h - 2* \ws) node (v3) {};
  \node[align=left, anchor=south west, color=\tc] at (v3) {\bf{solution} \\ \cref{sec:optimization}};
  \node[circle, draw, fill=white, inner sep=0pt, outer sep=0pt, minimum size=0.2cm, xshift=\dx cm, yshift=1.7cm] (C) at (v3) {};

  \node[draw=\fc, rectangle, fill=white, text width=2.5cm, align=center, xshift=5.85cm, yshift=1.7cm] (UD) at (v3)
  {user-defined \\
  algorithma
  };

  \node[draw=\fc, rectangle, fill=white, text width=2.9cm, align=center,
  xshift=10.45 cm , yshift=1.7cm
  ] (IF) at (v3)
  {AML interfaces: \\
  Pyomo,
  GAMS,
  \scriptsize \dots
  };

  \node[draw=\fc, rectangle, fill=white, text width=\tw cm, text height=0.28cm, align=left, xshift=\dx cm, yshift=.45cm] (SV) at (v3)
  {solver interfaces:
    {\small
        Gurobi,
        MAiNGO,
        BARON,
        \scriptsize \dots
    }
  };
  \node[anchor=west, text width=6cm] at (AFT.west|-SV)
    (SVT)
  {solver options \\ };
  \draw  (SVT.south west) -- (SVT.north west);

  \node[draw=\fc, rectangle, fill=white, text width=\tw cm, text height=0.28cm, align=center, xshift=\dx cm, yshift=-.6cm] (PP) at (v3)
  {optimal system design and operation};

  \def\bot{\top - 3 * \h - 2.5 * \ws}
  \draw[\fc, thick] (\lef, \top) rectangle (\rig, \bot);
  \draw[\fc, thick] (v1) ++(0, -\ws/2) -- ++(0:\rig - \lef);
  \draw[\fc, thick] (v2) ++(0, 0.05-\ws/2) -- ++(0:\rig - \lef);

  \draw[<-] (CM.east) -- (CM-|CMT.west);

  \draw[<-] (SM.east) -- (SM-|SMT.west);
  \draw[<-] (PF.east) -- (PF-|PFT.west);
  \draw[<-] (AF.east) -- (AF-|AFT.west);
  \draw[<-] (SV) -- (SVT);

  \draw[->] (CM) -- (SM);
  \draw[->] (SM) -- (PF);
  \draw[->] (PF) -- (AF);
  \draw[->] (AF) -- (C);
  \draw[->] (C) -- (C-|IF.west);
  \draw[->] (C.west) -- (C.west-|UD.east);
  \coordinate[right of=UD, xshift=22mm] (UDE);

  \draw[->] (C) -- (C|-SV.north);
  \draw[->] (IF.south) -- (IF.south|-SV.north);
  \draw[->, dashed] (UD.west) -- ++(-0.5cm, 0cm) |- node (b1) {} (AF);
  \draw[->, dashed] (b1) |- node (b2) {} (PF);
  \draw[->, dashed] (b2) |- node (b3) {} (SM);
  \draw[->, dashed] (b3) |- (CM);
  \draw[->] (SV) -- (PP);

  \begin{scope}[
		shift={(-4.5cm, 1.55cm)},
		scale=.7, transform shape, thin]
    \begin{scope}
      \color{black}
      \piece[fzjblue!10]{-1}{-1}{0}{0}
			\node at (0.7, 0.38) {\includegraphics[scale=0.04]{figures/Icons/CHP}};
			\node at (0.4, 0.7) {\includegraphics[scale=0.025]{figures/Icons/HP}};
			\node at (0.25, 0.3) {\includegraphics[scale=0.03]{figures/Icons/B}};
    \end{scope}
    \begin{scope}[xshift=1cm]
      \color{black}
      \piece[fzjblue!50]{1}{0}{0}{1}
			\node at (0.45, 0.7) {\includegraphics[scale=0.08]{figures/Icons/PG}};
			\node at (0.8, 0.5) {\includegraphics[scale=0.08]{figures/Icons/GG}};
    \end{scope}
    \begin{scope}[yshift=-1cm]
      \color{black}
      \piece[fzjblue!30]{0}{1}{1}{0}
			\node at (0.25, 0.5) {\includegraphics[scale=0.08]{figures/Icons/TES_H}};
			\node at (0.58, 0.3) {\includegraphics[scale=0.12]{figures/Icons/BAT}};
    \end{scope}
    \begin{scope}[xshift=1.15cm, yshift=-1.15cm]
      \color{black}
      \piece[fzjblue!10]{0}{0}{-1}{-1}
            \node at (0.44, 0.56) {\includegraphics[scale=0.06]{figures/Icons/WT}};
			\node at (0.7, 0.46) {\includegraphics[scale=0.04]{figures/Icons/PV}};
			\node at (0.5, 0.36) {\includegraphics[trim=122 0 0 0,clip, scale=0.08]{figures/Icons/CON}};
    \end{scope}
  \end{scope}

  \begin{scope} [shift={(-5.7cm, -3.55cm)}, scale=0.125, >={Stealth[black, scale=0.4]}, thin, transform shape]
    \node[circle, draw, minimum size=1.5cm] (c1) at (0, 7) {};
    \node[circle, draw, minimum size=1cm] (cc1) at (0, 7) {};
    \node (d1) at (0, 5) [draw, diamond, aspect=2, scale=4, fill=fzjorange] {};

    \node[rectangle, draw, fill=fzjblue!50, minimum width=2cm, minimum height=1cm] (r1) at (0, 3) {};
    \node[rectangle, draw, fill=fzjblue!50, minimum width=2cm, minimum height=1cm] (r2) at (-3, 3) {};
    \node[rectangle, draw, fill=fzjblue!50, minimum width=2cm, minimum height=1cm] (r3) at (0, 1) {};
    \node[rectangle, draw, fill=fzjblue!50, minimum width=2cm, minimum height=1cm] (r4) at (3, 1) {};

    \node (d2) at (-3, 1) [draw, diamond, aspect=2, scale=4, fill=fzjorange] {};

    \node (d3) at (0, -1) [draw, diamond, aspect=2, scale=4, fill=fzjorange] {};

     \node[circle, draw, minimum size=1.5cm] (c2) at (0, -3) {};
     \node[circle, draw, fill, minimum size=1cm] (cc2) at (0, -3) {};

     \draw[->] (c1) -- (d1);
     \draw[->] (d1) -- (r1);
     \draw[->] (r1) -- (r3);
     \draw[->] (r3) -- (d3);
     \draw[->] (d3) -- (c2);

     \draw[->] (d1) -| (r2);
     \draw[->] (r2) -- (d2);
     \draw[->] (d2) |- (c2);

     \draw[->] (d2) -- (r3);

     \draw[->] (d3) -| (r4);
     \draw[->] (r4) |- (r1);
  \end{scope}
\end{tikzpicture}

%% file: sections/3_COMANDO.tex

\section{The COMANDO ESMF} \label{sec:COMANDO}

The goal of COMANDO is to provide an open-source ESMF which allows to generate detailed models of energy system components, including differential-algebraic and nonlinear elements, and aggregate them to system models for the purpose of optimization.
Traditional ESMFs are typically oriented towards techno-economic analysis of systems at national or international scales, where (MI)LP formulations are an asset that ensures computational tractability.
In contrast, COMANDO is oriented towards the technical design and operation of small- to medium-scale systems, e.g., district energy systems, industrial sites, or energy conversion processes.
At these scales, investigation of realistic component and system behavior is possible via the consideration of technically relevant effects such as part-load and dynamic behavior.

Unlike most ESMFs, which are commonly based on an AML, COMANDO is implemented as a flat layer on top of the computer algebra system SymPy \citep{meurer2017sympy}.
This choice provides: i) data structures for the mathematical expressions used to describe components and systems, as well as ii) several routines useful for creating automatic reformulations and user-defined algorithms, such as automatic differentiation, substitution of expressions or solution of nonlinear systems of equations.
As a modeling framework, COMANDO itself does not provide any specialized solution methods.
Instead, it allows for component-oriented modeling at a high level of abstraction, while at the same time granting users access to low-level data structures.
This allows for both intuitive modeling, as well as advanced use cases such as problem reformulations and the development of user-defined algorithms, simplifying the development of tailored solution approaches.
During the development of COMANDO, we made an effort to maximize chances of its adoption by following best-practices for code development.
This includes the creation of automated unit and integration tests, provision of documentation (both in the source code and as a standalone document \citep{COMANDO_DOCS}), and the inclusion of the full code for running the case studies as detailed usage examples.
In order to further encourage adoption, we provide a generic parsing routine to translate individual COMANDO expressions to alternative textual or object-oriented representations, allowing users to easily link COMANDO to other software.
An overview of the structure of COMANDO and the typical workflow of modeling, problem formulation and optimization is given in \cref{fig:workflow}.

In \cref{sec:modeling} we describe the process of creating models for components and systems in COMANDO.
\cref{sec:problem_formulation} provides details on how optimization problems can be created from a system model and how alternative formulations of these problems can be obtained.
Finally, the different options for solving the formulated problems are given in \cref{sec:optimization}.

\subsection{Modeling process} \label{sec:modeling}
  The goal of the modeling phase is to generate a model describing the behavior of a given energy system.
  For the creation of such a system model, models for its constituting components as well as information on their connectivity are required.

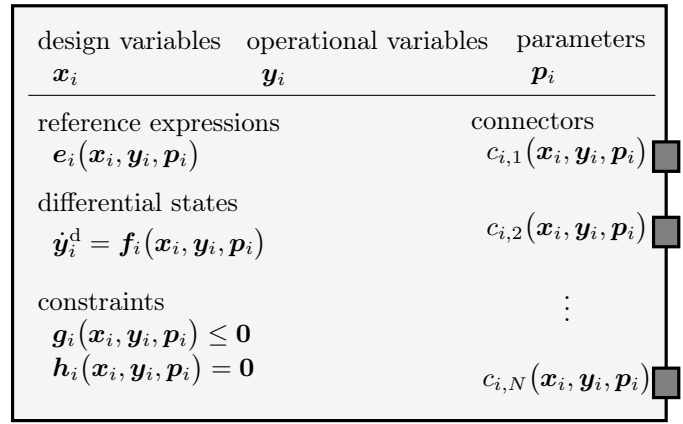
\begin{figure}[t!]
    \centering
    \input{figures/component}
    \vspace*{-1mm}
    \caption{
      Structure of a generic component $i$ in COMANDO.
      Mathematical expressions are specified based on symbols that are either parameters, design variables, or operational variables.
      These expressions can be kept for later reference, constitute the right-hand side of differential equations, form part of algebraic constraints, or describe possible in- and\slash or outputs through connectors.
    }
    \label{fig:component}
  \end{figure}

  We begin with the description of component models, which are used to represent elementary parts of an energy system.
  \cref{fig:component} depicts the structure of the \texttt{Component} class used for that purpose.
  A model of a component $i$ consists of several types of mathematical expressions, given in symbolic form.
  Following the notation introduced in \cref{sec:M_AND_O}, the expressions describing the component contain different symbols corresponding to quantities which are either \emph{parameters} ($\bm{p}_i$), \ie, placeholders for values that are assumed to be given before an optimization, or design or operational \emph{variables} ($\bm{x}_i$ and $\bm{y}_i$, respectively), \ie, placeholders for scalar and vector values that are to be determined during optimization.

  To instantiate a \texttt{Component}, a unique name must be provided, which serves as an identifier for the component.
  The names of parameters, variables, and constraints associated to the component are prepended with this identifier, in order to distinguish quantities from different instances of the same component model.
  The \texttt{Component} class can either be used directly or subclassed to specify specialized component classes with custom behavior.
  To create and add symbols to a component, the \texttt{Component} class provides three methods:
  \begin{itemize}[noitemsep]
    \item \texttt{make\_parameter},
    \item \texttt{make\_design\_variable}, and
    \item \texttt{make\_operational\_variable}.
  \end{itemize}
  All three methods require a name for the symbol that is used to represent the quantity.
  The methods for creating variables provide optional arguments for the specification of variable bounds, domain (integer/real) and a scalar value for initialization, while the parameter creation method only provides a single optional argument for the specification of its value.
  Note that time- and scenario-specific values for operational quantities are set in the problem formulation phase after the time and scenario structure, has been specified, see \cref{sec:problem_formulation}.

  Based on variables and parameters, mathematical expressions can be formed using the overloaded Python operators \texttt{+}, \texttt{-}, \texttt{*}, \texttt{/}, \texttt{**}, or any of the functions implemented in SymPy (\eg, exp, log, trigonometric, and hyperbolic functions).
  Any intermediate expressions $\bm{e}_i$ that are of interest can be assigned an identifier and stored in the component using the \texttt{add\_expression} method.
  These expressions can simply be used for evaluation or as parts of more complex expressions, \eg, system-level constraints, or an objective function, \cf \cref{sec:problem_formulation}.
  Vectors $\bm{g}_i$ and $\bm{h}_i$ contain inequality and equality constraints associated to the component $i$ and their elements can be specified using the methods
  \begin{itemize}[noitemsep]
    \item \texttt{add\_le\_constraint},
    \item \texttt{add\_eq\_constraint}, and
    \item \texttt{add\_ge\_constraint}.
  \end{itemize}
  Each of these methods takes two expressions and an optional name for the resulting relation as arguments.
  Explicit distinction into first and second stage expressions and constraints is not necessary and occurs automatically, based on the symbols present in the respective expressions.

  Dynamic behavior can be represented by specifying right-hand side expressions $\roc_i$ for the time derivatives of differential states $\bm{y}^\text{d}_i$ (recall that $\bm{y}^\text{d}_i$ constitutes a subset of the operational variables $\bm{y}_{i}$).
  Previously created operational variables may be declared differential states using the \texttt{declare\_state} method or differential states may be created directly using the \texttt{make\_state} method.
  The first method requires an existing variable and an expression, corresponding to entries of the vectors $\bm{y}^\text{d}_i$ and $\bm{f}_i$
  as mandatory arguments and allows for the specification of an initial state as well as bounds and an initial guess for the value of the derivative.
  The method results in the creation of a new operational variable, corresponding to an element in $\bm{\dot{y}}^\text{d}_i$, and an equality constraint, linking the time derivative with the given expression in $\bm{f}_i$.
  An explicit relation between the state and its derivative is not specified at this point, as it depends on the desired time-discretization which is handled by the solution interfaces, \cf \cref{sec:problem_formulation}.
  The \texttt{make\_state} method creates a new operational variable corresponding to the differential state and then calls \texttt{declare\_state}.

  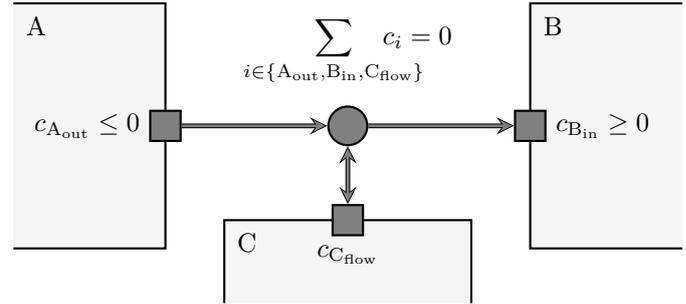
\begin{figure}[t!]
    \centering
    \input{figures/connector}
    \vspace*{-3mm}
    \caption{
      A connection formed by connecting three connectors to a bus: The components A, B, and C each define a connector for a particular quantity. The connectors of A and B are marked as outputs and inputs, respectively, restricting the sign of the associated expression, while the connector of C is not restricted.
      The connection of $\text{A}_\text{out}, \text{B}_\text{in}$, and $\text{C}_\text{flow}$ via a bus results in the creation of a balance constraint in the system model.
      This graphical notation is also used for the case studies in \cref{sec:case_study}.
    }
    \label{fig:connector}
  \end{figure}

  To allow for the aggregation of components to systems, individual expressions in $\bm{c}_{i}$ can be assigned to connectors (\cf \cref{fig:component}).
  Connectors are generally bidirectional, but may be specified to only allow for in-, or output.
  In- and output connectors restrict the assigned expression to a nonnegative or nonpositive range, respectively.

  A system model can be created as an instance of the \texttt{System} class, whose instantiation again requires a unique label that serves as an identifier.
  Optionally, a list of components and the connections between them can be passed to the constructor of the \texttt{System} class.
  Each connection is specified via a label and a list of associated connectors.
  The connectors are connected to a `bus' at which the quantities associated to them are balanced and a corresponding constraint is created automatically, see the graphical notation in \cref{fig:connector}, which is also used for the case studies in \cref{sec:case_study}.
  The elementary connections provided by COMANDO's \texttt{System} class, handle only simple balance equations.
  More complicated connectivities such as mixing streams with different temperatures, concentrations or other qualities are most naturally implemented as a dedicated component within COMANDO.
  Instead of specifying the complete structure during construction of a \texttt{System} instance, components and connections can also be added sequentially via corresponding methods, allowing for procedural model generation.
  As in DAMFs, a nested creation of systems from subsystems is possible by exposing connectors of individual components or extending existing connections via additional connectors.
  For instance, a neighborhood can be represented as a system composed of buildings as subsystems, which are in turn composed of heating, cooling and power equipment.
  As with component models, system models can be assigned their own variables, parameters, expressions and constraints describing their behavior.
  These two features are accomplished by letting the \texttt{System} class inherit from the \texttt{Component} class.
  The system superstructure can be considered explicitly by including appropriate design decisions within component models.
  More advanced approaches where the superstructure is not specified a-priori, e.g., superstructure-free synthesis \citep{voll2012superstructure}, or automated superstructure generation and expansion \citep{voll2013automated}, can be easily incorporated in the form of user-defined algorithms.

\subsection{Problem formulation} \label{sec:problem_formulation}
  Based on a system model, different kinds of optimization problems considering system design and\slash or operation can be created.
  To this end, COMANDO provides the \texttt{Problem} class, instances of which can be created by the \texttt{create\_problem} method of the \texttt{System} class.
  As the system model defines a constraint set which is parameterized by the parameters $\bm{p}$, only the objective terms $\dobj$ and $\moobj$ as well as a time and scenario structure and appropriate data (\ie, values for the parameters $\bm{p}$) need to be specified in the \texttt{create\_problem} method to obtain a complete problem formulation, corresponding to \eqref{prob:P}.
  Note that the user may decide which units to use for data and time steps but must ensure they match.
  Units given in the \nameref{sec:nomenclature} are those used for the case studies in \cref{sec:case_study}.

  To define the objective terms, the \texttt{System} class provides the \texttt{aggregate\_component\_ex\-pres\-sions} method.
  For a given expression identifier, it returns the sum of all expressions stored under that identifier in the individual components.
  The resulting expressions can be used for the objective terms $\dobj$ and $\moobj$, depending on whether they consist exclusively of first stage (\ie, scalar) quantities or not.
  A second use for the \texttt{aggregate\_component\_expressions} method is to create expressions for system-level constraints involving contributions from multiple components.

\begin{figure}[t!]
    \centering
    \resizebox{\columnwidth}{!}{%
      \input{figures/timesteps}
    }
    \vspace*{-5mm}
    \caption{
      Alternative ways to specify time steps for a particular scenario: For variable length an ordered mapping (left) and for constant length a list and the total length (right) can be specified.
      If multiple scenarios with different time structures are to be considered, one such description is given per scenario.
    }
    \label{fig:timesteps}
    \vspace*{-3mm}
  \end{figure}
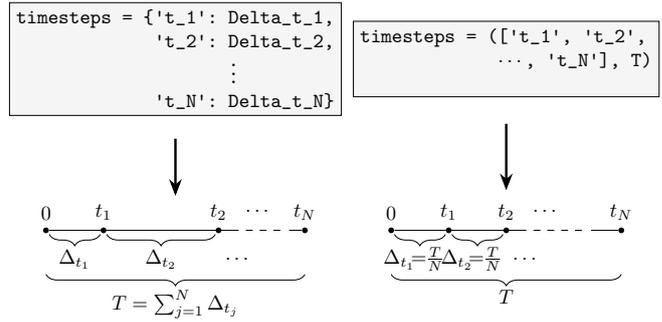

  The time and scenario structure is specified in terms of the considered scenarios $s \in \mathcal{S}$ and the corresponding discretized time horizons $\Tdisc_s$.
  The $\Tdisc_s$ are required by COMANDO's solver or AML interfaces for the automatic discretization of the differential equations.
  If more than one operational scenario is considered, the different scenarios can either be specified as a list of $M$ scenario identifiers, corresponding to scenarios with probability $1/M$, or by a series of scenario identifiers and associated weights $w_s$.
  In the latter case, the weights are not required to sum to one, allowing for a more general weighting.
  Similarly, individual time points for each time horizon are either specified via a mapping of time point labels $t$ to the corresponding lengths $\Delta_{s,t}$ or in the case of equidistant time steps via a list of labels and an end-time $T_s$, see \cref{fig:timesteps}.
  If the time horizons are identical for all scenarios, a single time horizon can be specified, otherwise, one specification per scenario is required.

  Parameter values corresponding to the resulting time and scenario structure can be specified during problem creation and may later be updated using the \texttt{data} attribute of the \texttt{Problem} instance.
  Similarly, design and operational variable values can be updated using the \texttt{design} and \texttt{operation} attributes, respectively.
  Values for design variables must be scalar while values for parameters and operational variables may be provided as scalars or as time- and\slash or scenario-dependent data.

  After the abovementioned steps, a problem in the form of \eqref{prob:P} is fully specified.
  However, it may be desirable to adapt the original problem formulation in different ways.
  Adaptations to the problem formulation range from simply adding further constraints to the reformulation of expressions in the problem.
  One generic reformulation routine implemented in COMANDO is the automatic linearization of arbitrary continuous multivariate expressions via convex-combination or multiple-choice linearization \citep{vielma2010mixed}.
  More generally, custom reformulations may be created making use of existing algorithms provided by SymPy \citep{meurer2017sympy}, \eg, for automatic differentiation, analytic solution of different kinds of nonlinear equation systems, or symbolic substitution of subexpressions.
  Note that reformulations do not have to result in approximations but can also be used to create alternative formulations that possess better properties than the original one, \eg, tighter relaxations for deterministic global optimization.

\subsection{Problem solution} \label{sec:optimization}
  A fully specified problem formulation can be directly passed to a suitable solver or to an AML.
  In this step, the problem structure and data are translated from the COMANDO representation to a new representation, matching the syntax of the target solver or AML.
  For this purpose, COMANDO contains a generic parsing routine that can be used to create new interfaces based on target-specific representations of the symbols and operations occurring within the different expressions of the problem formulation.
  Interfaces may be text-based, resulting in an input file for a solver or AML, or they can be object-oriented, resulting in a translation of the problem formulation using the target-API.
  Currently implemented interfaces are:
  \begin{itemize}[itemsep=0mm]
    \item text-based:
    \begin{itemize}[itemsep=0mm, topsep=0mm]
      \item BARON \citep{baron20_10_16} (solver)
      \item GAMS \citep{bussieck2004general} (AML)
      \item MAiNGO \citep{bongartz2018maingo} (solver)
    \end{itemize}
    \item API-based:
    \begin{itemize}[itemsep=0mm, topsep=0mm]
      \item Pyomo \citep{hart2011pyomo} (AML)
      \item Pyomo.DAE \citep{nicholson2018pyomo.dae} (AML)
      \item Gurobi \citep{gurobi9_1_1} (solver)
      \item MAiNGO \citep{bongartz2018maingo} (solver)
    \end{itemize}
  \end{itemize}
  All of these interfaces provide methods to solve the deterministic equivalent formulation of Problem \eqref{prob:P} with a given set of options, and to write back the obtained results to COMANDO.
  Note that a problem formulation may contain differential equations if states were defined in the component or system model.
  Since most solvers and AMLs do not support differential equations, the corresponding interfaces can specify different schemes for automatic time discretization.
  All existing interfaces implement implicit Euler discretization.
  More advanced schemes are available through the Pyomo.DAE interface.

  Instead of directly solving a problem, it can also be addressed with a user-defined algorithm.
  User-defined algorithms can range from simple preprocessing routines based on the system model and available data to more advanced methods, such as decomposition techniques, commonly used in stochastic programming \citep[see, \eg,][]{li2019generalized}.
  The architecture of COMANDO allows for manipulation at the level of component and system models as well as at the level of the resulting optimization problems.
  In particular the \texttt{Problem} class can also be used to specify the sub-problems that may occur within user-defined algorithms, allowing them to be passed to any of the available interfaces.

  \begin{table*}[htb!]
    \centering
    \caption{Overview of the presented case studies}
    \label{tab:case_studies}
      \scriptsize
    \begin{tabularx}{1\linewidth}{
      >{\raggedright\hsize=0.8\hsize\linewidth=\hsize}X
      >{\raggedright\hsize=1.0\hsize\linewidth=\hsize}X
      >{\centering\hsize=0.8\hsize\linewidth=\hsize}X
      >{\centering\hsize=0.4\hsize\linewidth=\hsize}X
      >{\raggedright\hsize=1.35\hsize\linewidth=\hsize}X
      >{\hsize=1.6\hsize\linewidth=\hsize}X
      }
      \toprule
      \thead{case study}
        & \thead{system structure}
        & \thead{problem class \\ (reformulations)}
        & \thead{problem \\ type}
        & \thead{operational horizon \\ representation}
        & \thead{demonstrated features}
        \\
      \midrule
      industrial energy system (\cref{sec:design_problem})
        & superstructure with CHP subsystem, 15 instances of 11 component classes
        & \makecell[ct]{MINLP \\ (MILP, NLP)}
        & design
        & 6 scenarios representing 4 typical days with 4 time steps of varying length, each, and two isolated time steps representing extreme demands, implicit Euler time discretization
        & \begin{mylist}
            \item superstructure optimization
            \item automatic linearization
            \item user-defined algorithm
            \item re-use of model for multiple problem formulations
          \end{mylist}
        \\
      \midrule
      building demand response (\cref{sec:dynamic_operation})
        & 9 instances of 4 component classes
        & MIDO (MILP)
        & operation
        & 24\,h horizon with 15\,min time steps, each with 4 collocation points
        & \begin{mylist}
            \item abstract components
            \item modeling with differential equations
            \item advanced time discretization via collocation
                  in Pyomo.DAE
                  \citep{nicholson2018pyomo.dae}
          \end{mylist}
        \\
      \midrule
      low-temperature district heating network (\cref{sec:HP_net})
        & superstructure with 9 instances of 2 subsystem classes (linking, consumer group), 26 instances of 6 component classes
        & MIQCQP
        & design
        & 11 scenarios, each representing a static operating point
        & \begin{mylist}
            \item modular model generation
            \item superstructure optimization
            \item stochastic programming
          \end{mylist}
        \\
      \midrule
      organic Rankine cycle (\cref{sec:ORC})
        & 8 instances of 4 component classes
        & NLP
        & operation
        & single operating point
        & \begin{mylist}
            \item hybrid modeling with ANNs
            \item reduced space formulation
            \item integration with different solver/AML interfaces
          \end{mylist}
        \\
      \bottomrule
    \end{tabularx}
  \end{table*}

%% file: figures/component.tex
\usetikzlibrary{positioning}
\def\myargs{\big(\bm{x}_i, \bm{y}_i, \bm{p}_i\big)}
\begin{tikzpicture}
  \coordinate (v1) at (-5, 5);
  \coordinate (v2) at (3.6, -0.5);
  \draw[very thick, fill=\bc] (v1) rectangle (v2);

  \node[below right of=v1, anchor=north west, xshift=-5mm, yshift=5mm, align=left] (DV) {
    design variables \\
    $\;\;\bm{x}_i$
  };
  \node[below right of=v1, anchor=north west, xshift=22.5mm, yshift=5mm, align=left] (OV) {
    operational variables \\
    $\;\;\bm{y}_i$
  };
  \node[below right of=v1, anchor=north west, xshift=58mm, yshift=5mm, align=left] (P) {
    parameters \\
    $\;\;\bm{p}_i$
   };

  \draw ($(DV.west)!(OV.south west)!(DV.south west)$) -- ($(P.east)!(OV.south west)!(P.south east)$);

  \node[below=of DV.west, anchor=north west, align=left, yshift=4mm] (E) {
    reference expressions \\
    $\;\;\bm{e}_i\myargs$
  };
  \node[below=of E.west, anchor=north west, yshift=4.5mm, align=left] (S) {
    differential states \\[1mm]
    $\;\;
    \dot{\bm{y}}^\text{d}_i = \roc_i\myargs $
  };
  \node[below=of S.west, anchor=north west, yshift=2.5mm, align=left] (CONS) {
    constraints \\
    $\;\;\bm{g}_i\myargs \leq \bm{0}$\\
    $\;\;\bm{h}_i\myargs = \bm{0}$
  };

  \node[rectangle, draw, very thick, fill=black!50, yshift=35mm, label={[xshift=-13mm, yshift=-5mm]$c_{i, 1}\myargs$}] (tc) at (v2) {
    \phantom{$\cdot$}
  };
  \node[rectangle, draw, very thick, fill=black!50, yshift=25mm, label={[xshift=-13mm, yshift=-5mm]$c_{i, 2}\myargs$}] (mc) at (v2) {
    \phantom{$\cdot$}
  };
   \node[yshift=15mm, label={[xshift=-13mm, yshift=-5mm]$\vdots$}] (dc) at (v2) {
     \phantom{$\cdot$}
   };
  \node[rectangle, draw, very thick, fill=black!50, yshift=5mm, label={[xshift=-13mm, yshift=-5mm]$c_{i, N}\myargs$}] (bc) at (v2) {
    \phantom{$\cdot$}
  };

  \node[anchor=north east, align=left, xshift=-7mm] (CONN) at ($(tc.center)!(E.north)!(v2)$) {
    connectors
  };
\end{tikzpicture}

%% file: figures/connector.tex
\usetikzlibrary[arrows.meta]
\begin{tikzpicture}[scale=1, >={Stealth[black]}]
  \def\depth{2cm}
  \def\dist{2.4cm}
  \def\cdist{1.25cm}
  \def\cdepth{1.1cm}
  \def\width{3.25 cm}
  \def\offset{0.3cm}

  \def\lw{0.4pt}
  \def\aw{1pt}
  \def\comm{black!50}
  \tikzset{arrow filled/.style args={#1}{
              -stealth, >=Stealth,line width=2*\lw+\aw, black,
              postaction={draw, -stealth, >=Stealth, #1,
                          line width=\aw,
                          shorten >=2*(\lw)
                          }
             }
  }
  \tikzset{double arrow filled/.style args={#1}{
              <->, >=stealth, line width=2*\lw+\aw, black,
              postaction={draw, <->, >=stealth, #1,
                          line width=\aw,
                          shorten <=2*(\lw),
                          shorten >=2*(\lw)}
             }
  }
  \tikzset{connector/.style={rectangle, draw, thick, fill=\comm, inner sep=0, minimum size=4mm}}
  \node[circle, draw, thick, fill=\comm, minimum size=5mm] (bus) at (0,0) {
    \phantom{$\cdot$}
  };
  \node[above of=bus] {$\displaystyle\sum\limits_{i \in \{\text{A}_\text{out}, \text{B}_\text{in}, \text{C}_\text{flow}\}} \hspace{-6mm} c_{i} = 0$};

  \coordinate (c1) at (-\dist - \depth, 0.5 * \width);
  \coordinate (c2) at (-\dist, -0.5 * \width);
  \fill[\bc]  (c1) rectangle (c2);
  \draw[thick] (c1) -- ++(\depth, 0) -- (c2) -- ++(-\depth, 0);
  \node[xshift=\offset, yshift=-\offset] at (c1) {A};
  \node[connector, label=left:{$c_{\text{A}_\text{out}} \leq 0$}] (l) at (-\dist,0) {
    \phantom{$\cdot$}
  };

  \coordinate (c3) at (\dist + \depth, -0.5 * \width);
  \coordinate (c4) at (\dist, 0.5 * \width);
  \fill[\bc]  (c3) rectangle (c4);
  \draw[thick] (c3) -- ++(-\depth, 0) -- (c4) -- ++(\depth, 0);
  \node[xshift=\offset, yshift=-\offset] at (c4) {B};
  \node[connector, label=right:{$c_{\text{B}_\text{in}} \geq 0$}] (r) at (\dist,0) {
    \phantom{$\cdot$}
  };

  \coordinate (c5) at (0.5 * \width, - \cdist - \cdepth);
  \coordinate (c6) at (-0.5 * \width, -\cdist);
  \fill[\bc]  (c5) rectangle (c6);
  \draw[thick] (c5) -- ++(0, \cdepth) -- (c6) -- ++(0, -\cdepth);
  \node[xshift=\offset, yshift=-\offset] at (c6) {C};
  \node[connector, label=below:{$c_{\text{C}_\text{flow}}$}] (b) at (0,-\cdist) {
    \phantom{$\cdot$}
  };
  \draw[arrow filled=gray] (l) -- (bus);
  \draw[arrow filled=gray] (bus) -- (r);
  \draw[double arrow filled=gray] (b) -- (bus);
\end{tikzpicture}

%% file: figures/timesteps.tex
\usetikzlibrary{arrows.meta, decorations.pathreplacing}
\begin{tikzpicture}[>={Stealth[black]}]
  \tikzset{time point/.style args={#1}{circle, fill=black, inner sep=0, minimum size=1mm, label=#1}}
  \begin{scope}
	  \node[circle, fill=black, inner sep=0, minimum size=1mm, label=0] (t0) at (0, 0) {};
	  \node[time point=$t_1$] (t1) at (1, 0)  {};
	  \node[time point=$t_2$] (t2)  at (3, 0) {};
	  \node[label=$\cdots$] (dots)  at (3.7, 0) {};
	  \node[time point=$t_N$] (tN)  at (4.5, 0) {};

	  \draw (t0) -- ($(t2)!0.5!(dots)$) node (ld) {};
	  \draw[dash pattern={on 1.2mm off 1.2mm}] (ld) -- ($(dots)!0.5!(tN)$) node (rd) {};
	  \draw (rd) -- (tN);
	  \draw[decoration={brace, raise=2pt, mirror, amplitude=2mm}, decorate] (t0) -- (t1);
	  \draw[decoration={brace, raise=2pt, mirror, amplitude=2mm}, decorate] (t1) -- (t2);

	  \draw[decoration={brace, raise=2pt, mirror, amplitude=2mm}, decorate] (t0) +(0, -0.7) -- node[yshift=-1mm, label=below:{$T = \sum_{j=1}^N \Delta_{t_j}$}] {} ($(tN) +(0, -0.7)$);

	  \node[yshift=-5mm] at ($(t0)!0.5!(t1)$) {$\Delta_{t_1}$};
	  \node[yshift=-5mm]  at ($(t1)!0.5!(t2)$) {$\Delta_{t_2}$};
	  \node[yshift=-5mm]  at (ld) {$\cdots$};
	  \coordinate (center) at ($(t0)!0.5!(tN)$);
	  \node[rectangle, draw, fill=\bc, yshift=3cm, align=left] (code) at (center)
	  {\texttt{timesteps = \{'t\_1':\ Delta\_t\_1,} \\
	   \texttt{\ \ \ \ \ \ \ \ \ \ \ \ \ 't\_2':\ Delta\_t\_2,} \\
	   \texttt{\ \ \ \ \ \ \ \ \ \ \ \ \ \ \ \ \ \ \ \ $\vdots$} \\
	   \texttt{\ \ \ \ \ \ \ \ \ \ \ \ \ 't\_N': Delta\_t\_N\}}
	  };
	  \draw[very thick, ->] ($(code.south)!0.2!(center)$) -- ($(code.south)!0.7!(center)$);
  \end{scope}
  \begin{scope}[xshift=6cm]
	  \node[circle, fill=black, inner sep=0, minimum size=1mm, label=0] (t0) at (0, 0) {};
	  \node[time point=$t_1$] (t1) at (1, 0)  {};
	  \node[time point=$t_2$] (t2)  at (2, 0) {};
	  \node[label=$\cdots$] (dots)  at (2.7, 0) {};
	  \node[time point=$t_N$] (tN)  at (4, 0) {};

	  \draw (t0) -- ($(t2)!0.5!(dots)$) node (ld) {};
	  \draw[dash pattern={on 1.2mm off 1.2mm}] (ld) -- ($(dots)!0.5!(tN)$) node (rd) {};
	  \draw (rd) -- (tN);
	  \draw[decoration={brace, raise=2pt, mirror, amplitude=2mm}, decorate] (t0) --  node[xshift=-1mm, yshift=-5mm] {$\Delta_{t_1}\!\!\!=\!\!\frac{T}{N}$} (t1);
	  \draw[decoration={brace, raise=2pt, mirror, amplitude=2mm}, decorate] (t1) -- node[xshift=-1mm, yshift=-5mm] {$\Delta_{t_2}\!\!\!=\!\!\frac{T}{N}$} (t2);
	  \draw[decoration={brace, raise=2pt, mirror, amplitude=2mm}, decorate] (t0) +(0, -0.7) -- node[yshift=-1mm, label=below:{$T$}] {} ($(tN) +(0, -0.7)$);

	  \node[yshift=-5mm]  at (ld) {$\cdots$};
	  \coordinate (center) at ($(t0)!0.5!(tN)$);
	  \node[rectangle, draw, fill=\bc, yshift=3cm, align=left] (code) at (center)
	  {\texttt{timesteps = (['t\_1', 't\_2',} \\
	   \texttt{\ \ \ \ \ \ \ \ \ \ \ \ \ $\cdots$, 't\_N'], T)} \\
	  };
	  \draw[very thick, ->] ($(code.south)!0.2!(center)$) -- ($(code.south)!0.7!(center)$);
  \end{scope}
\end{tikzpicture}

%% file: sections/4_case_study.tex

\def\BARONVERSION{\hyperlink{cite.baron20_10_16}{BARON 20.10.16}}
\def\GUROBIVERSION{\hyperlink{cite.gurobi9_1_1}{Gurobi 9.1.1}}

\section{Case Studies} \label{sec:case_study}
We now demonstrate key features of COMANDO in four case studies, which are illustrative of the kinds of design and operation problems we address with COMANDO.
The case studies focus on different aspects of energy systems and vary in their approaches for modeling the considered systems and their components.
The complete source code for all case studies can be found in the \texttt{examples} directory of the \citet{COMANDO_REPO}.

The first case study, based on our previous work \citep{voll2013automated, sass2019optimal}, consists of the greenfield design and operation of an industrial energy system considering both economic and environmental impact.
The component models account for nonlinearities in part-load behavior and investment cost, and differential equations for the state of charge of battery and thermal energy storage units, resulting in a mixed-integer dynamic optimization (MIDO) problem.
Here, the automatic implicit Euler discretization as well as the automatic linearization implemented in COMANDO are employed to obtain a MILP formulation, and a simple user-defined algorithm for multi-objective optimization is demonstrated.

In the second case study, the operation of a simple building energy system is optimized, considering forecasts for electricity price and ambient temperature.
The system model makes use of differential equations to describe the thermal behavior of the building, allowing to represent dynamic aspects of demand response via a MIDO problem.
The interface to Pyomo.DAE \citep{nicholson2018pyomo.dae} is used to apply orthogonal collocation on finite elements as an advanced time discretization method.

The third case study is a variation of the benchmark problem from \citep{saelens2020towards}, integrating low-temperature waste heat into a district heating network via heat pumps.
The explicit consideration of thermal losses and temperatures at different points of the network results in a nonconvex mixed-integer quadratically-constrained quadratic programming (MIQCQP) problem.
For the implementation in COMANDO, repeated structures within the system are abstracted via subsystems, allowing for re-use of the models and reducing modeling effort.
A stochastic formulation considering multiple operational scenarios based on clustered historical data is solved to obtain an optimal system design.

The fourth case study is a reimplementation of our previous work \citep{huster2019impact}, where the power production of an organic Rankine cycle is maximized.
The detailed thermodynamic behavior of the working fluid is described via artificial neural networks (ANNs), capable of predicting fluid properties with high accuracy.
The ANNs result in a highly nonconvex, but reduced-space NLP formulation that can be solved to global optimality with our inhouse solver MAiNGO \citep{bongartz2018maingo}.

All case studies are solved on a desktop PC with an i7-8700 CPU (3.20GHz), 32\,GB RAM, running Windows 10 Enterprise LTSC.
An overview of the presented case studies and their key characteristics is given in \cref{tab:case_studies}.

\subsection{Case study 1: Greenfield design of an industrial energy system} \label{sec:design_problem}
  This case study is inspired by our previous work \citep{sass2020model}.
  For demonstration, we consider a simpler system, allowing only up to one component of each type.
  We make use of inheritance to abstract common model aspects of conversion and storage components into generic classes and then derive specialized variants that implement more specific behavior.
  Furthermore, we take advantage of automatic linearization and discretization routines to obtain MILP problems from the originally dynamic and nonlinear component models of \citet{sass2020model}.

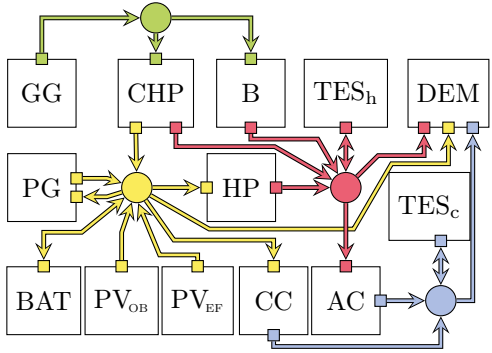
\begin{figure}[tb]
    \centering
    \graphicspath{{figures/}}
    \input{figures/IES_superstructure}
    \caption{
      Superstructure for the industrial energy system case study:
      gas-grid (GG), power-grid (PG), boiler (B), combined heat-and-power unit (CHP), compression chiller (CC), absorption chiller (AC), heat pump (HP), photovoltaic units on office buildings (PV\textsubscript{OB}) and on experimental facilities (PV\textsubscript{EF}), thermal energy storage for hot water (TES$_\text{h}$) and cooling water (TES$_\text{c}$), a battery (BAT), and a demand (DEM). Natural gas is shown in green, electricity in yellow, hot water in red, and cooling water in blue.
    }
    \label{fig:IESstructure}
  \end{figure}

  The industrial energy system needs to satisfy given time-dependent demands for heating, cooling, and electricity with minimal total annualized costs (TAC) and global warming impact (GWI).
  To satisfy these demands, multiple conversion and storage components are available in the superstructure of the system (\cref{fig:IESstructure}).
  For self-containment, we briefly repeat the description of the conversion and storage components here.
  More detailed information can be found in \citet{sass2020model} and in the source code for this case study, available in the \citet{COMANDO_REPO}.

  The conversion components $i\in\mathcal{I}^\text{conv} \{\text{AC, B, CC, CHP, HP}\}$ (\cf \cref{fig:IESstructure}) are modeled with nonlinear investment cost and part-load efficiency curves.
  Additionally, minimal part-load requirements are considered by introducing binary variables.
  The investment cost reflect decreasing marginal investment costs $C^\text{I}_i$ with increasing nominal component output $\dot{E}^\text{nom}_i$, i.e.,
  \begin{equation}
    C^\text{I}_i = C^\text{ref}_i \dot{E}^{\text{nom}^{M_i}}_i \quad\forall i\in\mathcal{I}^\text{conv},
    \label{eq:invIES}
  \end{equation} \noindent
  where $C^\text{ref}_{i}$ and $M_i$ are technology-specific parameters.
  The part-load efficiency $\eta_i$ is expressed via a base efficiency multiplied with a rational function of the part-load fraction ${\dot{E}_i^\text{out}}\slash{\dot{E}_i^{\text{nom}}}$, and describes the relationship of input $\dot{E}_i^\text{in}$ and output $\dot{E}_i^\text{out}$:
  \begin{align}
      \dot{E}_i^\text{out} = \eta_i \dot{E}_i^\text{in} \quad&\forall i\in\mathcal{I}^\text{conv}
  \label{eqn:inputoutputIES}
  \end{align}
  The HP and CHP models have variable base efficiencies that depend on temperatures and the nominal size, respectively.
  We create a generic conversion component class with an un\-pa\-ra\-me\-trized nonlinear efficiency and investment cost model (\cref{eq:invIES,eqn:inputoutputIES}).
  From this conversion component class, we derive the individual conversion technologies as subclasses.
  Three instances of the CHP model with different ranges for the nominal size are considered, accounting for the size-dependence of the conversion efficiencies for heat and electricity.
  The three CHP models are aggregated into a subsystem which enforces that at most one of them is built.
  The subsystem can then be incorporated into other system models like any other component.

  The storage components $i \in \mathcal{I}^\text{sto} = \{\text{BAT, TES}_\text{h}, \text{TES}_\text{c}\}$ are modeled with the differential equation
  \begin{equation}
    \label{eqn:storage}
    \der{E_i} = \eta_i^\text{in} \dot{E}_i^\text{in} - \frac{1}{\eta_i^\text{out}}\dot{E}_i^\text{out} - \frac{1}{\tau_i} E_i \quad \forall i\in\mathcal{I}^\text{sto},
  \end{equation} \noindent
  where the state $E_i$ is the stored energy, $\eta_i^\text{in}$ and $\eta_i^\text{out}$ are constant charging and discharging efficiencies, $\dot{E}_i^\text{in}$ and $\dot{E}_i^\text{out}$ are the charging and discharging rates, and $\tau_i$ is a time constant describing self-discharging.
  As with the conversion components, we create a generic storage component class and derive technology-specific sub-classes, e.g., batteries.
  For each component we additionally consider a binary variable and associated constraints, representing whether the component is built or not.

  We use the aggregated data from the supplementary material of \citet{sass2020model}, which originate from clustering a full year of data for demands, weather, prices, and global warming impacts via the method described in \cite{bahl2018typical}.
  The aggregated data represent the full year via four typical days, each with four time steps of varying lengths (between 1 and 17 hours), and two isolated time points of length zero, representing peak heating and cooling demands.
  For a similar design problem, \cite{bahl2018typical} showed that even coarse time resolutions such as this one provide optimal objective values, sufficiently close to those obtained with a full year at hourly resolution.
  In COMANDO we can consider such a time structure via six scenarios, corresponding to the four typical days and the two isolated time points for peak demands.
  The scenarios corresponding to typical days are weighted by number of days associated to them during clustering, and the scenarios for peak demands are assigned a weight of zero, \ie, they have no effect on the objective but are considered for feasibility, \cf formulation \eqref{prob:P}.

  Due to the storage dynamics \cref{eqn:storage}, problems derived from this system model will be MIDO problems.
  In our previous work \citep{sass2020model}, we manually implemented the MILP formulation resulting from explicit Euler discretization and a case-specific linearization in GAMS.
  As this process and subsequent changes are labor-intensive and error-prone, we instead make use of COMANDO's automatic routines for discretization and piecewise linearization.

  \setpgfexternalcounter{7}
  \begin{figure}[t!]
    \centering
    \vspace*{-4.5mm}
    \includegraphics{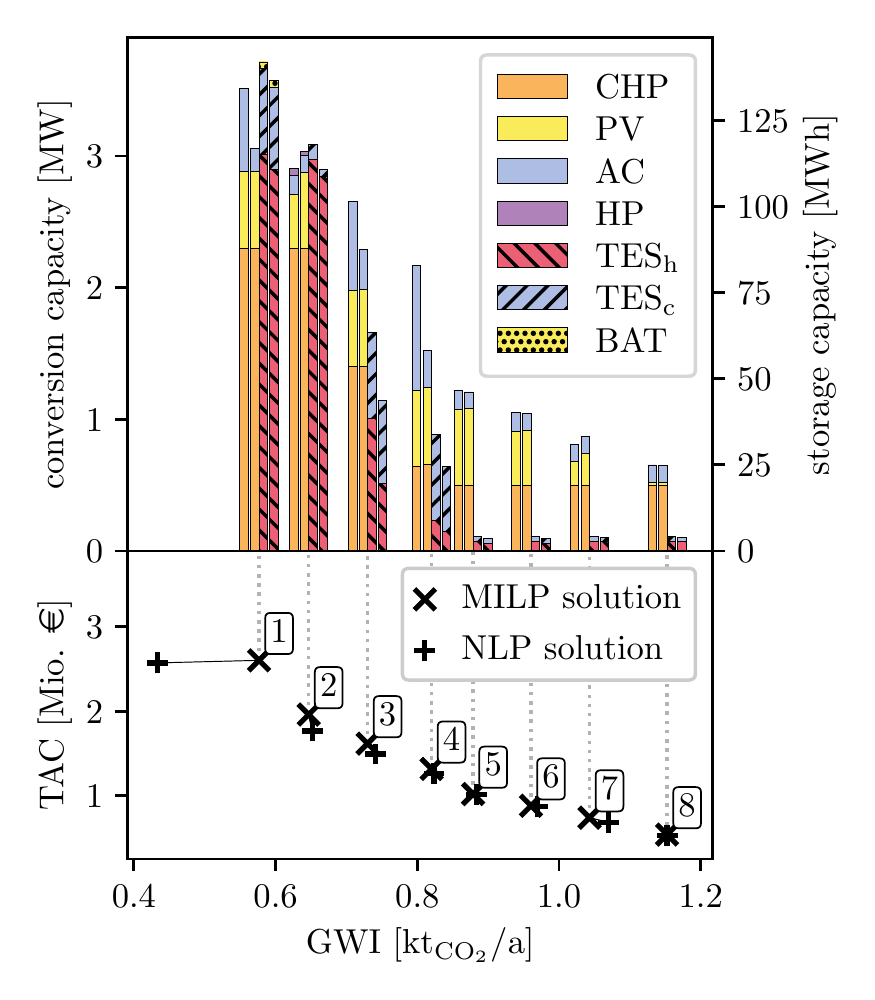}
    \vspace*{-5mm}
    \caption{
      Bottom: eight Pareto-optimal designs, determined from multi-objective optimization regarding total annualized cost (TAC) and global warming impact (GWI). Top: corresponding capacities of conversion (left) and storage components (right) from the MILP (inner bars) and NLP (outer bars) formulations.
      CHP: Combined heat and power unit, PV: photovoltaic array, AC: absorption chiller, HP: heat pump, TES$_\text{h}$: hot thermal energy storage, TES$_\text{c}$: cold thermal energy storage, BAT: battery.
      Note that boilers and     compression chillers are not part of any design and thus excluded from the legend.
    }
    \label{fig:paretofrontier}
  \end{figure}

  The augmented $\varepsilon$-constraint method \citep{mavrotas2009effective} is implemented as a user-defined algorithm, in which two design optimization problems with either TAC or GWI as objective function are repeatedly solved.
  For the solution of the two problems we use \GUROBIVERSION{} with a relative optimality tolerance of 1\%.
  Generating 8 designs from the Pareto front for TAC and GWI, shown in \cref{fig:paretofrontier}, takes about 3.6 hours.
  Note that a Pareto-optimal design can only improve upon one of the two objectives by worsening the other.

  The total GWI can be reduced by 50\% (from 1.152 to 0.577 kt/a) when accepting a fourfold increase in TAC (from 0.539 to 2.6~Mio. \euro) (\cref{fig:paretofrontier}, bottom).
  Solutions with lower TAC are characterized by small component capacities with lower investment costs, whereas solutions with lower GWI rely on large conversion and storage components (\cref{fig:paretofrontier} inner bars, top).
  As these results were obtained with a linearization of the original model, they are only approximate and the corrsponding designs may not be feasible with respect to the nonlinear model.

  However, correcting the infeasibilities is straightforward in COMANDO as the original, nonlinear model formulation is available.
  We first obtain the MINLP problem resulting from implicit Euler discretization of the original formulation with TAC as the objective.
  We then repeat the multi-objective optimization with the same algorithm but using the MINLP formulation.
  For each iteration, we set the appropriate upper bound on GWI and fix binary variables to the values of the corresponding MILP solution, obtaining an NLP formulation.
  The values of the remaining variables are used as an initial point and the resulting formulation is passed to \BARONVERSION{} using default options, except for a relative optimality tolerance of 1\% and a time limit of one hour for the subproblems.

  In three cases the subproblems are terminated due to the time limit (with 3.5\% relative gap for the TAC minimization of iteration 3 and 4, and 7.5\% relative gap for the GWI correction of iteration 3).
  The remaining subproblems take at most 78 s to be solved to the desired optimality.
  Thus, all cases result in a design and an operational strategy that are feasible with respect to the original nonlinear formulation.
  The resulting solutions exhibit slightly lower TAC values and slightly higher GWI values than their MILP counterparts, with the exception of iteration 1, where the GWI value is 25\% lower than for the MILP approach (433 t/a vs. 577 t/a).
  The corresponding designs can be seen in the outer bars in \cref{fig:paretofrontier} (top).
  While the MILP and NLP solutions of iterations 2 and 5--8 are similar, iterations 1, 3 and 4 exhibit larger conversion components and smaller storages in the NLP case.
  In summary, the approach provides MINLP-feasible system designs that allow a trade-off between the TAC and GWI of the resulting system.

\subsection{Case study 2: Demand response of a building energy system} \label{sec:dynamic_operation}
  To illustrate how to formulate and solve optimization problems with more pronounced dynamic effects in COMANDO, we model an illustrative building energy system.
  The system is heated by a heat pump (HP) and is capable to perform load shifting via concrete core activation, i.e., a concrete core with a high thermal inertia can be heated directly.
  We investigate a demand response (DR) case, where we optimize the operation of the building energy system over the horizon of one day with given profiles for electricity price and ambient temperature.

  The considered building energy system consists of three thermal zones: air, outside wall, and concrete core.
  Occupant comfort has to be ensured by maintaining the air temperature between minimal and maximal temperatures $T_\text{air}^\text{min}$ and $T_\text{air}^\text{max}$, respectively.
  To do so, the air in the room can be heated via a direct heat flow to the air $\dot Q_\text{air,in}$, or indirectly through the concrete core, which can be heated via the heat flow $\dot Q_\text{core,in}$.
  We consider a zero-dimensional model of each thermal zone. For instance, the energy balance of the air zone is given by
  \begin{align}
    \rho_\text{air}V_\text{air}c_{p,\text{air}}\frac{\text{d}T_\text{air}}{\text{d}t} = \dot Q_\text{core,air} - \dot Q_\text{air,wall} + \dot Q_\text{air,in},
    \label{eqn:air_mass}
  \end{align}
  where $T_\text{air}, V_\text{air}, \rho_\text{air}$, and $c_{p,\text{air}}$ are the air temperature, volume, density, and specific heat capacity, respectively, and $\dot Q_\text{core,air}$ and $\dot Q_\text{air,wall}$ are heat exchange flows with the adjacent zones.
  The heat flow $\dot Q_{A,B}$ between two zones $A$ and $B$ is calculated depending on the temperatures $T_A$ and $T_B$, the area $A_{A,B}$, and the heat transfer coefficient $U_{A,B}$:
  \begin{align}
    \dot Q_{A,B} = U_{A,B}A_{A,B}(T_A - T_B) \label{eqn:heat_transfer}
  \end{align}

  \begin{figure}[!t]
    \centering
    \input{figures/BES}
    \caption{
      Structure of the considered building energy system and implementation in COMANDO: three instances of the thermal mass class (M$_\text{air}$, M$_\text{core}$, M$_\text{wall}$), four instances of the heat transfer class (HT$_\text{air,wall}$, HT$_\text{wall,E}$, HT$_\text{air,core}$, HT$_\text{core,wall}$), heat pump (HP), and power grid (PG). Red arrows represent heat flows and yellow arrows electric power flows.
    }
    \label{fig:BES}
    \vspace*{-2mm}
  \end{figure}
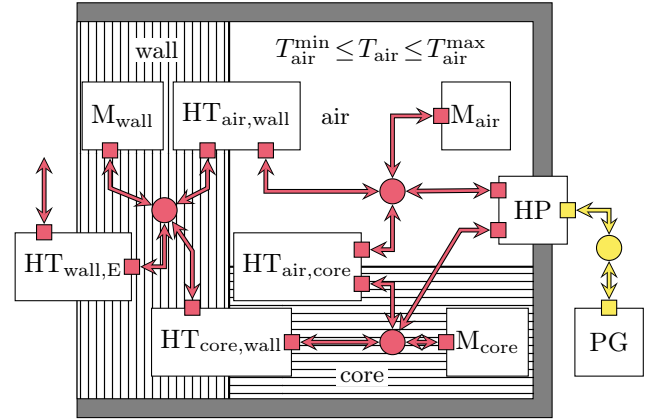

  The structure of the model is shown in \cref{fig:BES}.
  To model thermal masses, we introduce a component M, which is instantiated by specifying volume, density, and specific heat capacity, and optionally allows to specify minimal and maximal temperatures.
  The heat transfer is abstracted as a component HT, implementing \cref{eqn:heat_transfer}, and the heat pump is again modeled with a temperature-dependent efficiency, but with the option of splitting the output to multiple connectors.

 \begin{figure}[!t]
  \centering
  \input{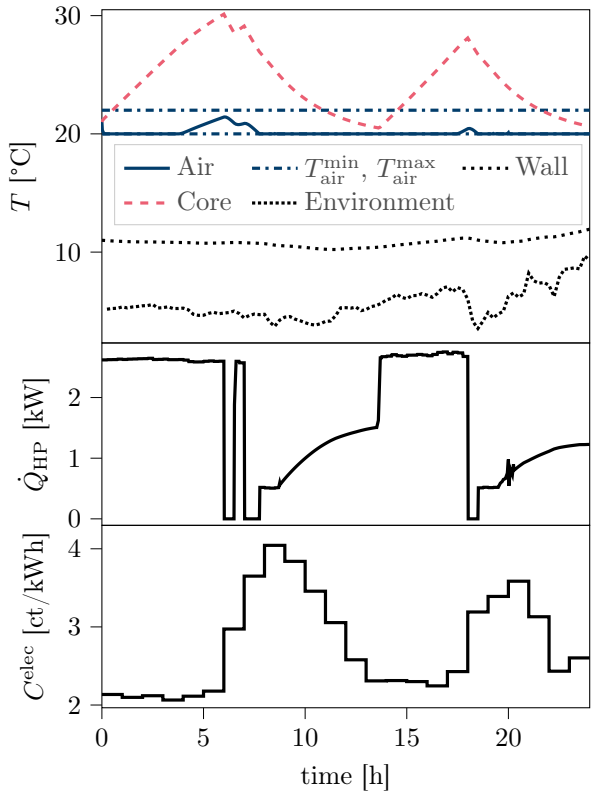}
  \vspace*{-3mm}
  \caption{
    Results of the demand response optimization for the building energy system: the temperatures of the different zones together with the air temperature comfort bounds $T_\text{air}^\text{min}$ and $T_\text{air}^\text{max}$ (top), the heat flow supplied by the heat pump $\dot{Q}_\text{HP}$ (center), and the electricity costs $C^\text{elec}$ (bottom).
  }
  \label{fig:BES_DSM}
  \vspace*{-3mm}
  \end{figure}

  Based on the model of the building energy system, we define a DR optimization problem, i.e., we minimize the integral over the electricity costs for a given electricity price profile.
  The resulting operational objective function is thus chosen as $\dot{F}_\operation = C^\text{elec} P_\text{HP}$, where $C^\text{elec}$ and $P_\text{HP}$ are the electricity costs and electric input power of the heat pump, respectively.

  As we consider a minimum part-load constraint for the heat pump, the resulting problem is a MIDO problem.
  The time horizon is a 24 hour period considered at quarter-hourly resolution and the input data consists of hourly electricity price data and ambient temperature data at quarter-hourly resolution.
  We use a full discretization approach \citep{cuthrell1987optimization} via the COMANDO interface to Pyomo.DAE \citep{nicholson2018pyomo.dae}.
  Specifically, we use Legendre-Radau collocation with four elements per hour and fourth-order polynomials.
  Since the model contains exclusively linear expressions and we use collocation with a fixed time grid, we obtain a MILP problem after discretization.
  The resulting formulation has 6931 constraints and 6257 variables, 96 of which are binary.
  The problem can be solved with \GUROBIVERSION{} to global optimality in less than one second of CPU time.
  Results are visualized in \cref{fig:BES_DSM}, where the temperatures of the three thermal zones, the ambient temperature, the heat flow supplied by the heat pump, and the electricity price are shown.
  During times of low prices, the concrete core is heated to store energy.
  During times of high prices, the concrete core transfers the stored heat to the air zone and cools down such that the heat pump has to supply less heat.
  Thus, load is shifted to times of favorable prices, while the air temperature remains within the comfort range.

  Using the introduced component models for general thermal masses and heat transfers, the extension to a larger building energy system with several rooms, thermal masses, and heat transfers is straightforward.
  We note that it is also possible to perform rolling horizon optimization in COMANDO by defining an appropriate user-defined algorithm, e.g., as in our previous publication \citep{shu2019optimal}, where a preliminary version of COMANDO was used.

\subsection{Case Study 3: Design of a low-temperature district heating network} \label{sec:HP_net}
    In this case study, we extend components of previous work \citep{hering2020design} to describe a district heating network and apply them to a design optimization of the network described by \cite{saelens2020towards}.
    The system comprises a source of waste heat, a distribution network, and 16 consumers.
    We aggregate the 16 consumers into four consumer groups, comprising four consumers each, and assume linear heating curves for the flow temperature $T^\text{fl}$.
    The heating curves are described by the flow temperatures $T^\text{fl,max}$ and $T^\text{fl,min}$, at ambient air temperatures of -$12\si{\celsius}$ and $20\si{\celsius}$, respectively, see \cref{tab:building_groups}.

    \begin{table}[tb]
      \centering
      \caption{
        Clustering of neighbouring buildings into consumer groups. Buildings within a group are assumed to have identical heating curves.
      }
      \begin{tabular}{ccc}
        \toprule
         Consumer & $T^\text{fl,max}$
                        & $T^\text{fl,min}$ \\
         group    & $(T_\text{air}\!=\!-12\si{\celsius})$
                        & $ (T_\text{air}\!=\!20\si{\celsius})$ \\
         \midrule
         CG$_{40}$ & $40\si{\celsius}$ & $35\si{\celsius}$ \\
         CG$_{50}$ & $50\si{\celsius}$ & $40\si{\celsius}$ \\
         CG$_{70}$ & $70\si{\celsius}$ & $50\si{\celsius}$ \\
         CG$_{85}$ & $85\si{\celsius}$ & $60\si{\celsius}$ \\
         \bottomrule
      \end{tabular}
      \label{tab:building_groups}
    \end{table}

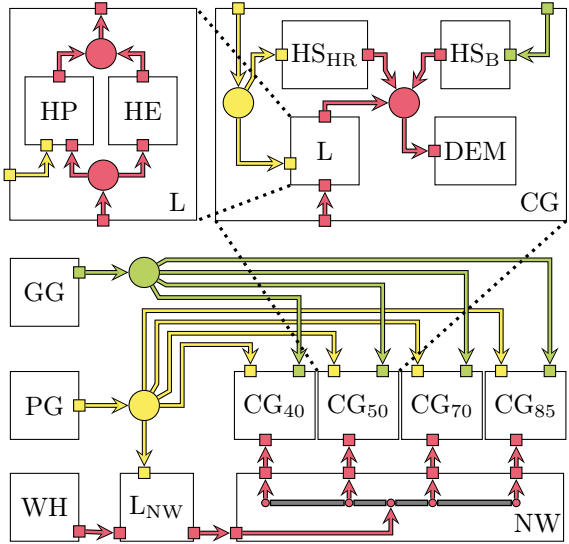
\begin{figure}[t]
      \centering
      \input{figures/DESTEST}
      \caption{
        Superstructure with components for the gas grid (GG), power grid (PG), waste heat source (WH) and network (NW) as well as subsystems for linking (L) and consumer groups (CG), see top.
        The superstructure of the linking subsystem contains a heat pump (HP) and a heat exchanger (HE) and that of the consumer group subsystems contains two heat source (HS) components parameterized as a heating rod (HS$_\text{HR}$) and a boiler (HS$_\text{B}$) and demand (DEM) as well as a decentral linking subsystem.
        To connect the different consumer groups, the necessary pipe segments (depicted as gray bars within NW) need to be built.
      }
      \label{fig:ltdh}
    \end{figure}

    The source of waste heat supplies heat to a network to which each consumer group may be connected or not.
    Both waste heat and consumer groups are linked to the network via a heat exchanger or a heat pump, and connecting a consumer group additionally requires the necessary pipes to be built.
    Independently of whether a consumer group is connected or not, it may also be equipped with a gas-fired boiler or an electric heating rod.
    The superstructure of the heating network is shown in \cref{fig:ltdh}.

    The system is modeled using components for a source of waste heat (WH), the distribution network (NW), the power grid (PG) and the gas grid (GG).
    As both the linking unit and the consumer groups are composed of multiple components and occur more than once, they are modeled as subsystems.
    The linking subsystem (L) contains a heat pump (HP) and a heat exchanger (HE) and the consumer group subsystem (CG) contains a linking subsystem, a demand (DEM), and two instances of a generic heat source with different parametrizations, representing a boiler (HS$_\text{B}$) and a heating rod (HS$_\text{HR}$).

    The design decisions comprise binary variables for the type of linking component (heat exchanger, heat pump, or none) and the type of additional heat source (gas boiler, electric heater, or none) to be built, as well as continuous variables for component sizing and the maximum and minimum return temperature of the network $T^\text{re,max}_\text{NW}$ and $T^\text{re,min}_\text{NW}$, respectively.
    Finally, four pipe segments can be added to the network model separately using the decision variables, $b_\text{NW}^\text{s}$.
    The linking components for the consumer groups can only be built if all necessary pipe segments of the network are built.
    The demand component has a parameter for the required heat demand and computes the required flow temperature based on the ambient air temperature.
    The heat demand is based on \citet{saelens2020towards}, while the flow temperature is assumed to depend linearly on the ambient air temperature (\cf \cref{tab:building_groups}).
    The network return temperature $T^\text{re}_\text{NW}$ also depends linearly on the ambient air temperature $T_\text{air}$ and the design variables $T^\text{re,max}_\text{NW}$ and $T^\text{re,min}_\text{NW}$, while the network flow temperature $T^\text{fl}_\text{NW}$ is assumed to be 15\,K higher than $T^\text{re}_\text{NW}$.

    We aggregate the whole network into one pipe network with two branches, \cf \cref{fig:ltdh}.
    The central linking component is connected to the center of the network with $T^\text{fl}_\text{NW}$ and $T^\text{re}_\text{NW}$.
    Despite being located at different distances from the center, we assume that all consumer groups receive and reject water at the same flow and return temperatures, $T^\text{fl}_\text{NW} - \Delta T^\text{fl,loss}_\text{NW}$ and $T^\text{re}_\text{NW} + \Delta T^\text{re,loss}_\text{NW}$, respectively.
    For this simplification to be conservative, we use the total length of the network, $l_\text{NW}$, calculated as
    \begin{equation}
      \label{eqn:network_length}
      l_\text{NW} = \sum b_\text{s} l_\text{s},
    \end{equation}
    to calculate the temperature drops, where $b_\text{s}$ is the build decision and $l_\text{s}$ is the length of each network segment $s$, \cf \cref{fig:ltdh}.
    To obtain the temperature differences in the flow and return pipes, $\Delta T^\text{fl,loss}_\text{NW}$ and $\Delta T^\text{re,loss}_\text{NW}$, respectively, we consider energy balances of the water for both the flow (fl) and return (re) pipe of the network, i.e.,
    \begin{align}
        \dot{m}_\text{NW} \, c_p \, \Delta T^\text{fl,loss}_\text{NW} &= U_\text{NW} \, l_\text{NW} \, (T^\text{fl}_\text{NW} - T_\text{gr})
        \label{eq:thloss3} \\
        \dot{m}_\text{NW} \, c_p \, \Delta T^\text{re,loss}_\text{NW} &= U_\text{NW} \, l_\text{NW} \, (T^\text{re}_\text{NW} + \Delta T^\text{re,loss}_\text{NW}  - T_\text{gr})
    \label{eq:thloss4}
    \end{align}
    where $\Delta T^\text{fl,loss}_\text{NW}$ and $\Delta T^\text{re,loss}_\text{NW}$ are operational variables describing the temperature drop in the respective pipe, $c_p$ is the constant specific heat capacity of water, $U_\text{NW} = 0.035 \frac{\text{W}}{\text{mK}}$ is the specific heat transfer coefficient and $l_\text{NW}$ is the pipe network length, and $T_\text{gr} = 8\si{\celsius}$ is the average ground temperature.

    \def\QHE{\dot{Q}_{\text{HE}}}
    \def\QHP{\dot{Q}_{\text{HP}}}
    \def\meva{\dot{m}_{\text{eva}}}
    \def\mcon{\dot{m}_{\text{con}}}
    \def\Tevafl{T^\text{fl}_{\text{eva}}}
    \def\Tevare{T^\text{re}_{\text{eva}}}
    \def\Tconfl{T^\text{fl}_{\text{con}}}
    \def\Tconre{T^\text{re}_{\text{con}}}

    The heat pump model in each linking component is modeled via the following set of equations:
    \begin{align}
      &\QHP \leq b_\text{HP} \, 400\,\text{kW}   \label{eq:QHP} \\
      &P_\text{HP} \, \Tconre \, \eta_\text{COP} = \QHP (\Tconre - \Tevare) \label{eq:carnot} \\
      &\meva c_p \, (\Tevafl - \Tevare) + P_\text{HP} \nonumber \\
      & \quad = \mcon c_p \, (\Tconre - \Tconfl) \label{eq:energy_balance}
    \end{align}

    Here,
    $\QHP$,
    $P_\text{HP}$,
    $\meva$,
    $\mcon$,
    $\Tevafl$,
    $\Tevafl$,
    $\Tconre$
    and
    $\Tconfl$
    are operational variables, and $\eta_\text{COP} = 0.6$ is the heat pump efficiency relative to the carnot efficiency.
    The outgoing heat flow for the heat pump ($\QHP$) is bounded by zero or the maximum allowable nominal size of 400 kW through \cref{eq:QHP}.
    The input power $P_\text{HP}$ is coupled to $\QHP$ via \cref{eq:carnot}.
    In the energy balance \cref{eq:energy_balance}, enthalpy differences at the evaporator and condenser side are described by the associated mass flows
    $\meva$
    and
    $\mcon$,
    and flow and return temperatures
    $\Tevafl$,
    $\Tevafl$,
    $\Tconfl$
    and
    $\Tconre$.

    For the investment cost, we assume linear cost correlations with a specific cost c$_\text{spec}$ and a fixed cost c$_\text{fix}$ according to \cref{tab:my_label}.

    \begin{table}[t]
      \centering
      \defcitealias{bundesministeriumverkehr2012ermittlung}{BMVBS, 2012}
      \defcitealias{bundesinstitutbau2014kosten}{BBSR, 2014}
      \caption{Specific and fixed costs for heating equipment according to \citepalias{bundesministeriumverkehr2012ermittlung} and \citepalias{bundesinstitutbau2014kosten}}
      \label{tab:my_label}
      \begin{tabular}{lrr}
      \toprule
        Component & c$_\text{spec}$ & c$_\text{fix}$ \\
        \midrule
        central HP & \SI{500}{\euro\per\kilo\watt} & \SI{0}{\euro} \\
        decentral HP & \SI{620}{\euro\per\kilo\watt} & \SI{0}{\euro} \\
        HE & \SI{90}{\euro\per\kilo\watt} & \SI{0}{\euro} \\
        HS$_\text{HR}$ & \SI{10}{\euro\per\kilo\watt} & \SI{100}{\euro} \\
        HS$_\text{B}$ & \SI{111}{\euro\per\kilo\watt} & \SI{4300}{\euro} \\
        \bottomrule
      \end{tabular}
    \end{table}

    Additionally, we consider the costs for each pipe segment of the network based on \cite{jentsch2008handbuch}.
    Thus, the total investment costs of the system includes the investments into heating components and piping.

\begin{figure}[!t]
      \centering
      \input{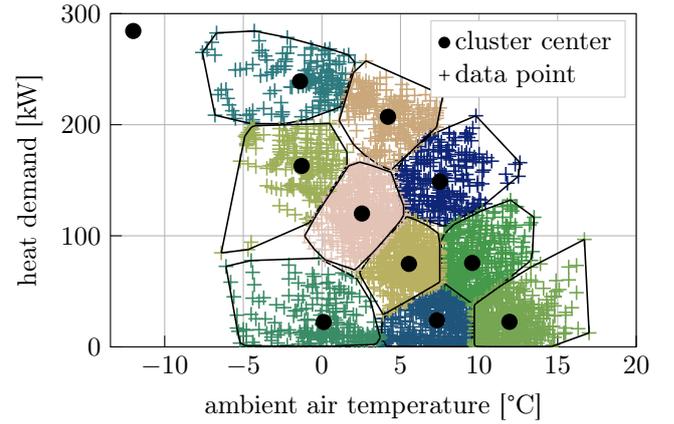}
      \vspace*{-2mm}
      \caption{
        Heat demand clusters: Each cross represents one pair of measurements of total daily mean heat demand and daily mean ambient air temperature. Colors and boundaries are used to aid visual distinction of the clusters whose centers are mean values depicted as black dots.
      }
      \label{fig:heatcluster}
    \end{figure}

    To obtain an economical design, we minimize TAC.
    We use k-means clustering \citep{pedregosa2011scikit} to aggregate the original set of ambient temperatures and heat demands into representative clusters.
    Each resulting cluster center is a pair of daily mean values for temperature and heat demand and can be considered as a representative operating scenario.
    To reduce computational demand, the data is clustered into 11 such scenarios, including one scenario representing the maximum heat demand.
    Demand data with zero heat demand are dropped from the dataset.
    \cref{fig:heatcluster} shows the resulting 11 clusters.

\begin{figure}[!bt]
      \centering
      \input{figures/ltdh_result}
      \caption{
        Optimal system structure: A central heat pump HP supplies waste heat from WH to the network NW.
        Consumer groups CG$_{40}$ and CG$_{50}$ are connected to NW via heat exchangers (HE) and use heating rods (HR) for peak demands.
        Consumer groups CG$_{70}$ and CG$_{85}$ are not connected and satisfy their heat demand via boilers (B).
      }
      \label{fig:ltdh_result}
    \end{figure}
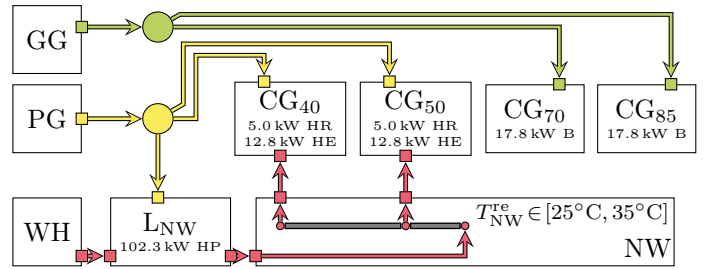

    We use the clusters as scenarios in the COMANDO framework, with the fraction of data points in each cluster as the corresponding scenario weight.
    Considering the data in this way ensures that the final design is feasible for all considered scenarios and is optimized with regards to the expected value of TAC.
    The resulting problem is a MIQCQP with 526 continuous variables, 147 binary variables and 275 quadratic constraints.
    An optimal design with 0\% optimality gap is obtained within six minutes of CPU time, using the Gurobi API interface with \GUROBIVERSION{}
    and 12 threads.
    The global optimal solution corresponds to the system shown in \cref{fig:ltdh_result}.

    The network is designed with a variable return temperature between $25\si{\celsius}$ and $35\si{\celsius}$ and is connected to the waste heat source using a \SI{102}{\kilo\watt} heat pump.
    Consumer groups CG$_{40}$ and CG$_{50}$ are connected to the network using heat exchangers and have additional electric heating rods installed.
    Consumer groups CG$_{70}$ and CG$_{85}$ are not connected but satisfy their heat demand using gas-fired boilers instead.
    The TAC of this design are \SI{22095}{\euro}.
    At an annual heat demand of 322.7\,MWh this corresponds to a specific heating cost of \SI{68.5}{\euro\per\mega\watt\hour}.
    In order to assure that the obtained system design is feasible for the original demand data, we perform a second optimization for which we fix the design (\ie, system structure and component sizes) and perform a purely operational optimization using the full set of demands.
    The design proves to be feasible, with the corrected TAC increasing by less than 2\% to \SI{22414}{\euro}.

\subsection{Case study 4: Optimal operating point of an organic Rankine cycle (ORC)} \label{sec:ORC}
    Finally, we consider a case study from our previous work \citep{huster2019impact}, where an optimal operating point of an organic Rankine cycle (ORC) with respect to net power production is sought.
    With this case study, we demonstrate how COMANDO can handle complex modeling features such as accurate fluid properties via artificial neural networks \citep{schweidtmann2018deterministic, schweidtmann2019deterministic} and a sequential modeling approach that gives rise to reduced-space formulations beneficial for global optimization \citep{bongartz2017deterministic}.

    Again, we give a short overview of the case study for self-containment.
    In the considered process, the working fluid isobutane (ib) is first pressurized by a pump and then preheated in a recuperator before being
    heated to evaporation temperature, evaporated and superheated by cooling geothermal brine (gb) from 408\,K to 357\,K.
    After expanding in a turbine, the working fluid is used in the recuperator to preheat the pressurized fluid and is finally condensed and cooled to its original state using cooling water at 288\,K.
    The heat passed from the condenser to the cooling water (cw) is dissipated by a cooling system consisting of multiple fans.

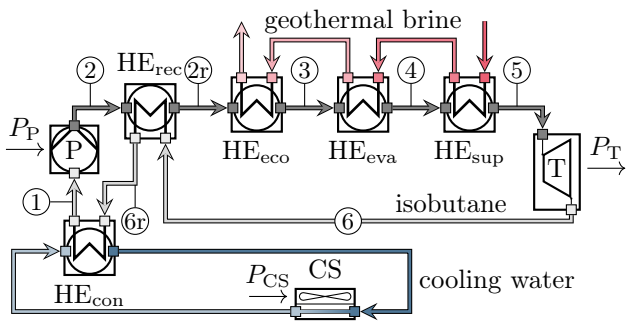
\begin{figure}[!b]
      \centering
      \vspace*{-1mm}
      \input{figures/ORC}
      \caption{
        System model of the ORC process from \citet{huster2019impact}.
        The components are a pump (P), a recuperator (HE$_\text{rec}$), an economizer (HE$_\text{eco}$), an evaporator (HE$_\text{eva}$), a superheater (HE$_\text{sup}$), a turbine (T), a condenser (HE$_\text{con}$), and a cooling system (CS).
        Flows of geothermal brine, the working fluid isobutane, and cooling water are depicted in red, gray, and blue, respectively.
        Electrical power is consumed by pump ($P_\text{P}$) and cooling system ($P_\text{CS}$) and produced by the turbine ($P_\text{T}$).
      }
      \label{fig:ORC}
      \vspace*{-3mm}
      \end{figure}

    The ORC is modeled as a system consisting of 4 types of components, \ie, a pump (P), a turbine (T), a cooling system  (CS), and five heat exchangers (condenser HE$_\text{con}$, recuperator HE$_\text{rec}$, economizer HE$_\text{eco}$, evaporator HE$_\text{eva}$, and superheater HE$_\text{sup}$).
    All components have connectors for enthalpy in- and out-flows that are connected as depicted in \cref{fig:ORC} to obtain the system model.

    As discussed in \citet{bongartz2017deterministic}, reduced-space formulations, \ie, formulations in which a large number of variables and constraints are eliminated by substitution, are well suited for global optimization of power cycles such as the present ORC.
    To obtain a reduced-space formulation, model generation begins with an empty system model to which different component models are added sequentially.
    First, the decision variables are specified at the system level as follows:
    The mass flow $\dot{m}$ of the working fluid, the pressures $p_1$ and $p_2$ before and after the pump, and the specific enthalpy after the recuperator $h_\text{2r}$, as well as the isentropic specific enthalpy after the turbine $h^\text{is}_6$.
    All other quantities of interest are defined in terms of these five variables.

    In our previous work \citep{schweidtmann2018deterministic, schweidtmann2019deterministic}, the use of artificial neural networks (ANNs) in combination with our inhouse global MINLP solver MAiNGO \citep{bongartz2018maingo} has been shown to result in tight relaxations, beneficial for deterministic global optimization.
    In \citet{huster2019impact}, we trained several ANNs to learn the relations between various quantities of different thermodynamic phases of the working fluid isubutane, using data generated from the equations of state implemented in the thermophysical property library CoolProp \citep{bell2014pure}.
    The ANNs are used as data-driven surrogte models for the equations of state, which cannot be used directly within the optimization, as they are not available as analytical expressions \citep{schweidtmann2019deterministic}.
    The validity of the ANNs used for this case study was extensively analyzed and discussed in the original publication \citep{huster2019impact}.
    Each ANN
    has two hidden layers with six neurons each, all of which use tanh as the activation function.
    The ANNs express individual output quantities in terms of either pressure $p$, pressure and specific enthalpy $h$, or pressure and specific entropy $s$, as inputs.
    As a result of training, we thus obtain explicit analytical expressions for various quantities.
    In this case study, eight of the ANNs from \citet{huster2019impact} are used as analytical surrogate models for the following quantities:
    \begin{itemize}[itemindent=2cm, noitemsep]
      \item [$h^\text{liq}(p,s)$] liquid enthalpy
      \item [$T^\text{liq}(p,h)$] liquid temperature
      \item [$h^\text{sat,liq}(p)$] enthalpy of saturated liquid
      \item [$s^\text{sat,liq}(p)$] entropy of saturated liquid
      \item [$T^\text{sat}(p)$] saturation temperature
      \item [$h^\text{sat,vap}(p)$] enthalpy of saturated vapor
      \item [$s^\text{vap}(p,h)$] vapor entropy
      \item [$T^\text{vap}(p,h)$] vapor temperature
    \end{itemize}

    The enthalpy flows of pump and turbine are described via mass flow and specific enthalpies, and the electrical power consumed by the pump ($P_\text{P}$) and provided by the turbine ($P_\text{T}$) are modeled as
    \begin{align}
      P_\text{P} &= \dot{m} \, \frac{h^\text{is,out}_\text{P} - h^\text{in}_\text{P}}{\eta^\text{is}_\text{P}}, \label{eqn:P}\\
      P_\text{T} &= \dot{m} \, (h^\text{in}_\text{T} - h^\text{is,out}_\text{T}) \, \eta^\text{is}_\text{T},
    \end{align}
    where $\eta^\text{is}_\text{P}$ and $\eta^\text{is}_\text{T}$ are known, constant isentropic efficiencies and the required specific enthalpies $h$ are determined via the appropriate ANNs.

    For each heat exchanger, the differences of enthalpy flows at the hot (h) and cold (c) side are either defined in terms of a mass flow and specific enthalpies (ib) or in terms of a specific heat capacity flow $\dot{m}c_p$ and temperatures (cw and gb):
    \begin{align}
      \dot{Q}_\text{h} &= \begin{cases}
        \dot{m}_\text{h} \, (h^\text{in}_\text{h} - h^\text{out}_\text{h}), & \text{h = ib} \\
        (\dot{m}c_p)_\text{h} \, (T^\text{in}_\text{h} - T^\text{out}_\text{h}), & \text{h $\in$ \{cw, gb\}}
      \end{cases} \label{eqn:Qh}\\
      \dot{Q}_\text{c} &= \begin{cases}
        \dot{m}_\text{c} \, (h^\text{out}_\text{c} - h^\text{in}_\text{c}), & \text{c = ib} \\
        (\dot{m}c_p)_\text{c} \, (T^\text{out}_\text{c} - T^\text{in}_\text{c}), & \text{c $\in$ \{cw, gb\}}
      \end{cases} \label{eqn:Qc}
    \end{align}
    As heat losses are neglected, the energy balance reduces to
    \begin{equation}
      \dot{Q}_\text{h} = \dot{Q}_\text{c}. \label{eqn:EB}
    \end{equation}

    Since we aim for a reduced-space formulation, no variables are introduced for the left-hand sides of \crefrange{eqn:P}{eqn:Qc} and the corresponding right-hand side expressions are used directly, avoiding the addition of constraints.
    In particular, where possible, \cref{eqn:EB} is automatically reformulated to obtain a definition for one of the temperatures or specific enthalpies in the right-hand sides of \cref{eqn:Qh,eqn:Qc} in terms of the other quantities.
    The heat-exchanger model is configured to perform the appropriate reformulation automatically, based on the provided quantities.

    A pinch point is assumed in the condenser, \ie, the temperature of the cooling water at the pinch point, $T_\text{pinch}$, is assumed to lie $\Delta T_\text{min}=10\,$K below the evaporation temperature $T^\text{sat}(p_1)$.
    Through this assumption, it is possible to compute the heat capacity flow of the cooling water, $(\dot{m}c_p)_\text{cw}$, as
    \begin{align}
      (\dot{m}c_p)_\text{cw} &= \frac{\dot{m}\,(h_\text{pinch} - h_1)}{\max(10^{-5}\,\text{K},\;T_\text{pinch} - T^\text{in}_\text{cw})} \nonumber \\
        &= \frac{\dot{m}\,\Big(h^\text{sat,vap}(p_1) - h^\text{sat,liq}(p_1)\Big)}{\max\Big(10^{-5}\,\text{K},\;T^\text{sat}(p_1) - 10\,\text{K} - 288\,\text{K}\Big)}. \label{eqn:mcp_cw}
    \end{align}
    Note that the max function and the constant $10^{-5}$ in \cref{eqn:mcp_cw} are introduced to avoid division by zero.
    The electrical power $P_\text{CS}$, required to run the fans of the cooling system, is modeled to be proportional to the specific heat capacity flow of the air $(\dot{m}c_p)_\text{air}$ passing through them and is computed as
    \begin{align}
      P_\text{CS} &= \frac{\dot{V}_\text{air}\,\Delta p_\text{fan}}{\eta_\text{fan}} \nonumber \\
                  &= \frac{(\dot{m}c_p)_\text{air}\,\Delta p_\text{fan}}{c_{p,\text{air}}\,\rho_\text{air}\,\eta_\text{fan}},
    \end{align}
    where $\Delta p_\text{fan} = 170\,$Pa and $\eta_\text{fan} = 0.65$ are the pressure drop and efficiency of the fan, $\dot{V}_\text{air}$, $c_{p,\text{air}} = 1000\,\frac{\text{J}}{\text{kg\,K}}$ and $\rho_\text{air} = 1.2\,\frac{\text{kg}}{\text{m}^3}$
    are the volume flow, specific heat capacity and density of the air, respectively.
    With the assumption that
    \begin{equation}
      (\dot{m}c_p)_\text{air} = (\dot{m}c_p)_\text{cw},
    \end{equation}
    the power of the cooling system is fully determined.
    For the complete formulation, the reader is referred to the model source code.

    The reduced space formulation results in a system model with relatively few expressions, however, since several quantities that are described by ANNs are themselves inputs to other ANNs or used in reformulations within the heat exchangers, the model expressions become deeply nested.
    For this particular use case, the standard SymPy backend (implemented in pure Python) proved to be inefficient as model generation takes about 45 minutes.
    \hyphenation{Sym-En-gine}
    Therefore, SymEngine \citep{certik2019symengine}, a C++ implementation of a subset of SymPy, was integrated as an alternative backend for COMANDO.
    Although SymEngine has a reduced feature set compared to SymPy, all functionality relevant for the presented case study is provided.
    The use of SymEngine reduces the model generation time to about 0.1 seconds.
    Nevertheless, the nested expressions in the model result in very large input files that can take substantial time when written to disk.
    For instance, when using only a single scenario and operating point and maximizing the net power production
    \begin{equation}
        P_\text{net} = P_\text{T} - P_\text{P} - P_\text{CS},
    \end{equation}
    the resulting optimization problem has only 5 variables and 32 constraints.

    In order to solve this problem with BARON \citep{baron20_10_16}, the nonsmooth max function in \cref{eqn:mcp_cw} is approximated with $\text{max}(a, b) \approx 0.5 (a + b + [(a - b + 10^{-4})^2]^{0.5})$ and the $\tanh(x)$ function present in the ANNs is equivalently expressed as $1 - 2/[\exp(2x) + 1]$.
    Generating the BARON input file takes around 1 minute and results in a file size of about 40\,MB.
    This input file is passed to \BARONVERSION{} with absolute and relative optimality tolerances set to 1e-3.
    BARON reports finding a feasible solution with an objective value of $P_\text{net} = 16.48\,\text{MW}$ during preprocessing and terminates after the first iteration and 8 s of CPU time.
    Although a lower bound within the optimality tolerance is given in the log file, BARON states that it cannot guarantee global optimality due to missing bounds for certain nonlinear subexpressions.

    \begin{figure}[!t]
      \centering
      \vspace*{-2.5mm}
      \hspace*{-6mm}
      \input{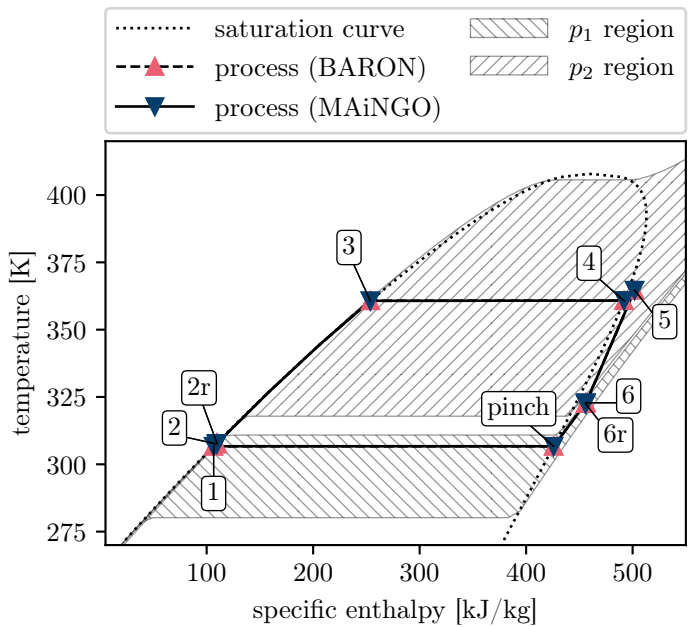}
      \vspace*{-7mm}
      \caption{
        Processes resulting from the optimization using BARON and MAiNGO and boundaries of pressure variables $p_1$ and $p_2$.
      }
      \label{fig:ORC_res}
      \vspace*{-5mm}
    \end{figure}

    To prove the global optimality of this solution, we use the COMANDO interface to the API of our inhouse solver MAiNGO \citep{bongartz2018maingo}.
    MAiNGO automatically provides relaxations of the nested expressions by propagating McCormick relaxations through subexpressions \citep{mitsos2009mccormick}.
    The COMANDO interface uses a SymEngine implementation of \emph{common subexpression elimination} to find subexpressions that occur more than once within the problem description.
    By creating intermediate variables and replacing all occurrences of these subexpressions, a small (21\,kB) input file for MAiNGO can be created.
    Since MAiNGO is capable of propagating McCormick relaxations, the user does not need to provide bounds on these intermediate variables and they are not treated as decision variables, maintaining the reduced-space formulation.
    Solving the resulting problem via MAiNGO version 0.3 with the solution returned by BARON as an initial point takes 22\,s and confirms its global optimality (see \cref{fig:ORC_res}), matching the results reported in \citet{huster2019impact}.

%% file: figures/IES_superstructure.tex
\begin{tikzpicture}
  \tikzset{component/.style={rectangle, minimum size=9mm, draw, fill=white}}
  \tikzset{bus/.style args={#1}{circle, minimum size=4mm, draw, fill=#1, text=white}}
  \tikzset{connector/.style  args={#1}{rectangle, inner sep=0pt, minimum size=1.5mm, draw, fill=#1}}
  \def\lw{0.4pt}
  \def\aw{1pt}
  \tikzset{arrow filled/.style args={#1}{
            -stealth, >=Stealth,line width=2*\lw+\aw, black,
            postaction={draw, -stealth, >=Stealth, #1,
                        line width=\aw,
                        shorten >=2*(\lw)
                        }
           }
  }

  \tikzset{double arrow filled/.style args={#1}{
            <->, >=stealth, line width=2*\lw+\aw, black,
            postaction={draw, <->, >=stealth, #1,
                        line width=\aw,
                        shorten <=2*(\lw),
                        shorten >=2*(\lw)
                        }
           }
  }
	\def\gascolor{fzjgreen}
	\def\eleccolor{fzjyellow}
	\def\heatcolor{fzjred}
	\def\coolcolor{fzjlightblue}
	\def\griddist{3}

	\node[component] (GG) at (-\griddist, 0) {GG};
	\node[connector=\gascolor] (GG_OUT) at (GG.north) {};

  \def\dx{0.1}
  \def\dy{0.05}

	\node[component] (CHP) at (-1.5,0) {CHP};
	\node[connector=\gascolor] (CHP_GAS) at (CHP.north) {};
	\node[connector=\eleccolor] (CHP_ELEC) at (CHP.240) {};
	\node[connector=\heatcolor] (CHP_HEAT) at (CHP.300) {};

	\node[component] (B) at (-0.25, 0) {B};
	\node[connector=\gascolor] (B_GAS) at (B.north) {};
	\node[connector=\heatcolor] (B_HEAT) at (B.south) {};

    \node[component] (STH) at (1,0) {TES\textsubscript{h}};
    \node[connector=\heatcolor] (STH_HEAT) at (STH.south) {};

    \node[component] (HP) at (-0.375, -1.25) {HP};
    \node[connector=\eleccolor] (HP_ELEC) at (HP.west) {};
    \node[connector=\heatcolor] (HP_HEAT) at (HP.east) {};

	\node[component] (PG) at (-\griddist, -1.25) {PG};
	\node[connector=\eleccolor] (PG_BUY) at (PG.15) {};
	\node[connector=\eleccolor] (PG_SELL) at (PG.-15) {};

    \node[component] (BAT) at (-2.975, -2.75) {BAT};
    \node[connector=\eleccolor] (BAT_ELEC) at (BAT.north) {};

	\node[component] (PVOB) at (-1.95, -2.75) {PV$_{\!\scaleto{\text{OB}}{3pt}}$};
	\node[connector=\eleccolor] (PVOB_ELEC) at (PVOB.north) {};

    \node[component] (PVEF) at (-0.95, -2.75) {PV$_{\!\scaleto{\text{EF}}{3pt}}$};
    \node[connector=\eleccolor] (PVEF_ELEC) at (PVEF.north) {};

	\node[component] (TC) at (0.05, -2.75) {CC};
	\node[connector=\eleccolor] (TC_ELEC) at (TC.north) {};
	\node[connector=\coolcolor] (TC_COOL) at (TC.south) {};

  \node[component] (AC) at (1, -2.75) {AC};
  \node[connector=\heatcolor] (AC_HEAT) at (AC.north) {};
  \node[connector=\coolcolor] (AC_COOL) at (AC.east) {};

    \node[component] (STC) at (2.1, -1.5) {TES\textsubscript{c}};

  \node[component] (CON) at (2.35, 0) {DEM};
	\node[connector=\heatcolor] (CON_HEAT) at (CON.235) {};
	\node[connector=\eleccolor] (CON_ELEC) at (CON.270) {};
	\node[connector=\coolcolor] (CON_COOL) at (CON.305) {};

		\node[bus=\gascolor] (GAS) at (-1.5, 1) {};
	\draw[arrow filled=\gascolor] (GG_OUT) |- (GAS);
	\draw[arrow filled=\gascolor] (GAS) -- (CHP_GAS);
	\draw[arrow filled=\gascolor] (GAS) -| (B_GAS);
	  \node[bus=\eleccolor] (ELEC) at (-1.75, -1.25) {};
  \draw[arrow filled=\eleccolor] (CHP_ELEC) -- (ELEC);
  \draw[double arrow filled=\eleccolor] (BAT_ELEC) -- ++(90:0.45) -- ++(0:0.5) -- (ELEC);
  \draw[arrow filled=\eleccolor] (PVOB_ELEC) -- ++(90:0.45) -- (ELEC);
  \draw[arrow filled=\eleccolor] (PVEF_ELEC) -- ++(90:0.15) -- ++(180:0.5) -- (ELEC.275);
  \draw[arrow filled=\eleccolor] (PG_BUY) -- ++(\griddist-2.6, 0) -- (ELEC);
  \draw[arrow filled=\eleccolor] (ELEC) -- ($(PG_SELL) + (\griddist-2.6, 0)$) -- (PG_SELL);
  \draw[arrow filled=\eleccolor] (ELEC) -- (HP_ELEC);
  \draw[arrow filled=\eleccolor] (ELEC) -- ++(-55:0.85) -| (TC_ELEC);
  \draw[arrow filled=\eleccolor] (ELEC) -- ++(-30:1.1) -- ++(0:1.9) -- ++(60:1) -| (CON_ELEC);
  \node[bus=\heatcolor] (HEAT) at (1, -1.25) {};
  \draw[arrow filled=\heatcolor] (CHP_HEAT) -- ++(-90:0.25) -- ++(0:1.3) -- (HEAT.150);
  \draw[arrow filled=\heatcolor] (B_HEAT)-- ++(-90:0.15) -- ++(0:0.88) -- (HEAT.120);
  \draw[double arrow filled=\heatcolor] (STH_HEAT) -- (HEAT);
  \draw[arrow filled=\heatcolor] (HP_HEAT) -- (HEAT);
  \draw[arrow filled=\heatcolor] (HEAT) -- (AC_HEAT);
  \draw[arrow filled=\heatcolor] (HEAT) -- ++(45:0.6) -| (CON_HEAT);
	\node[bus=\coolcolor] (COOL) at (2.25, -2.75) {};
	\draw[arrow filled=\coolcolor] (AC_COOL) -- (COOL);
	\draw[arrow filled=\coolcolor] (TC_COOL) -- ++(-90:0.15) -| (COOL);
    \node[connector=\coolcolor] (STC_COOL) at ($(STC.south west)!(COOL)!(STC.south east)$) {};
  \draw[double arrow filled=\coolcolor] (COOL) -- (STC_COOL);
  \draw[arrow filled=\coolcolor] (COOL) -| (CON_COOL);
\end{tikzpicture}

%% file: figures/BES.tex
\usetikzlibrary{patterns}
\begin{tikzpicture}[x=1cm]
  \tikzset{component/.style={rectangle, minimum size=9mm, draw,fill=white}}
  \tikzset{bus/.style args={#1}{circle, minimum size=3mm, draw, fill=#1, text=white}}
  \tikzset{connector/.style  args={#1}{rectangle, inner sep=0pt, minimum size=2mm, draw, fill=#1}}

  \def\lw{0.4pt}
  \def\aw{1pt}
  \tikzset{arrow filled/.style args={#1}{
            -stealth, >=Stealth,line width=2*\lw+\aw, black,
            postaction={draw, -stealth, >=Stealth, #1,
                        line width=\aw,
                        shorten >=2*(\lw)
                        }
           }
  }
    \tikzset{double arrow filled/.style args={#1}{
            <->, >=stealth, line width=2*\lw+\aw, black,
            postaction={draw, <->, >=stealth, #1,
                        line width=\aw,
                        shorten <=2*(\lw),
                        shorten >=2*(\lw)}
           }
    }

  \def\iso{0.25}
  \def\top{4.25}
  \def\rightt{7}
  \def\bottom{-0.75}

  \draw[fill=gray] (1, \top + \iso) rectangle (\rightt + \iso, \bottom - \iso);
  \draw[fill=white] (1, \top) rectangle (7, \bottom);
  \draw[thick] (3, \top) -- (3, \bottom);
  \fill[pattern=vertical lines] (1, \top) rectangle (3, \bottom);
  \draw[thick] (3, 1) -- (7, 1);
  \fill[pattern=horizontal lines] (3,1) rectangle (7, \bottom);
  \node[component] (HTWE) at (0.95,1) {HT$_\text{wall,E}$};
  \node[connector=fzjred] (HTWEE) at (HTWE.east) {};
  \node[connector=fzjred] (HTWEN) at (HTWE.130) {};
  \node[component] (W) at (1.6, 3) {M$_\text{wall}$};
  \node[connector=fzjred] (WS) at (W.250) {};
  \node[component] (HTAW) at (3.1, 3) {HT$_\text{air,wall}$};
  \node[connector=fzjred] (HTAWE) at (HTAW.-50) {};
  \node[connector=fzjred] (HTAWW) at (HTAW.-130) {};
  \node[component] (A) at (6.25, 3) {M$_{\text{air}}$};
  \node[connector=fzjred] (AE) at (A.west) {};
  \node[component] (HTCW) at (2.9, 0) {HT$_\text{core,wall}$};
  \node[connector=fzjred] (HTCWE) at (HTCW.east) {};
  \node[connector=fzjred] (HTCWW) at (HTCW.130) {};
  \node[component] (C) at (6.4, 0) {M$_\text{core}$};
  \node[connector=fzjred] (CW) at (C.west) {};
  \node[component] (HTAC) at (3.9, 1) {HT$_\text{air,core}$};
  \node[connector=fzjred] (HTACT) at (HTAC.15) {};
  \node[connector=fzjred] (HTACB) at (HTAC.-15) {};
  \node[component] (HP) at (7, 1.75) {HP};
  \node[connector=fzjred] (HPA) at (HP.150) {};
  \node[connector=fzjred] (HPC) at (HP.210) {};
  \node[connector=fzjyellow] (HPP) at (HP.east) {};
  \node[component] (PG) at (8, 0) {PG};
  \node[connector=fzjyellow] (PGN) at (PG.north) {};
  \node[bus=fzjyellow] (PB) at (8,1.25) {};
  \draw[double arrow filled=fzjyellow] (PGN) -- (PB);
  \draw[double arrow filled=fzjyellow] (PB) |- (HPP);
  \draw[double arrow filled=fzjred] (HTWEN) -- +(0, 1);
  \node[bus=fzjred] (WB) at (2.15, 1.75) {};
  \draw[double arrow filled=fzjred] (WS) -- ++(0, -0.5) -- (WB);
  \draw[double arrow filled=fzjred] (WB) |- (HTWEE);
  \draw[double arrow filled=fzjred] (HTAWW)  -- ++(0, -0.5) -- (WB);
  \draw[double arrow filled=fzjred] (HTCWW)  -- ++(0, 0.75) -- (WB);
  \node[bus=fzjred] (AB) at (5.15, 2) {};
  \draw[double arrow filled=fzjred] (AE) -| (AB);
  \draw[double arrow filled=fzjred] (AB) -| (HTAWE);
  \draw[double arrow filled=fzjred] (HTACT) -| (AB);
  \draw[double arrow filled=fzjred] (AB) -- (HPA);
  \node[bus=fzjred] (CB) at (5.15, 0) {};
  \draw[double arrow filled=fzjred] (CW) -- (CB);
  \draw[double arrow filled=fzjred] (HTCWE) -- (CB);
  \draw[double arrow filled=fzjred] (HTACB) -| (CB);
  \draw[double arrow filled=fzjred] (HPC) -- ++(-0.5, 0) -- (CB);
  \node[fill=white, inner sep=0.5mm] at (4.75, -0.465) {\small core };
  \node[fill=white, inner sep=0.5mm] at (2.05, 3.875) {\small wall };
  \node[fill=white, inner sep=0.5mm] at (4.4, 3) {\small air };
  \node[fill=white, inner sep=0.5mm] at (5, 3.875) {$T_\text{air}^\text{min}\!\leq\!T_\text{air}\!\leq\!T_\text{air}^\text{max}$};
\end{tikzpicture}

%% file: figures/DESTEST.tex
\begin{tikzpicture}
  \def\dist{9.85mm}
  \def\bdist{15mm}
  \def\hc{fzjred}
  \def\cc{fzjlightblue}
  \def\pc{fzjyellow}
  \def\gc{fzjgreen}

  \usetikzlibrary{arrows.meta}
  \def\lw{0.4pt}
  \def\aw{1pt}

  \tikzset{arrow filled/.style args={#1}{
            -stealth, >=Stealth,line width=2*\lw+\aw, black,
            postaction={draw, -stealth, >=Stealth, #1,
                        line width=\aw,
                        shorten >=2*(\lw)
                        }
           }
  }
  \tikzset{pipe/.style args={#1}{
            ,line width=2.5 * \lw + 1.5 * \aw, black,
            postaction={draw, #1,
                        line width=1.5 * \aw,
                        shorten <=1.25 * (\lw),
                        shorten >=1.25 * (\lw)
                        }
           }
  }
  \def\cw{9mm}
  \tikzset{component/.style={rectangle, minimum size=\cw, draw, fill=white}}
  \tikzset{optional component/.style={rectangle, minimum size=8mm, draw, fill=white}}
  \tikzset{bus/.style args={#1}{circle, minimum size=4mm, inner sep=0, draw, fill=#1, text=white}}
  \tikzset{minibus/.style args={#1}{circle, minimum size=1mm, inner sep=0, draw, fill=#1}}
  \tikzset{connector/.style  args={#1}{rectangle, inner sep=0pt, minimum size=1.5mm, draw, fill=#1}}
  \tikzset{optional connector/.style  args={#1}{rectangle, inner sep=0pt, minimum size=1.5mm, draw, fill=#1}}

  \def\baseheight{0.1 * \bdist}
  \node[component] (WH) at (0, \baseheight) {WH};
  \node[connector=\hc] (WHOUT) at (WH.-35) {};

  \node[component] (LNW) at (1.5 * \dist, \baseheight) {L$_\text{NW}$};
  \node[connector=\hc] (LNWIN) at (LNW.215) {};
  \node[connector=\hc] (LNWOUT) at (LNW.-35) {};
  \node[connector=\pc] (LNWP) at (LNW.110) {};

  \def\ppos{1 * \bdist}
  \node[component] (PG) at (0, \ppos) {PG};
  \node[connector=\pc] (POUT) at (PG.east) {};

  \node[bus=\pc] (PB) at (POUT -| LNWP) {};

  \def\gpos{2 * \bdist}
  \node[component] (GG) at (0, \gpos) {GG};
  \node[connector=\gc] (GOUT) at (GG.30) {};

  \node[bus=\gc] (GB) at (GOUT -| LNWP) {};

  \def\cs{2mm}
  \def\nwpos{4.75 * \dist}

  \node[component] (CG40) at (\nwpos - 1.5 * \cs - 1.5 * \cw, 1 * \bdist) {CG$_{40}$};
  \node[connector=\hc] (CG40IN) at (CG40.255) {};
  \node[connector=\gc] (CG40G) at (CG40.55) {};
  \node[connector=\pc] (CG40P) at (CG40.125) {};

  \node[component] (CG50) at (\nwpos - 0.5 * \cs - 0.5 * \cw, 1 * \bdist) {CG$_{50}$};
  \node[connector=\hc] (CG50IN) at (CG50.255) {};
  \node[connector=\gc] (CG50G) at (CG50.55) {};
  \node[connector=\pc] (CG50P) at (CG50.125) {};

  \node[component] (CG70) at (\nwpos + 0.5 * \cs + 0.5 * \cw, 1 * \bdist) {CG$_{70}$};
  \node[connector=\hc] (CG70IN) at (CG70.255) {};
  \node[connector=\gc] (CG70G) at (CG70.55) {};
  \node[connector=\pc] (CG70P) at (CG70.125) {};

  \node[component] (CG85) at (\nwpos + 1.5 * \cs + 1.5 * \cw, 1 * \bdist) {CG$_{85}$};
  \node[connector=\hc] (CG85IN) at (CG85.255) {};
  \node[connector=\gc] (CG85G) at (CG85.55) {};
  \node[connector=\pc] (CG85P) at (CG85.125) {};

  \node[component, align=center, minimum width=1mm + 4 * \cw + 3 * \cs] (NW) at (\nwpos, \baseheight) {
  };
  \node[connector=\hc] (NWIN) at ($(NW.south west)!(LNWOUT)!(NW.north west)$) {};
  \node[connector=\hc] (NWOUT40) at ($(NW.north west)!(CG40IN)!(NW.north east)$) {};
  \node[connector=\hc] (NWOUT50) at ($(NW.north west)!(CG50IN)!(NW.north east)$) {};
  \node[connector=\hc] (NWOUT70) at ($(NW.north west)!(CG70IN)!(NW.north east)$) {};
  \node[connector=\hc] (NWOUT85) at ($(NW.north west)!(CG85IN)!(NW.north east)$) {};
  \coordinate (NWMID) at ($(NWOUT50)!0.5!(NWOUT70)$);
  \node at ($(NW.south east) + (-0.35, 0.25)$) {NW};

  \draw[arrow filled=\pc] (POUT) -- (PB);
  \draw[arrow filled=\pc] (PB) -- +(0:0.5) |- ($(CG40P) + (0, 0.35)$) -- (CG40P);
  \draw[arrow filled=\pc] (PB) -- +(30:0.385) |- ($(CG50P) + (0, 0.5)$) -- (CG50P);
  \draw[arrow filled=\pc] (PB) -- +(60:1/3) |- ($(CG70P) + (0, 0.65)$) -- (CG70P);
  \draw[arrow filled=\pc] (PB) -- +(90:0.5) |- ($(CG85P) + (0, 0.8)$) -- (CG85P);

  \draw[arrow filled=\gc] (GOUT) -- (GB);
  \draw[arrow filled=\gc] (GB) -- +(-60:0.385) -| (CG40G);
  \draw[arrow filled=\gc] (GB) -- +(-30:1/3) -| (CG50G);
  \draw[arrow filled=\gc] (GB) -- +(0:0.5) -| (CG70G);
  \draw[arrow filled=\gc] (GB) -- +(30:1/3) -| (CG85G);

  \draw[arrow filled=\pc] (PB) -- (LNWP);
  \draw[arrow filled=\hc] (WHOUT) -- (LNWIN);
  \draw[arrow filled=\hc] (LNWOUT) -- (NWIN);
  \draw[arrow filled=\hc] (NWOUT40) -- (CG40IN);
  \draw[arrow filled=\hc] (NWOUT50) -- (CG50IN);
  \draw[arrow filled=\hc] (NWOUT70) -- (CG70IN);
  \draw[arrow filled=\hc] (NWOUT85) -- (CG85IN);

  \node[minibus=\hc] (NWCENTER) at ($(NWIN -| NWMID) + (0, 0.4)$) {};
  \node[minibus=\hc] (PIPE40) at (NWOUT40 |- NWCENTER) {};
  \node[minibus=\hc] (PIPE50) at (NWOUT50 |- NWCENTER) {};
  \node[minibus=\hc] (PIPE70) at (NWOUT70 |- NWCENTER) {};
  \node[minibus=\hc] (PIPE85) at (NWOUT85 |- NWCENTER) {};

  \draw[arrow filled=\hc] (NWIN) -| (NWCENTER);
  \draw[pipe=gray] (NWCENTER) -- (PIPE50);
  \draw[pipe=gray] (PIPE50) -- (PIPE40);
  \draw[pipe=gray] (NWCENTER) -- (PIPE70);
  \draw[pipe=gray] (PIPE70) -- (PIPE85);

  \draw[arrow filled=\hc] (PIPE40) -- (NWOUT40);
  \draw[arrow filled=\hc] (PIPE50) -- (NWOUT50);
  \draw[arrow filled=\hc] (PIPE70) -- (NWOUT70);
  \draw[arrow filled=\hc] (PIPE85) -- (NWOUT85);

  \coordinate (LINELEFT) at ($(GG.north west) + (0, 0.25)$);
  \coordinate (LINERIGHT) at ($(NW.east)!(LINELEFT)!(CG85.east)$);
  \coordinate (LBOTTOMLEFT) at ($(LINELEFT) + (0, 0.25)$);

  \node[component] (HP) at ($(LBOTTOMLEFT) + (0.2cm + 0.5 * \cw, 1cm + 0.5 * \cw)$) {HP};
  \node[connector=\hc] (HPIN) at (HP.290) {};
  \node[connector=\pc] (HPP) at (HP.250) {};
  \node[connector=\hc] (HPOUT) at (HP.north) {};

  \node[component] (HX) at ($(LBOTTOMLEFT) + (1.3cm + 0.5 * \cw, 1cm + 0.5 * \cw)$) {HE};
  \node[connector=\hc] (HXIN) at (HX.south) {};
  \node[connector=\hc] (HXOUT) at (HX.north) {};

  \coordinate (LBOTTOMRIGHT) at ($(HX.south east) + (0.25, -1)$);
  \coordinate (LTOPLEFT) at ($(HP.north west) + (-0.25, 0.9)$);
  \coordinate (LTOPRIGHT) at ($(HX.north east) + (0.25, 0.9)$);

  \draw (LBOTTOMLEFT) rectangle (LTOPRIGHT);
  \node[connector=\hc] (LIN) at ($(LBOTTOMLEFT)!0.5!(LBOTTOMRIGHT)$) {};
  \node[connector=\hc] (LOUT) at ($(LTOPLEFT)!0.5!(LTOPRIGHT)$) {};
  \node[connector=\pc] (LP) at ($(LBOTTOMLEFT) + (0, 0.6)$) {};

  \node[bus=\hc] (LHI) at ($(LIN) + (0, 0.6)$) {};
  \node[bus=\hc] (LHO) at ($(LOUT) + (0, -0.6)$) {};

  \draw[arrow filled=\hc] (LIN) -- (LHI);
  \draw[arrow filled=\hc] (LHI) -| (HPIN);
  \draw[arrow filled=\hc] (LHI) -| (HXIN);
  \draw[arrow filled=\hc] (HPOUT) |- (LHO);
  \draw[arrow filled=\hc] (HXOUT) |- (LHO);
  \draw[arrow filled=\hc] (LHO) -- (LOUT);
  \draw[arrow filled=\pc] (LP) -| (HPP);

  \node at ($(LBOTTOMRIGHT) + (-0.25, 0.25)$) {L};

  \coordinate (CGTOPLEFT) at ($(LTOPRIGHT) + (0.25, 0)$);
  \coordinate (CGBOTTOMLEFT) at (LBOTTOMRIGHT -| CGTOPLEFT);

    \coordinate (CGBOTTOMRIGHT) at (LBOTTOMRIGHT -| CG85.east);
  \coordinate (CGTOPRIGHT) at (CGTOPLEFT -| CGBOTTOMRIGHT);

  \draw (CGBOTTOMLEFT) rectangle (CGTOPRIGHT);

  \def\usedCG{CG50}
  \draw[very thick, dotted] (\usedCG.north west) -- (CGBOTTOMLEFT);
  \draw[very thick, dotted] (\usedCG.north east) -- (CGBOTTOMRIGHT);

  \def\boty{3.25 * \bdist}
  \def\topy{4.1 * \bdist}
  \def\lefx{\nwpos - 1.1 * \cw}
  \def\rigx{\nwpos + 1.1 * \cw}
  \node[component] (CGL) at (\lefx, \boty) {L};
  \node[connector=\hc] (CGLIN) at (CGL.south) {};
  \node[connector=\pc] (CGLP) at (CGL.200) {};
  \node[connector=\hc] (CGLOUT) at (CGL.north) {};

  \node[component] (HR) at (\lefx, \topy) {HS$_\text{HR}$};
  \node[connector=\pc] (HRIN) at (HR.180) {};
  \node[connector=\hc] (HROUT) at (HR.0) {};

  \node[component] (B) at (\rigx, \topy) {HS$_\text{B}$};
  \node[connector=\gc] (BIN) at (B.0) {};
  \node[connector=\hc] (BOUT) at (B.180) {};

  \node[component] (DEM) at (\rigx, \boty) {DEM};
  \node[connector=\hc] (DEMIN) at (DEM.west) {};

  \node[connector=\hc] (CGIN) at ($(CGBOTTOMLEFT)!(CGLIN)!(CGBOTTOMRIGHT)$) {};
  \node[connector=\pc] (CGP) at ($(CGTOPLEFT) + (0.3, 0)$) {};
  \node[connector=\gc] (CGG) at ($(CGTOPRIGHT) + (-0.25, 0)$) {};

  \coordinate (CGMIDHEIGHT) at ($(HR.south)!0.5!(CGL.north)$);
  \coordinate (CGMIDWIDTH) at ($(HR.east)!0.5!(B.west)$);

  \node[bus=\pc] (CGPB) at (CGMIDHEIGHT -| CGP) {};
  \node[bus=\hc] (CGHB) at (CGMIDHEIGHT -| CGMIDWIDTH) {};

  \draw[arrow filled=\hc] (CGIN) -- ++(0, 0.2) -| (CGLIN);
  \draw[arrow filled=\hc] (CGLOUT) |- (CGHB);
  \draw[arrow filled=\hc] (HROUT) -| ($(CGHB) + (120:0.5)$) -- (CGHB);
  \draw[arrow filled=\hc] (BOUT) -| ($(CGHB) + (60:0.5)$) -- (CGHB);
  \draw[arrow filled=\hc] (CGHB) |- (DEMIN);

  \draw[arrow filled=\pc] (CGP) -- (CGPB);
  \draw[arrow filled=\pc] (CGPB) |- (CGLP);
  \draw[arrow filled=\pc] (CGPB) -- ++(60:0.385) |- (HRIN);

  \draw[arrow filled=\gc] (CGG) |- (BIN);

  \def\usedL{CGL}
  \draw[very thick, dotted] (\usedL.south west) -- (LBOTTOMRIGHT);
  \draw[very thick, dotted] (\usedL.north west) -- (LTOPRIGHT);

  \node at ($(CGBOTTOMRIGHT) + (-0.35, 0.25)$) {CG};
\end{tikzpicture}

%% file: figures/ltdh_result.tex
\begin{tikzpicture}
  \def\dist{10mm}
  \def\bdist{15mm}
  \def\hc{fzjred}
  \def\cc{fzjlightblue}
  \def\pc{fzjyellow}
  \def\gc{fzjgreen}

  \usetikzlibrary{arrows.meta}
  \def\lw{0.4pt}
  \def\aw{1pt}

  \tikzset{arrow filled/.style args={#1}{
            -stealth, >=Stealth,line width=2 * \lw + \aw, black,
            postaction={draw, -stealth, >=Stealth, #1,
                        line width=\aw,
                        shorten >=2*(\lw)
                        }
           }
  }
  \tikzset{pipe/.style args={#1}{
            ,line width=2.5 * \lw + 1.5 * \aw, black,
            postaction={draw, #1,
                        line width=1.5 * \aw,
                        shorten <=1.25 * (\lw),
                        shorten >=1.25 * (\lw)
                        }
           }
  }
  \def\cw{9mm}
  \tikzset{component/.style={rectangle, minimum size=\cw, draw, fill=white}}
  \tikzset{optional component/.style={rectangle, minimum size=8mm, draw,
  fill=white}}
  \tikzset{bus/.style args={#1}{circle, minimum size=4mm, inner sep=0, draw, fill=#1, text=white}}
  \tikzset{minibus/.style args={#1}{circle, minimum size=1mm, inner sep=0, draw, fill=#1}}
  \tikzset{connector/.style  args={#1}{rectangle, inner sep=0pt, minimum size=1.5mm, draw, fill=#1}}
  \tikzset{optional connector/.style  args={#1}{rectangle, inner sep=0pt, minimum size=1.5mm, draw,
  fill=#1}}

  \node[component] (WH) at (0, 0) {WH};
  \node[connector=\hc] (WHOUT) at (WH.-35) {};

  \def\cws{12mm}
  \node[component, align=center, minimum width=\cws] (LNW) at (1.625 * \dist, 0) {L$_\text{NW}$\\[-1.5mm]
  \tiny 102.3\,kW HP};
  \node[connector=\hc] (LNWIN) at ($(LNW.south west)!(WHOUT)!(LNW.north west)$) {};
  \node[connector=\hc] (LNWOUT) at ($(LNW.south east)!(LNWIN)!(LNW.north east)$) {};
  \node[connector=\pc] (LNWP) at (LNW.110) {};

  \def\ppos{1 * \bdist}
  \node[component] (PG) at (0, \ppos) {PG};
  \node[connector=\pc] (POUT) at (PG.east) {};

  \node[bus=\pc] (PB) at (POUT -| LNWP) {};

  \def\gpos{1.7 * \bdist}
  \node[component] (GG) at (0, \gpos) {GG};
  \node[connector=\gc] (GOUT) at (GG.20) {};

  \node[bus=\gc] (GB) at (GOUT -| LNWP) {};

  \def\cs{2mm}
  \def\nwpos{5.5 * \dist}

  \node[component, align=center, minimum width=\cws] (CG40) at (\nwpos - 2.5 * \cs - 1.5 * \cws, 1 * \bdist) {CG$_{40}$\\[-1.5mm]
  \tiny 5.0\,kW HR\\[-2mm]
  \tiny 12.8\,kW HE};
  \node[connector=\hc] (CG40IN) at (CG40.255) {};
  \node[connector=\pc] (CG40P) at (CG40.125) {};

  \node[component, align=center, minimum width=\cws] (CG50) at (\nwpos - 1 * \cs - 0.5 * \cw, 1 * \bdist) {CG$_{50}$\\[-1.5mm]
  \tiny 5.0\,kW HR\\[-2mm]
  \tiny 12.8\,kW HE};
  \node[connector=\hc] (CG50IN) at (CG50.255) {};
  \node[connector=\pc] (CG50P) at (CG50.125) {};

  \node[component, align=center, minimum width=\cws] (CG70) at (\nwpos + 1 * \cs + 0.8 * \cw, 1 * \bdist) {CG$_{70}$\\[-1.5mm]
  \tiny 17.8\,kW B};
  \node[connector=\gc] (CG70G) at (CG70.55) {};

  \node[component, align=center, minimum width=\cws] (CG85) at (\nwpos + 2.5 * \cs + 2.1 * \cw, 1 * \bdist) {CG$_{85}$\\[-1.5mm]
  \tiny 17.8\,kW B};
  \node[connector=\gc] (CG85G) at (CG85.55) {};

  \node[component, align=center, minimum width=1mm + 4 * \cws + 3 * \cs] (NW) at (\nwpos, 0) {};
  \node[connector=\hc] (NWIN) at ($(NW.south west)!(LNWOUT)!(NW.north west)$) {};
  \node[connector=\hc] (NWOUT40) at ($(NW.north west)!(CG40IN)!(NW.north east)$) {};
  \node[connector=\hc] (NWOUT50) at ($(NW.north west)!(CG50IN)!(NW.north east)$) {};
  \coordinate (NWMID) at ($(NW.west)!0.5!(NW.east)$);
  \node at ($(NW.south east) + (-0.35, 0.25)$) {NW};
  \node[anchor=north east] at ($(NW.north east) + (0.05, 0.05)$) {
    \footnotesize
    $T^\text{re}_\text{NW}\!\in\![25\si{\celsius}, 35\si{\celsius}]$
  };

  \draw[arrow filled=\pc] (POUT) -- (PB);
  \draw[arrow filled=\pc] (PB) -- (LNWP);
  \draw[arrow filled=\pc] (PB) -- +(0:0.5) |- ($(CG40P) + (0, 0.35)$) -- (CG40P);
  \draw[arrow filled=\pc] (PB) -- +(30:0.385) |- ($(CG50P) + (0, 0.5)$) -- (CG50P);

  \draw[arrow filled=\gc] (GOUT) -- (GB);
  \draw[arrow filled=\gc] (GB) -- +(0:0.5) -| (CG70G);
  \draw[arrow filled=\gc] (GB) -- +(30:1/3) -| (CG85G);

  \draw[arrow filled=\hc] (WHOUT) -- (LNWIN);
  \draw[arrow filled=\hc] (LNWOUT) -- (NWIN);
  \draw[arrow filled=\hc] (NWOUT40) -- (CG40IN);
  \draw[arrow filled=\hc] (NWOUT50) -- (CG50IN);

  \node[minibus=\hc] (NWCENTER) at ($(NWIN -| NWMID) + (0, 0.4)$) {};
  \node[minibus=\hc] (PIPE40) at (NWOUT40 |- NWCENTER) {};
  \node[minibus=\hc] (PIPE50) at (NWOUT50 |- NWCENTER) {};

  \draw[arrow filled=\hc] (NWIN) -| (NWCENTER);
  \draw[pipe=gray] (NWCENTER) -- (PIPE50);
  \draw[pipe=gray] (PIPE50) -- (PIPE40);

  \draw[arrow filled=\hc] (PIPE40) -- (NWOUT40);
  \draw[arrow filled=\hc] (PIPE50) -- (NWOUT50);
\end{tikzpicture}

%% file: figures/ORC.tex
\begin{tikzpicture}
  \usetikzlibrary{arrows.meta}
  \def\lw{0.4pt}
  \def\aw{1pt}

  \tikzset{arrow filled/.style args={#1}{
            -stealth, >=Stealth,line width=2*\lw+\aw, black,
            postaction={draw, -stealth, >=Stealth, #1,
                        line width=\aw,
                        shorten >=2*(\lw)
                        }
           }
  }
  \tikzset{state/.style={circle, fill=white, inner sep=0, minimum size=3.5mm, draw}}
  \tikzset{component/.style={rectangle, minimum size=9mm, draw}}
  \tikzset{bus/.style args={#1}{circle, minimum size=3mm, draw, fill=#1, text=white}}
  \tikzset{connector/.style  args={#1}{rectangle, inner sep=0pt, minimum size=1.3mm, draw, fill=#1}}

      \newcommand{\Pump}[4][0]
      {
          \def\cen{(#2, #3)}
          \node[circle, draw, thick, minimum size=2*\r, inner sep=0] (#4) at \cen {#4};
          \node[rectangle, draw, thick, minimum size=2.1*\r, inner sep=0] at \cen {\phantom{(#4)}};
          \draw[line width=1pt, rotate around={#1:\cen}]
              (#2 cm - \r, #3 cm) -- (#2 cm, #3 cm + \r) -- (#2 cm + \r, #3 cm);
      }
    \newcommand{\HE}[5][0]
      {
          \def\cen{(#2, #3)}
          \node[circle, draw, thick, minimum size=2*\r, label=90+#1:#4] (#5) at \cen {};
          \node[rectangle, draw, thick, minimum width=2.2*\r, minimum height=2.4*\r, yshift=-0.15*\r, rotate around={#1:(#5)}] at \cen {};
          \draw[line width=1pt,sharp corners, rotate around={#1:\cen}]
              (#2 cm - 2/3 * \r, #3 cm - 4/3 * \r) -- (#2 cm - 2/3 * \r, #3 cm + \r/3) -- (#2 cm, #3 cm - \r/3) -- (#2 cm + 2/3 * \r, #3 cm + \r/3) -- (#2 cm + 2/3 * \r, #3 cm - 4/3 * \r);
      }
    \newcommand{\Turbine}[3]
      {
          \def\cen{(#1, #2)}
          \draw[draw, thick] (#1 - 0.4,#2 - 0.4) coordinate (#3bl) rectangle ++(0.8, 0.8) coordinate (#3tr);
          \node[trapezium, thick, minimum height=0.35cm, trapezium angle=60, draw, rotate=90, label=center:#3] (#3) at (\cen) {};
          \node[connector=fzjgray] (#3in) at ($(#3.top right corner)!(#3tr)!(#3.top left corner)$) {};
          \node[connector=fzjgray] (#3out) at (#3.bottom left corner) {};
      }

  \def\cw{fzjblue!30}
  \def\ww{fzjblue!70}

  \def\r{0.3 cm}
  \Pump{0.7}{.95}{P};
  \HE[180]{0.9}{-.4}{HE$_\text{con}$}{CON};
  \HE{1.7}{1.5}{HE$_\text{rec}$}{REC};
  \HE[180]{3.1}{1.5}{HE$_\text{eco}$}{ECO};
  \HE[180]{4.5}{1.5}{HE$_\text{eva}$}{EVA};
  \HE[180]{5.9}{1.5}{HE$_\text{sup}$}{SUP};
  \def\Tx{7.05}
  \def\Ty{0.65}
  \draw[draw, thick] (\Tx - 0.3, \Ty - 0.5) coordinate (Tbl) rectangle ++(0.6, 1) coordinate (Ttr);
  \node[trapezium, draw, thick, minimum height=0.35cm, trapezium angle=60, yshift=0.5mm, rotate=90, label=center:T] (T) at (\Tx, \Ty) {};
  \node[connector=gray] (Tin) at ($(T.top right corner)!(Ttr)!(T.top left corner)$) {};
  \draw[] (Tin) -- (T.top right corner) node[circle, draw, inner sep = 0, minimum size = 0.4] {};
  \node[connector=fzjgray] (Tout) at (T.bottom left corner)at ($(Tbl)!(T.bottom left corner)!($(Tbl) + (1, 0)$)$) {};
  \draw[] (T.bottom left corner) node[circle, draw, inner sep = 0, minimum size = 0.4] {} -- (Tout);

  \def\bottomright{2/3 * \r, - 4/3 * \r}
  \def\bottomleft{-2/3 * \r, - 4/3 * \r}
  \def\topright{2/3 * \r, 4/3 * \r}
  \def\topleft{-2/3 * \r, 4/3 * \r}

  \node[connector=gray] (Pout) at (P.north) {};
  \node[connector=fzjgray] (Pin) at (P.south) {};
  \node[connector=gray] (RECcin) at (REC.west) {};
  \node[connector=gray] (RECcout) at (REC.east) {};
  \node[connector=fzjgray] (REChin) at ($(REC) + (\bottomright)$) {};
  \node[connector=fzjgray] (REChout) at ($(REC) + (\bottomleft)$) {};

  \node[connector=gray] (ECOcout) at (ECO.east) {};
  \node[connector=gray] (ECOcin) at (ECO.west) {};
  \node[connector=fzjred!50] (ECOhin) at ($(ECO) + (\topright)$) {};
  \node[connector=fzjred!30] (ECOhout) at ($(ECO) + (\topleft)$) {};

  \node[connector=gray] (EVAcout) at (EVA.east) {};
  \node[connector=gray] (EVAcin) at (EVA.west) {};
  \node[connector=fzjred!80] (EVAhin) at ($(EVA) + (\topright)$) {};
  \node[connector=fzjred!50] (EVAhout) at ($(EVA) + (\topleft)$) {};

  \node[connector=gray] (SUPcout) at (SUP.east) {};
  \node[connector=gray] (SUPcin) at (SUP.west) {};
   \node[connector=fzjred] (SUPhin) at ($(SUP) + (\topright)$) {};
  \node[connector=fzjred!80] (SUPhout) at ($(SUP) + (\topleft)$) {};

  \node[connector=\ww] (CONcout) at (CON.east) {};
  \node[connector=\cw] (CONcin) at (CON.west) {};
  \node[connector=fzjgray] (CONhin) at ($(CON) + (\topright)$) {};
  \node[connector=fzjgray] (CONhout) at ($(CON) + (\topleft)$) {};

  \draw[arrow filled=gray] (Pout) |- node (Pup) {} (RECcin);
  \draw[arrow filled=gray] (RECcout) -- (ECOcin);
  \draw[arrow filled=gray] (ECOcout) -- (EVAcin);
  \draw[arrow filled=gray] (EVAcout) -- (SUPcin);

  \draw[arrow filled=gray] (SUPcout) -| (Tin);
  \draw[arrow filled=gray!30] (Tout) -- ++(0, -0.15) coordinate (Tdown) -| coordinate (Tleft) (REChin);

  \draw[arrow filled=gray!30] (REChout) -- ++(0, -0.6) coordinate (RECdown) -| (CONhin);
  \draw[arrow filled=gray!30] (CONhout) -- (Pin);

  \draw[arrow filled=fzjred] (SUPhin) ++(0, 0.75 cm) coordinate (GTr) -- (SUPhin);
  \draw[arrow filled=fzjred!80] (SUPhout) -- ++(0, 0.5 cm) -| (EVAhin);
  \draw[arrow filled=fzjred!50] (EVAhout) -- ++(0, 0.5 cm) -| (ECOhin);
  \draw[arrow filled=fzjred!30] (ECOhout) -- ++(0, 0.75cm) coordinate (GTl);
  \node at ($(GTl)!0.5!(GTr)$) {geothermal brine};

  \node[state] (s2) at ($(Pup)!0.21!(REC) + (0, 0.5)$) {\small 2};
  \draw (s2) -- ++(0, -0.48);
  \node[state](s2r)  at ($(REC)!0.45!(ECO) + (0, 0.5)$) {\small 2r};
  \draw (s2r) -- ++(0, -0.48);
  \node[state] (s3) at ($(ECO)!0.45!(EVA) + (0, 0.5)$) {\small 3};
    \draw (s3) -- ++(0, -0.48);
  \node[state] (s4) at ($(EVA)!0.45!(SUP) + (0, 0.5)$) {\small 4};
    \draw (s4) -- ++(0, -0.48);
  \node[state] (s5) at ($(SUP)!0.65!(SUP-|T.top right corner) + (0, 0.5)$) {\small 5};
    \draw (s5) -- ++(0, -0.48);
  \node[state] (IB) at ($(Tdown)!0.55!(Tleft)$) {\small 6};
  \node[anchor=south west, xshift=5mm] at (IB) {isobutane};
  \node[state] (s6r) at ($(RECdown) + (0, -0.5)$) {\small 6r};
  \draw (s6r) -- ++(0, 0.48);
  \node[state] (s1) at ($(CONhout)!0.4!(Pin) + (-0.5, 0)$) {\small 1};
  \draw (s1) -- ++(0.48, 0);

  \node[label=above:{CS}] (F) at (4, -1) {};
   \def\radius{0.05}
   \def\angle{10}
   \pgfmathsetmacro{\length}{\radius/tan(\angle)}
   \pgfmathsetmacro{\width}{\radius*(1 + 1/sin(\angle))}

   \draw ($(F) + (-\width-\radius, -2*\radius)$)  rectangle ($(F) + (\width+\radius, +2*\radius)$);

   \draw[thick] ($(F) + (-\width-\radius, -6*\radius)$)  rectangle ($(F) + (\width+\radius, +2*\radius)$);
   \node[connector=\ww] (Fin) at ($(F) + (\width+\radius, -4*\radius)$) {};
   \node[connector=\cw] (Fout) at ($(F) + (-\width-\radius, -4*\radius)$) {};

  \draw (F.center) -- +(180-\angle:\length) arc (90-\angle:270+\angle:\radius) -- (F.center)
	-- +(-\angle:\length) arc (-90-\angle:90+\angle:\radius) -- (F.center);

  \draw[arrow filled=\cw] (Fout) -| ($(CON) + (-1cm, 0 cm)$) -- (CONcin);
  \draw[arrow filled=\ww] (CONcout) -- ++(3.9, 0) coordinate (c1) |- coordinate (c2) (Fin);
  \node[anchor=west] at ($(c1)!0.5!(c2)$) {cooling water};
  \draw[ultra thick] (Fin) -- (Fout);
  \node[left color=\cw, right color=\ww, minimum width=7mm, inner sep = 0] at ($(Fin)!0.5!(Fout)$) {};

  \draw[->] ($(P) - (0.9, 0)$) -- node[above]{$P_\text{P}$} ($(P) - (0.4, 0)$);
  \draw[->] ($(T) + (0.4, 0)$) -- node[above]{$P_\text{T}$} ($(T) + (0.9, 0)$);
  \draw[->] ($(F) - (1, 0)$) -- node[above]{$P_\text{CS}$} ($(F) - (0.5, 0)$);
\end{tikzpicture}

%% file: sections/5_conclusion.tex

\section{Conclusion} \label{sec:conclusion}
  We present COMANDO, our flexible open-source framework for \acronym of energy systems.
  COMANDO combines desirable features of existing tools and provides layers of abstraction suitable for structured model generation and flexible problem formulation.
  The behavior of individual components can be represented with detailed models, including dynamic and nonlinear effects based on mechanistic, data-driven or hybrid modeling approaches.
  The component models are then aggregated to energy system models, based on which different optimization problems concerning the design and/or operation of the energy system can be formulated.
  COMANDO natively allows to consider multiple operating scenarios via stochastic programming formulations, allowing to find system designs that are suitable for operation under uncertainty.
  The resulting problem formulations can either be manipulated in user-defined algorithms, or be passed to algebraic modeling languages or directly to solvers.

  COMANDO allows for flexible model creation beyond the capabilities of existing MILP-based energy-system modeling tools and provides a wide range of options for problem formulation.
  Contrary to classical algebraic modeling frameworks, it allows for modular component and system representations, and is dedicated to energy system design and operation.

  In four case studies, we demonstrate how COMANDO can be used to create modular and reusable component and system models of various types of energy systems. Further, we formulate and solve associated optimization problems.
  With COMANDO, we facilitate and enhance workflows of computer-based analysis of future integrated energy systems.
  We plan to continuously improve and expand COMANDO's capabilities, with future versions being published via the \citet{COMANDO_REPO}.

%% file: sections/6_acknowledgement.tex

\section*{Author Contribution}
\begin{itemize}[noitemsep]
  \item ML developed COMANDO, and wrote \cref{sec:introduction,sec:M_AND_O,sec:COMANDO,sec:ORC,sec:conclusion} in close collaboration with MD and help
  \& guidance from AM.
  \item DS and ML contributed the code for automatic linearization, the first case study and wrote \cref{sec:design_problem} with help and guidance from MD and AB.
  \item FB incorporated Pyomo.DAE into the Pyomo interface, contributed the second case study and wrote \cref{sec:dynamic_operation} with help and guidance from ML, MD and AB.
  \item DH and ML contributed the third case study and wrote \cref{sec:HP_net} with help and guidance from MD, AX and DM.
  \item UB and MD gave conceptual input for the creation of COMANDO.
  \item MD supervised the writing process.
  \item All authors reviewed and edited the manuscript
\end{itemize}

\section*{Declaration of Competing Interest}
We have no conflict of interest.

\section*{Acknowledgements}
We would like to thank Alexander Holtwerth (Forschungszentrum Jülich GmbH, Institute of Energy and Climate Research, Energy Systems Engineering (IEK-10)) for providing an initial version of the Gurobi interface, and the anonymous reviewer for valuable comments and suggestions that have considerably improved this manuscript.
This work was funded by the Helmholtz Association of German Research Centers through program-oriented funding, the Joint Initiative ``Energy System 2050 -- A Contribution of the Research Field Energy'', and the Initiative ``Energy System Integration''.

%% file: sections/7_appendices.tex
\section*{Nomenclature} \label{sec:nomenclature}
\noindent\textbf{Acronyms} \\
\noindent
\begin{tabularx}{\columnwidth}{lX}
    AML & algebraic modeling language\\
    ANN & artificial neural network\\
    API & application programming interface\\
    COP & coefficient of performance\\
    DAMF & differential-algebraic modeling frame\-work\\
    ESMF & energy system modeling framework \\
    GWI & global warming impact \\
    LP & linear programming \\
    MIDO & mixed-integer dynamic optimization \\
    MILP & mixed-integer linear programming \\
    MINLP & mixed-integer nonlinear programming \\
    MIQCQP & mixed-integer quadratically constrained quadratic programming\\
    NLP & nonlinear programming\\
    ORC & organic Rankine cycle\\
    TAC & total annualized costs
\end{tabularx} \\[5mm]
%
\noindent\textbf{Component labels} \\
\noindent
\begin{tabularx}{\columnwidth}{lX}
  AC & absorption chiller \\
  B & boiler \\
  BAT & battery \\
  CG & consumer group subsystem\\
  CC & compression chiller \\
  CHP & combined heat-and-power unit \\
  CS & cooling system \\
  DEM & demand \\
  GG & gas grid \\
  HE & heat exchanger \\
  HP & heat pump \\
  HR & heating rod \\
  HS & heat source \\
  L & linking subsystem \\
\end{tabularx} \\[5mm]
\noindent\textbf{Component labels (cont.)} \\
\noindent
\begin{tabularx}{\columnwidth}{lX}
  NW & network \\
  P & pump \\
  PG & power grid \\
  PV & photovoltaic unit \\
  T & turbine \\
  TES & thermal energy storage \\
  WH & waste heat
\end{tabularx} \\[5mm]
%
\noindent\textbf{Latin symbols} \\
\noindent
\begin{tabularx}{\columnwidth}{lX}
  $A$ & contact area [m$^2$]\\
  $b$ & build decision (1: build, 0: do not build) \\
  $c$ & generic connector expression \\
  $c_p$ & heat capacity [J/kg/K]\\
  $C$ & cost [\euro]\\
  $\bm{e}$ & generic algebraic expression \\
  $E, \dot{E}$ & generalized energy, energy flow [J], [W] \\
  $F$ & generic objective function \\
  $\bm{g}$ & left-hand side of generic inequality constraints \\
  $h$ & specific enthalpy [J/kg] \\
  $\bm{h}$ & left-hand side of generic equality constraints \\
  $\mathcal{I}$ & set of components \\
  $\dot{m}$ & mass flow rate [kg/s] \\
  $M$ & investment cost exponent \\
  $p$ & pressure [Pa]\\
  $\bm{p}$ & generic parameters \\
  $P$ & electric power [W]\\
  $\dot{Q}$ & heat transfer rate [W]\\
  $s$ & specific entropy [W/kg/K] \\
  $t$ & time point \\
  $T$ & temperature [K] \\
  $U$ & heat transfer coefficient [W/m$^2$/K] \\
  $V, \dot{V}$ & volume, volumetric flow [m$^3$], [m$^3$/s] \\
  $\bm{x}$ & generic design variables \\
  $\bm{y}$ & generic operational variables \\
  $\mathcal{T}$ & set of all considered time points \\
  $\mathcal{X}$ & host-set of generic design variables \\
  $\mathcal{Y}$ & host-set of generic operational variables \\
  $w$ & scenario weight
\end{tabularx} \\[5mm]
%
\noindent\textbf{Greek symbols} \\
\noindent
\begin{tabularx}{\columnwidth}{lX}
  $\Delta_{s,t}$ & time step [h]\\
  $\Delta$T & temperature difference [K]\\
  $\eta$ & efficiency \\
  $\rho$ & density [kg/m$^3$]\\
  $\tau$ & self-discharge of storage component [h]
\end{tabularx} \\[5mm]
\newpage
\noindent\textbf{Subscripts} \\
\noindent
\begin{tabularx}{\columnwidth}{lX}
  0 & initial point \\
  \makecell[l]{1, 2, 2r, \\
               3, 4, 5, 6, \\
               6r, pinch} &
  working fluid states in the fourth case study \\
  \makecell[l]{40, 50, \\
               70, 85} & design temperatures in the third case study \\
  A,B & thermal zones A and B in the second case study \\
  c & cool \\
  con & condenser \\
  core & concrete core \\
  cw & cooling water \\
  eva & evaporator \\
  eco & economizer \\
  gb & geothermal brine \\
  h & hot \\
  $i$ & generic system component  \\
  ib & isobutane \\
  I, II & first- and second-stage quantities \\
  rec & recuperator \\
  $s$ & scenario \\
  sup & superheater \\
\end{tabularx} \\[5mm]
%
\noindent\textbf{Superscripts} \\
\noindent
\begin{tabularx}{\columnwidth}{lX}
  conv & conversion components \\
  d & differential states \\
  elec & electricity \\
  fl & flow \\
  gr & ground \\
  I & investment \\
  in & input, in-flowing stream \\
  is & isentropic \\
  liq & liquid \\
  max & maximum value \\
  min & minimum value \\
  nom & nominal value \\
  out & output, out-flowing stream \\
  re & return \\
  ref & reference value \\
  sat & saturation \\
  sto & storage components \\
  vap & vapor \\
\end{tabularx}